\documentclass[11pt,twoside]{amsart}
\usepackage{graphicx}
\usepackage{amssymb}
\usepackage{amsmath}
\usepackage{amsthm}
\usepackage{amsfonts}
\usepackage{amsmath, amsfonts, amssymb}
\usepackage{multirow}
\usepackage{rotating}
\usepackage{multirow}
\usepackage{rotating}
\usepackage{multirow}
\usepackage{rotating}
\usepackage{longtable}
\usepackage{color}
\newtheorem{theorem}{Theorem}[section]
\newtheorem{corollary}[theorem]{Corollary}

\newtheorem{proposition}[theorem]{Proposition}
\newtheorem{definition}[theorem]{Definition}
\newtheorem{remark}[theorem]{\bf{Remark}}

\numberwithin{equation}{section}

\begin{document}
	\small \begin{center}{\textit{ In the name of
				Allah, the Beneficent, the Merciful}}\end{center}
	\vspace{0.5cm}
	\title{Subalgebras, idempotents, ideals and quasi-units of two-dimensional algebras}
	
	\author{H.Ahmed$^1$, U.Bekbaev$^2$, I.Rakhimov$^3$}
	
	\thanks{{\scriptsize
			emails: $^1$houida\_m7@yahoo.com; $^2$uralbekbaev@gmail.com;
			$^3$rakhimov@upm.edu.my.}}
	\maketitle
	\begin{center}
		{\scriptsize \address{$^{1}$Depart. of Math., Faculty of Science, Taiz University,
				Taiz, Yemen $\&$ \\ Department of Math., Faculty of Science, UPM,
				Selangor, Malaysia}\\
			\address{$^2$Department of Natural and Mathematical Sciences, TTPU, Tashkent, Uzbekistan}\\
			\address{$^3$Laboratory of Cryptography, Analysis and Structure,\\
				Institute for Mathematical Research (INSPEM), UPM, Serdang, Selangor, Malaysia}}
	\end{center}

\begin{abstract}
All subalgebras, idempotents, left(right) ideals and left quasi-units of two-dimensional algebras are described. Classification of algebras with given number of subalgebras, left(right) ideals are provided. In particular, a list of isomorphism classes of all simple two-dimensional algebras is given.
In the study of ideals and subalgebras the number of them depend on roots of certain system of polynomials at structure constants of the algebra. We also give explicit forms of the polynomials.
\end{abstract}

keywords: {Subalgebra; idempotent; ideal; quasi-unit.}\\

ccode: {Mathematics Subject Classification 2010: 15A72; 16D25; 17A30; 17A99. }

\section{Introduction}
	It is noticed that well known classification theorems of algebras, for example Lie, associative, Jordan and etc, in fixed $n$-dimensional case cover only small parts of all $n$-dimensional algebras. In fact, the part out of the consideration is much bigger and it is a dense in Zariski topology subset of the set of all $n$-dimensional algebras. Therefore, instead of the  classification of some classes of algebras one can try to classify all algebras in fixed dimensions. Such an attempt has been made in \cite{P2000} for all two-dimensional algebras over an infinite field $\mathbb{F}.$ Here the author used basis (coordinate) free, invariant approach. Unlike to the case of \cite{P2000} in \cite{A} we have given a coordinate based classification of such algebras over any field $\mathbb{F},$ where  quadratic and cubic polynomials have a root. The case $\mathbb{F}=\mathbb{R}$ has been considered in \cite{B}. The paper \cite{K} contains another coordinate based approach to the classification problem of two-dimensional algebras. The main advantage of the approach, proposed in \cite{B0} and applied in \cite{A}, is the fact that it reduces the classification and some other problems of finite-dimensional algebras to the investigation of a system of equations. Particularly, in two-dimensional case the problems are written as systems of equations for the twelve classes of canonical representatives of all nontrivial two-dimensional algebras found in \cite{A}. By the use of the results of \cite{A} one can classify, up to isomorphism, nearly any class of two-dimensional algebras. Particularly, the approach has been used in \cite{A1, A2} to get the classification of some classes of two-dimensional algebras which were unknown before. In \cite{A3} it was used for the description of automorphism groups and derivation algebras of all two-dimensional algebras.
	
	 Note that in two-dimensional case the number of classes of canonical representatives is small. However, increasing the dimension of the algebra produces a large number of such isomorphism classes, calculations show that even in three dimensional algebras case the corresponding number is expected to be about one thousand \cite{B1}.
	
	The present work grew out of an attempt to respond a comment made by a reviewer of our paper. The reviewer drown our attention to the importance of the description of all subalgebras, idempotents, left(right) ideals and
	quasi-units of two-dimensional algebras through those twelve classes of canonical representatives in \cite{A}. Recall that we consider only nontrivial two-dimensional algebras and by a subalgebra (left(right) ideal) we mean nontrivial, that is one dimensional, subalgebra (respectively, left(right) ideal). It is shown that every such an algebra may have only one, two, three or infinitely many subalgebras and may have zero, one, two or infinitely many ideals. We provide a classification for each class of algebras with $n$ different subalgebras (left(right) ideals), where $n=1,2,3,\infty$ (respectively,  $n=0,1,2,\infty$) and subsequently, give a classification of two-dimensional simple algebras. Another result of the paper is the descriptions of idempotents and left quasi-units of two-dimensional algebras.
	
	Here we cite only some of the papers related to the classification problem of all two-dimensional algebras \cite{Al,Al1,Berm,Berm1,Bu,C,D,G,G1,Wallace}. There are many others dealing with some particular cases of the problem.
	
	The paper is organized as follows. Section 2 contains preliminary results, Section 3 deals with subalgebras  and in Section 4 we list idempotents of two-dimensional algebras. In Sections 5, 6 and 7 we study the left, right and two-sided ideals of two-dimensional algebras, respectively. The left quasiunits are considered in Section 8. The results of the paper we tabulate in Section 9.
	
	 Throughout the paper $\mathbb{F}$ stands for any field over that every quadratic and cubic polynomial has a root in $\mathbb{F}$, in particular, $\mathbb{F}$ may be any algebraically closed field.
	
	\section{Preliminaries} \label{Sec1}
	
	Further we need the following simple result on roots of a polynomial $p(y)=ay^3+by^2+cy+d$ over a field $\mathbb{F}$.
	\begin{proposition}\emph{} \label{Pr21}

\begin{itemize}
		\item[$\bullet$]\emph{}Let $Char(\mathbb{F})\neq 2,3$.

\begin{itemize}
\item[$\star \  \ \mbox{If}\ a\neq 0$] the polynomial $p(y)$ has only
\begin{itemize}
\item one root if and only if $b^2-3ac=0$ and $bc-9ad=0;$
\item  it has only two different roots  if and only if $b^2-3ac\neq 0$ and $$p\left(\frac{-b+\sqrt{b^2-3ac}}{3a}\right)p\left(\frac{-b-\sqrt{b^2-3ac}}{3a}\right)=0;$$
\item it has three different roots if and only if  $$p\left(\frac{-b+\sqrt{b^2-3ac}}{3a}\right)p\left(\frac{-b-\sqrt{b^2-3ac}}{3a}\right)\neq 0.$$
\end{itemize}
\item[$\star \  \ \mbox{If}\  a= 0$] the polynomial $p(y)$ has
\begin{itemize}
\item no root if and only if $b=c=0$ and $d\neq 0;$
\item it has only one root if and only if $b\neq 0$ and $c^2-4bd=0$ or $b=0$ and $c\neq 0;$
\item has only two different roots if and only if $b\neq 0$ and $c^2-4bd\neq 0;$
\item it has infinitely many different roots  if and only if $b=c=d=0.$
\end{itemize}
\end{itemize}		
		\item[$\bullet$] \emph{}Let $Char(\mathbb{F})=2$.

\begin{itemize}
\item[$\star \  \ \mbox{If}\ a\neq 0$] the polynomial $p(y)$ has
\begin{itemize}
\item one root if and only if $ad-bc=0$ and $ac-b^2=0;$
\item  it has only two different roots  if and only if $ad-bc=0$ and $ac-b^2\neq 0;$
\item it has three different roots if and only if $ad-bc\neq 0.$
\end{itemize}
\item[$\star \  \ \mbox{If}\  a= 0$] the polynomial $p(y)$
\begin{itemize}
\item has no root if and only if $b=c=0$ and $d\neq 0$;
\item it has only one root if and only if $b\neq 0$ and $c=0$ or $b= 0$ and $c\neq 0$;
\item  it has only two different roots if and only if $b\neq 0$ and $c\neq 0$;
\item  it has infinitely many different roots  if and only if $b=c=d=0$.
\end{itemize}
\end{itemize}
\item[$\bullet$] \emph{}Let $Char(\mathbb{F})=3.$

\begin{itemize}
\item[$\star \  \ \mbox{If}\ a\neq 0$] the polynomial $p(y)$ has
\begin{itemize}
\item one root if and only if $b=c=0;$
\item  it has only two different roots  if and only if $b\neq 0$ and $p\left(\frac{c}{b}\right)=0;$
\item it has three different roots if and only if $b\neq 0$ and $p\left(\frac{c}{b}\right)\neq 0$ or $b=0$ and $c\neq 0.$
\end{itemize}
\item[$\star \  \ \mbox{If}\  a= 0$] the polynomial $p(y)$
\begin{itemize}
\item has no root if and only if $b=c=0$ and $d\neq 0$;
\item it has only one root if and only if $b\neq 0$ and $c^2-bd= 0$ or $b=0$ and $c\neq 0$;
\item  it has only two different roots if and only if $b\neq 0$ and $c^2-bd\neq 0$;
\item  it has infinitely many different roots  if and only if $b=c=d=0$.
\end{itemize}
\end{itemize}
\end{itemize}
	\end{proposition}
	
	\begin{proof} The proofs can be derived from the following facts.
\begin{itemize}
\item If $Char(\mathbb{F})\neq 2,3$ under $a\neq 0$ we have  $$p(y)=\frac{3ay+b}{9a}p'(y)+\frac{6ac-2b^2}{9a}y+\frac{9ad-bc}{9a},\ \ p'(y)=\frac{y}{2}p''(y)+by+c\ \mbox{and}\ p''(y)=6ay+2b;$$
\item $Char(\mathbb{F})=2.$ If $a\neq 0$ then one has $$p(y)=\left(y+\frac{b}{a}\right)p'(y)+\frac{ad-bc}{a}, \ \ p'(y)=ay^2+c \ \mbox{and}\ p''(y)=0;$$

\item $Char(\mathbb{F})=3.$ If $b\neq 0$ then $$p(y)=-\left(\frac{a}{b}y^2+\frac{b^2+ac}{b^2}y+\frac{c(2b^2+ac)}{b^3}\right)p'(y)+p\left(\frac{c}{b}\right), \ \ p'(y)=2by+c\ \mbox{and}\ p''(y)=2b.$$
\end{itemize}

The case $a=0$ is immediate.
\end{proof}	

	Let $\mathbb{A}$ be a two-dimensional algebra over $\mathbb{F}$ and $\{e_1,e_2\}$ be its basis and  \[A=\left(
	\begin{array}{cccc}
	\alpha _1 & \alpha _2 & \alpha _3 & \alpha _4 \\
	\beta _1 & \beta _2 & \beta _3 & \beta _4
	\end{array}
	\right)\] be the matrix of structure constants (MSC) of $\mathbb{A}$ in $\{e_1,e_2\}$ arranged as $$e_1e_1=\alpha_1e_1+\beta_1e_2,\ \  e_1e_2=\alpha_2e_1+\beta_2e_2,\ \
	e_2e_1=\alpha_3e_1+\beta_3e_2, \ \ e_2e_2=\alpha_4e_1+\beta_4e_2.$$
	
	Further it is assumed that a basis $\mathrm{e}=\{e_1,e_2\}$ is fixed and we do not distinguish between $\mathbb{A}$ and its MSC $A$.

	In the paper we make use the following results and notations from \cite{A}.
	
	Any nontrivial two-dimensional algebra over $\mathbb{F}$ is isomorphic to only one of the following algebras given by their MSC:
	
	In   $Char(\mathbb{F})\neq 2,3$ case:
	\begin{itemize}
		\item $A_{1}(\mathbf{c})=\left(
		\begin{array}{cccc}
		\alpha_1 & \alpha_2 &\alpha_2+1 & \alpha_4 \\
		\beta_1 & -\alpha_1 & -\alpha_1+1 & -\alpha_2
		\end{array}\right),\ \mbox{where}\ \mathbf{c}=(\alpha_1, \alpha_2,
		\alpha_4, \beta_1)\in \mathbb{F}^4,$
		\item $A_{2}(\mathbf{c})=\left(
		\begin{array}{cccc}
		\alpha_1 & 0 & 0 & 1 \\
		\beta _1& \beta _2& 1-\alpha_1&0
		\end{array}\right)\simeq \left(
		\begin{array}{cccc}
		\alpha_1 & 0 & 0 & 1 \\
		-\beta _1& \beta _2& 1-\alpha_1&0
		\end{array}\right),\ \mbox{where}\ \mathbf{c}=(\alpha_1, \beta_1,
		\beta_2)\in \mathbb{F}^3,$
		\item $A_{3}(\mathbf{c})=\left(
		\begin{array}{cccc}
		0 & 1 & 1 & 0 \\
		\beta _1& \beta _2 & 1&-1
		\end{array}\right),\ \mbox{where}\ \mathbf{c}=(\beta_1, \beta_2)\in
		\mathbb{F}^2,$
		\item $A_{4}(\mathbf{c})=\left(
		\begin{array}{cccc}
		\alpha _1 & 0 & 0 & 0 \\
		0 & \beta _2& 1-\alpha _1&0
		\end{array}\right),\ \mbox{where}\ \mathbf{c}=(\alpha_1, \beta_2)\in
		\mathbb{F}^2,$
		\item $A_{5}(\mathbf{c})=\left(
		\begin{array}{cccc}
		\alpha_1& 0 & 0 & 0 \\
		1 & 2\alpha_1-1 & 1-\alpha_1&0
		\end{array}\right),\ \mbox{where}\ \mathbf{c}=\alpha_1\in \mathbb{F},$
		\item $A_{6}(\mathbf{c})=\left(
		\begin{array}{cccc}
		\alpha_1 & 0 & 0 & 1 \\
		\beta _1& 1-\alpha_1 & -\alpha_1&0
		\end{array}\right)\simeq \left(
		\begin{array}{cccc}
		\alpha_1 & 0 & 0 & 1 \\
		-\beta _1& 1-\alpha_1 & -\alpha_1&0
		\end{array}\right),\ \mbox{where}\ \mathbf{c}=(\alpha_1, \beta_1)\in
		\mathbb{F}^2,$
		\item $A_{7}(\mathbf{c})=\left(
		\begin{array}{cccc}
		0 & 1 & 1 & 0 \\
		\beta_1& 1& 0&-1
		\end{array}\right),\ \mbox{where}\ \mathbf{c}=\beta_1\in \mathbb{F},$
		\item $A_{8}(\mathbf{c})=\left(
		\begin{array}{cccc}
		\alpha_1 & 0 & 0 & 0 \\
		0 & 1-\alpha_1 & -\alpha_1&0
		\end{array}\right),\ \mbox{where}\ \mathbf{c}=\alpha_1\in \mathbb{F},$
		\item $A_{9}=\left(
		\begin{array}{cccc}
		\frac{1}{3}& 0 & 0 & 0 \\
		1 & \frac{2}{3} & -\frac{1}{3}&0
		\end{array}\right),\ \ A_{10}=\left(
		\begin{array}{cccc}
		0 & 1 & 1 & 0 \\
		0 &0& 0 &-1
		\end{array}
		\right),\ \ A_{11}=\left(
		\begin{array}{cccc}
		0 & 1 & 1 & 0 \\
		1 &0& 0 &-1
		\end{array}
		\right),\\ A_{12}=\left(
		\begin{array}{cccc}
		0 & 0 & 0 & 0 \\
		1 &0&0 &0\end{array}
		\right).$
	\end{itemize}
	
	In   $Char(\mathbb{F})= 2$ case:
	\begin{itemize}
		\item $A_{1,2}(\mathbf{c})=\left(
		\begin{array}{cccc}
		\alpha_1 & \alpha_2 &\alpha_2+1 & \alpha_4 \\
		\beta_1 & -\alpha_1 & -\alpha_1+1 & -\alpha_2
		\end{array}\right),\ \mbox{where}\ \mathbf{c}=(\alpha_1, \alpha_2,
		\alpha_4, \beta_1)\in \mathbb{F}^4,$
		\item $A_{2,2}(\mathbf{c})=\left(
		\begin{array}{cccc}
		\alpha_1 & 0 & 0 & 1 \\
		\beta _1& \beta_2 & 1-\alpha_1&0
		\end{array}\right),\ \mbox{where}\ \mathbf{c}=(\alpha_1, \beta_1,
		\beta_2)\in \mathbb{F}^3,$
		\item $A_{3,2}(\mathbf{c})=\left(
		\begin{array}{cccc}
		\alpha_1 & 1 & 1 & 0 \\
		0& \beta_2 & 1-\alpha_1&1
		\end{array}\right),\ \mbox{where}\ \mathbf{c}=(\alpha_1, \beta_2)\in
		\mathbb{F}^2,$
		\item $A_{4,2}(\mathbf{c})=\left(
		\begin{array}{cccc}
		\alpha _1 & 0 & 0 & 0 \\
		0 & \beta_2 & 1-\alpha _1&0
		\end{array}\right),\ \mbox{where}\ \mathbf{c}=(\alpha_1,\beta_2)\in
		\mathbb{F}^2,$
		\item $A_{5,2}(\mathbf{c})=\left(
		\begin{array}{cccc}
		\alpha_1 & 0 & 0 & 0 \\
		1 & 1 & 1-\alpha_1&0
		\end{array}\right),\ \mbox{where}\ \mathbf{c}=\alpha_1\in \mathbb{F},$
		\item $A_{6,2}(\mathbf{c})=\left(
		\begin{array}{cccc}
		\alpha_1 & 0 & 0 & 1 \\
		\beta _1& 1-\alpha_1 & -\alpha_1&0
		\end{array}\right),\ \mbox{where}\ \mathbf{c}=(\alpha_1, \beta_1)\in
		\mathbb{F}^2,$
		\item $A_{7,2}(\mathbf{c})=\left(
		\begin{array}{cccc}
		\alpha_1 & 1 & 1 & 0 \\
		0& 1-\alpha_1& -\alpha_1&-1
		\end{array}\right),\ \mbox{where}\ \mathbf{c}=\alpha_1\in \mathbb{F},$
		\item $A_{8,2}(\mathbf{c})=\left(
		\begin{array}{cccc}
		\alpha_1 & 0 & 0 & 0 \\
		0 & 1-\alpha_1 & -\alpha_1&0
		\end{array}\right),\ \mbox{where}\ \mathbf{c}=\alpha_1\in \mathbb{F},$
		\item $A_{9,2}=\left(
		\begin{array}{cccc}
		1 & 0 & 0 & 0 \\
		1 & 0 & 1&0
		\end{array}\right),\ \ A_{10,2}=\left(
		\begin{array}{cccc}
		0 & 1 & 1 & 0 \\
		0 &0& 0 &-1
		\end{array}
		\right),\ \ A_{11,2}=\left(
		\begin{array}{cccc}
		1 & 1 & 1 & 0 \\
		0 &-1& -1 &-1
		\end{array}
		\right),\\ A_{12,2}=\left(
		\begin{array}{cccc}
		0 & 0 & 0 & 0 \\
		1 &0&0 &0\end{array}
		\right).$
	\end{itemize}
		
	In   $Char(\mathbb{F})=3$ case:
	
	\begin{itemize}
		\item $A_{1,3}(\mathbf{c})=\left(
		\begin{array}{cccc}
		\alpha_1 & \alpha_2 &\alpha_2+1 & \alpha_4 \\
		\beta_1 & -\alpha_1 & -\alpha_1+1 & -\alpha_2
		\end{array}\right),\ \mbox{where}\ \mathbf{c}=(\alpha_1, \alpha_2,
		\alpha_4, \beta_1)\in \mathbb{F}^4,$
		\item $A_{2,3}(\mathbf{c})=\left(
		\begin{array}{cccc}
		\alpha_1 & 0 & 0 & 1 \\
		\beta _1& \beta _2& 1-\alpha_1&0
		\end{array}\right)\simeq\left(
		\begin{array}{cccc}
		\alpha_1 & 0 & 0 & 1 \\
		-\beta _1& \beta _2& 1-\alpha_1&0
		\end{array}\right),\ \mbox{where}\ \mathbf{c}=(\alpha_1, \beta_1,
		\beta_2)\in \mathbb{F}^3,$
		\item $A_{3,3}(\mathbf{c})=\left(
		\begin{array}{cccc}
		0 & 1 & 1 & 0 \\
		\beta _1& \beta _2 & 1&-1
		\end{array}\right),\ \mbox{where}\ \mathbf{c}=(\beta_1, \beta_2)\in
		\mathbb{F}^2,$
		\item $A_{4,3}(\mathbf{c})=\left(
		\begin{array}{cccc}
		\alpha _1 & 0 & 0 & 0 \\
		0 & \beta _2& 1-\alpha _1&0
		\end{array}\right),\ \mbox{where}\ \mathbf{c}=(\alpha_1, \beta_2)\in
		\mathbb{F}^2,$
		\item $A_{5,3}(\mathbf{c})=\left(
		\begin{array}{cccc}
		\alpha_1& 0 & 0 & 0 \\
		1 & -1-\alpha_1 & 1-\alpha_1&0
		\end{array}\right),\ \mbox{where}\ \mathbf{c}=\alpha_1\in \mathbb{F},$
		\item $A_{6,3}(\mathbf{c})=\left(
		\begin{array}{cccc}
		\alpha_1 & 0 & 0 & 1 \\
		\beta _1& 1-\alpha_1 & -\alpha_1&0
		\end{array}\right)\simeq \left(
		\begin{array}{cccc}
		\alpha_1 & 0 & 0 & 1 \\
		-\beta _1& 1-\alpha_1 & -\alpha_1&0
		\end{array}\right),\ \mbox{where}\ \mathbf{c}=(\alpha_1, \beta_1)\in
		\mathbb{F}^2,$
		\item $A_{7,3}(\mathbf{c})=\left(
		\begin{array}{cccc}
		0 & 1 & 1 & 0 \\
		\beta_1& 1& 0&-1
		\end{array}\right),\ \mbox{where}\ \mathbf{c}=\beta_1\in \mathbb{F},$
		\item $A_{8,3}(\mathbf{c})=\left(
		\begin{array}{cccc}
		\alpha_1 & 0 & 0 & 0 \\
		0 & 1-\alpha_1 & -\alpha_1&0
		\end{array}\right),\ \mbox{where}\ \mathbf{c}=\alpha_1\in \mathbb{F},$
		\item $A_{9,3}=\left(
		\begin{array}{cccc}
		0 & 1& 1& 0 \\
		1 &0&0 &-1\end{array}
		\right),\ \ A_{10,3}=\left(
		\begin{array}{cccc}
		0 & 1 & 1 & 0 \\
		0 &0&0 &-1\end{array}
		\right),$
		\item $A_{11,3}=\left(
		\begin{array}{cccc}
		1 & 0 & 0 & 0 \\
		1 &-1&-1 &0\end{array}
		\right),\ \ A_{12,3}=\left(
		\begin{array}{cccc}
		0 &0 &0 & 0 \\
		1 &0& 0 &0
		\end{array}\right).$
	\end{itemize}
	
	\section{Subalgebras} \label{Sec3}
Any nontrivial subalgebra of a two-dimensional algebra $\mathbb{A}$ is one dimensional, it is generated by a nonzero element $\mathbf{u}=u_1e_1+u_2e_2\in \mathbb{A}$ such that  $\mathbf{u}^2=\lambda \mathbf{u}$ for some $\lambda\in \mathbb{F}$. If $u_1=0$ then, due to $\mathbb{F}\mathbf{u}=\mathbb{F}k\mathbf{u}$ for any $k\neq 0$, $\mathbf{u}$ can be represented as $\mathbf{u}=e_2$, if $u_1\neq 0$ then $\mathbf{u}$ can be represented as $\mathbf{u}=e_1+y_0e_2$.
	
	Let polynomial $p_A(y)$ stand for  \[p_A(y)=y^3 \alpha _4+y^2(\alpha _2+\alpha _3-\beta _4) +y(\alpha _1-\beta _2-\beta _3)-\beta _1.\]
	\begin{proposition} \label{P22}
The set of all nontrivial subalgebras of $\mathbb{A}$ is given as follows
\begin{itemize}		
\item $\{\mathbb{F}(e_1+y_0e_2): y_0 \in \mathbb{F},\ p_A(y_0)=0\} \ \mbox{if}\ \alpha_4\neq 0;$
\item $\{\mathbb{F}(e_1+y_0e_2): y_0\in \mathbb{F},\ p_A(y_0)=0\}\cup \{\mathbb{F}e_2\}\ \mbox{if}\ \alpha_4= 0.$
\end{itemize}
	\end{proposition}
	\begin{proof} Due to $\mathbf{u}\cdot \mathbf{v}=\mathrm{e}A(u\otimes v),$ where  $\mathbf{u}=\mathrm{e}u,\mathbf{v}=\mathrm{e}v$ in terms of $A$ and $u$ the equality $\mathbf{u}^2=\lambda \mathbf{u}$ is written
	\begin{equation}\label{2}
	A(u\otimes u)-\lambda u=0.
	\end{equation}
	Assuming $\mathbf{u}=xe_1+ye_2=\mathrm{e}\left(
	\begin{array}{c}x\\ y\end{array}
	\right)$ and applying (\ref{2}) we obtain the system of equations:
	\begin{equation}\label{SEs}
	\begin{array}{rr}
	-x \lambda +x^2 \alpha _1+x y \alpha _2+x y \alpha _3+y^2 \alpha _4& =0, \\
	-y \lambda +x^2 \beta _1+x y \beta _2+x y \beta _3+y^2 \beta _4&=0.
	\end{array}
	\end{equation}
	
	If $x= 0$ then one can put $y=1$ and therefore in $\alpha_4\neq 0$ case it has no solution, whereas if $\alpha_4= 0$ then $u=\left(
	\begin{array}{c}
	0 \\
	1
	\end{array}
	\right)$ is a solution and therefore there is at least one nontrivial subalgebra $\mathbb{F}e_2$.
	
	If $x\neq 0$ one can assume $x=1$ and the system (\ref{SEs}) becomes equivalent to \[p_A(y)=y^3 \alpha _4+y^2(\alpha _2+\alpha _3-\beta _4) +y(\alpha _1-\beta _2-\beta _3)-\beta _1=0.\]
	Therefore, each root $y_0$ of $p_A(y)$ provides the corresponding subalgebra $\mathbb{F}(e_1+y_0e_2)$.
	\end{proof}
	So  the following result holds true.
	\begin{corollary} Any two-dimensional algebra may have only one, two, three or infinitely many nontrivial subalgebras. In the case of infinitely many subalgebras every $1$-dimensional subspace is a subalgebra.\end{corollary}
	
	Note that Proposition \ref{P22} also is true for any two-dimensional algebra over a field $\mathbb{F}$.

	In this section we describe all nontrivial subalgebras of the classes of canonical representatives of two-dimensional algebras given in Section \ref{Sec1}.
	\begin{theorem} Let  $Char(\mathbb{F})\neq 2,3.$ Two-dimensional algebras over $\mathbb{F}$ with respect to the number of non-trivial subalgebras are distributed as follows.
\begin{itemize}
\item Any two-dimensional algebra with only one subalgebra is isomorphic to one of the following algebras
\begin{itemize}
		\item[$\ast$] $A_{1}(\alpha_1,\alpha_2,\alpha_4,\beta_1),$ where $\alpha_4\neq 0,$ $\alpha_1=\frac{(3\alpha_2+1)^2}{9\alpha_4}+\frac{1}{3},$ $ \beta_1=\frac{-(3\alpha_2+1)^3}{27\alpha_4^2};$
\item[$\ast$] $A_{1}(\frac{1}{3},-\frac{1}{3},0,\beta_1)$, where $\beta_1\neq
		0;$
\item[$\ast$] $A_{2}(\alpha_1,0, 2\alpha_1-1);$
\item[$\ast$] $A_{5}(\alpha_1)$;
\item[$\ast$] $A_{6}(\frac{1}{3},0)$;
\item[$\ast$] $A_{9}$;
\item[$\ast$] $A_{12}.$
\end{itemize}
\item Any two-dimensional algebra with two subalgebras is isomorphic to one of the following algebras
\begin{itemize}
		\item[$\ast$] $A_{1}(\alpha_1,\alpha_2,\alpha_4,\beta_1)$, where $\alpha_4\neq 0,$ $\alpha_1\neq\frac{(3\alpha_2+1)^2}{9\alpha_4}+\frac{1}{3},$ $p_{A_1}(y_{-1})p_{A_1}(y_{+1})=0$, and \[y_{-1}=\frac{-(3\alpha_2+1)-\sqrt{(3\alpha_2+1)^2-\alpha_4(3\alpha_1-1)}}{3\alpha_4},\] \[y_{+1}=\frac{-(3\alpha_2+1)+\sqrt{(3\alpha_2+1)^2-\alpha_4(3\alpha_1-1)}}{3\alpha_4};\]
		\item[$\ast$] $A_{1}(\alpha_1,\alpha_2,0,\beta_1)$, where
		$\alpha_2\neq -\frac{1}{3},
		(3\alpha_1-1)^2+4(\alpha_2+1)\beta_1= 0$;
		\item[$\ast$] $A_{1}(\alpha_1,-\frac{1}{3},0,\beta_1)$, where
		$\alpha_1\neq\frac{1}{3}$;
		\item[$\ast$] $A_{2}(\alpha_1,\beta_1, \beta_2)$, where $\beta_2\neq 2\alpha_1-1$, $p_{A_2}(y_{-2})p_{A_2}(y_{+2})=0$,
		\[y_{-2}=\frac{-\sqrt{3(1+\beta_2-2\alpha_1)}}{3},\ y_{+2}= \frac{\sqrt{3(1+\beta_2-2\alpha_1)}}{3};\]
\item[$\ast$] $A_{3}(-\frac{(\beta_2+1)^2}{12}, \beta_2)$;
\item[$\ast$] $A_{4}(\alpha_1, \beta_2)$, where $\beta_2\neq 2\alpha_1-1$;
\item[$\ast$] $A_{6}(\alpha_1,\beta_1)$, where $\alpha_1\neq \frac{1}{3}$, $ p_{A_6}(y_{-6})p_{A_6}(y_{+6})=0$,
		\[y_{-6}=\frac{\sqrt{3(1-3\alpha_1)}}{3},\ y_{+6}=\frac{-\sqrt{3(1-3\alpha_1)}}{3};\]
\item[$\ast$] $A_{7}(-\frac{1}{12})$;
\item[$\ast$] $A_{8}(\alpha_1)$, where $\alpha_1\neq \frac{1}{3}$;
\item[$\ast$] $A_{10}$.
\end{itemize}
\item Any two-dimensional algebra with three subalgebras is isomorphic to one of the following algebras
\begin{itemize}
		\item[$\ast$] $A_{1}(\alpha_1,\alpha_2,\alpha_4,\beta_1)$, where $\alpha_4p_{A_1}(y_{-1})p_{A_1}(y_{+1})\neq 0$;
\item[$\ast$] $A_{1}(\alpha_1,\alpha_2,0,\beta_1)$, where $(3\alpha_2+1)((3\alpha_1-1)^2+4(\alpha_2+1)\beta_1)\neq 0$;
	\item[$\ast$] $A_{2}(\alpha_1,\beta_1, \beta_2)$, where $p_{A_2}(y_{-2})p_{A_2}(y_{+2})\neq 0$;
		\item[$\ast$] $A_{3}(\beta_1, \beta_2)$, where $(\beta_2+1)^2+12\beta_1\neq 0$;
		\item[$\ast$] $A_{6}(\alpha_1,\beta_1)$, where $ p_{A_6}(y_{-6})p_{A_6}(y_{+6})\neq 0;$
		\item[$\ast$] $A_{7}(\beta_1)$, where $12\beta_1+1\neq 0$;
\item[$\ast$] $A_{11}$.
\end{itemize}
\item Any two-dimensional algebra with infinitely many subalgebras is isomorphic to one of the following algebras
\begin{itemize}
		\item[$\ast$] $A_{1}(\frac{1}{3},-\frac{1}{3},0,0)$;
\item[$\ast$] $A_{4}(\alpha_1,2\alpha_1-1)$;
\item[$\ast$] $A_{8}(\frac{1}{3})$.
	\end{itemize}
\end{itemize}
	\end{theorem}
	
	\begin{proof}
	In $A=A_{1}(\alpha_1,\alpha_2,\alpha_4,\beta_1)=\left(
	\begin{array}{cccc}
	\alpha _1 & \alpha _2 & \alpha _2+1 & \alpha _4 \\
	\beta _1 & -\alpha _1 & -\alpha _1+1 & -\alpha _2
	\end{array}
	\right)$ case one has \[p_A(y)=p_{A_1}(y)=\alpha_4 y^3+(3\alpha_2+1)y^2+(3\alpha_1-1)y-\beta_1.\] Therefore, due to  Propositions \ref{Pr21} and \ref{P22} to know the number of nontrivial subalgebras one should study the number of roots of $p_A(y).$
	In $\alpha_4\neq 0$ case the algebra $A_{1}(\alpha_1,\alpha_2,\alpha_4,\beta_1)$ has only \textit{one nontrivial subalgebra} if and only if $$3\alpha_4(3\alpha_1-1)-(3\alpha_2+1)^2=0\ \mbox{and}\  (3\alpha_2+1)(3\alpha_1-1)+9\beta_1\alpha_4=0,$$ i.e., $\alpha_1=\frac{(3\alpha_2+1)^2}{9\alpha_4}+\frac{1}{3},
	\beta_1=\frac{-(3\alpha_2+1)^3}{27\alpha_4^2},$
	the algebra $A_{1}(\alpha_1,\alpha_2,\alpha_4,\beta_1)$ has only \textit{two nontrivial subalgebras} if and only if \[3\alpha_4(3\alpha_1-1)-(3\alpha_2+1)^2\neq 0\] and one of  \[p_{A_1}\left(\frac{-(3\alpha_2+1)+\sqrt{(3\alpha_2+1)^2-\alpha_4(3\alpha_1-1)}}{3\alpha_4}\right), \ \ p_{A_1}\left(\frac{-(3\alpha_2+1)-\sqrt{(3\alpha_2+1)^2-\alpha_4(3\alpha_1-1)}}{3\alpha_4}\right)\] is zero, and finally,
	the algebra $A_{1}(\alpha_1,\alpha_2,\alpha_4,\beta_1)$ has only \textit{three nontrivial subalgebras} if and only if none of \[p_{A_1}\left(\frac{-(3\alpha_2+1)+\sqrt{(3\alpha_2+1)^2-\alpha_4(3\alpha_1-1)}}{3\alpha_4}\right), \ \ p_{A_1}\left(\frac{-(3\alpha_2+1)-\sqrt{(3\alpha_2+1)^2-\alpha_4(3\alpha_1-1)}}{3\alpha_4}\right)\] is zero.
	
	Similarly, in $\alpha_4= 0$ case the algebra $A_{1}(\alpha_1,\alpha_2,0,\beta_1)$ has only
	\emph{one nontrivial subalgebra} if and only if \[3\alpha_1-1=3\alpha_2+1=0\mbox{ and } \beta_1\neq 0,\]
	it has \emph{two nontrivial subalgebras} if and only if \[3\alpha_2+1\neq 0 \ , (3\alpha_1-1)^2+4(\alpha_2+1)\beta_1=0\] or \[3\alpha_2+1= 0 \mbox{ and } 3\alpha_1-1\neq 0,\]
	it has \emph{three nontrivial subalgebras} if and only if \[3\alpha_2+1\neq 0 \ , (3\alpha_1-1)^2+4(\alpha_2+1)\beta_1\neq 0\]
	and it has \emph{infinitely many nontrivial $1$-dimensional subalgebras} if and only if \[3\alpha_2+1=3\alpha_1-1=\beta_1=0.\]

	In $A_{2}(\alpha_1,\beta_1,\beta_2)=\left(
	\begin{array}{cccc}
	\alpha _1 & 0 & 0 & 1 \\
	\beta _1 & \beta _2 & -\alpha _1+1 & 0
	\end{array}
	\right)$ case $\alpha_4=1\neq 0$ and
	\[ p_{A_2}(y)=y^3+(2\alpha_1-\beta_2-1)y-\beta_1.\]
	
	Therefore, the algebra $A_{2}(\alpha_1,\beta_1, \beta_2)$ has only
	one nontrivial subalgebra if and only if \[2\alpha_1-\beta_2-1=0\mbox{ and } \beta_1= 0,\]
	it has two nontrivial subalgebras if and only if \[2\alpha_1-\beta_2-1\neq 0 \mbox{ and one of } p_{A_2}\left(\frac{\sqrt{3(1+\beta_2-2\alpha_1)}}{3}\right), \ \ p_{A_2}\left(\frac{-\sqrt{3(1+\beta_2-2\alpha_1)}}{3}\right)\] is zero,
	it has three nontrivial subalgebras if and only if none of \[p_{A_2}\left(\frac{\sqrt{3(1+\beta_2-2\alpha_1)}}{3}\right), \ \ p_{A_2}\left(\frac{-\sqrt{3(1+\beta_2-2\alpha_1)}}{3}\right)\] is zero.
	
	In $A_{3}(\beta_1,\beta_2)=\left(
	\begin{array}{cccc}
	0 & 1 & 1 & 0 \\
	\beta _1 & \beta _2 & 1 & -1
	\end{array}
	\right)$ case  $\alpha_4= 0,$ and one has
	\[p_{A_3}(y)=3y^2-y(1+\beta_2)-\beta_1.\]
	Therefore, $A_{3}(\beta_1,\beta_2)$ has only
	two nontrivial subalgebras if and only if \[(\beta_2+1)^2+12 \beta_1=0,\]
	it has three nontrivial subalgebras if and only if \[(\beta_2+1)^2+12 \beta_1 \neq 0.\]
	
	In $A_{4}(\alpha_1,\beta_2)=\left(
	\begin{array}{cccc}
	\alpha _1 & 0 & 0 & 0 \\
	0 & \beta _2 & -\alpha _1+1 & 0
	\end{array}
	\right)$  case $\alpha_4= 0,$ and we have
	\[p_{A_4}(y)=(2\alpha_1-\beta_2-1)y.\]
	Therefore, $A_{4}(\alpha_1,\beta_2)$ has only
	two nontrivial subalgebras if and only if \[2\alpha_1-\beta_2-1 \neq 0,\]
	it has infinitely many nontrivial subalgebras if and only if \[2\alpha_1-\beta_2-1 = 0.\]
	
	For $A_{5}(\alpha_1)=\left(
	\begin{array}{cccc}
	\alpha _1 & 0 & 0 & 0 \\
	1 & 2\alpha _1-1 & -\alpha _1+1 & 0
	\end{array}
	\right)$ case $\alpha_4= 0$ and
	$p_{A_5}(y)=-1$. Therefore, it has only one nontrivial subalgebra, namely, $\mathbb{F}e_2$.
	
	For $ A_{6}(\alpha_1,\beta_1)=\left(
	\begin{array}{cccc}
	\alpha _1 & 0 & 0 & 1 \\
	\beta _1 & -\alpha _1+1 & -\alpha _1 & 0
	\end{array}
	\right)$ case $\alpha_4=1\neq 0$ and
	\[p_{A_6}(y)= y^3+(3\alpha_1-1)y-\beta_1.\]
	Therefore, $A_{6}(\alpha_1,\beta_1)$ has only
	one nontrivial subalgebra if and only if \[3\alpha_1-1 = 0 \mbox{ and } \beta_1=0,\]
	it has two nontrivial subalgebras if and only if \[3\alpha_1-1 \neq 0 \mbox{ and one of } p_{A_6}\left(\frac{\sqrt{3(1-3\alpha_1)}}{3}\right),\ \ p_{A_6}\left(\frac{-\sqrt{3(1-3\alpha_1)}}{3}\right)\] is zero,
	it has three nontrivial subalgebras if and only if none of \[p_{A_6}\left(\frac{\sqrt{3(1-3\alpha_1)}}{3}\right),\ \ p_{A_6}\left(\frac{-\sqrt{3(1-3\alpha_1)}}{3}\right)\] is zero.

	For $ A_{7}(\beta_1)=\left(
	\begin{array}{cccc}
	0 & 1 & 1 & 0 \\
	\beta _1 & 1 & 0 & -1
	\end{array}
	\right)$ case $\alpha_4=0$ and
	\[p_{A_7}(y)=3y^2-y-\beta_1.\]
	Therefore, $A_{7}(\beta_1)$ has only
	two nontrivial subalgebras if and only if \[1+12\beta_1=0,\]
	it has three nontrivial subalgebras if and only if \[1+12\beta_1 \neq 0.\]
	
	In $ A_{8}(\alpha_1)=\left(
	\begin{array}{cccc}
	\alpha _1 & 0 & 0 & 0 \\
	0 & 1-\alpha _1 & -\alpha _1 & 0
	\end{array}
	\right)$ case  $\alpha_4=0$ and
	\[p_{A_8}(y)=(3\alpha_1-1)y.\]
	Therefore, $A_{8}(\alpha_1)$ has only
	two nontrivial subalgebras if and only if \[3\alpha_1-1 \neq 0,\]
	it has infinitely many nontrivial subalgebras if and only if \[3\alpha_1-1 = 0.\]
	
	For $A_{9}=\left(
	\begin{array}{cccc}
	\frac{1}{3} & 0 & 0 & 0 \\
	1 & \frac{2}{3} & \frac{-1}{3} & 0
	\end{array}
	\right) $ we have  $\alpha_4=0$ and
	$p_{A_9} (y)=-1\neq 0$. The algebra $A_{9}$ has only one nontrivial subalgebra, namely,
	$\mathbb{F}e_2.$
	
	For $A_{10}=\left(
	\begin{array}{cccc}
	0 & 1 & 1 & 0 \\
	0 & 0 & 0 & -1
	\end{array}
	\right) $ one has  $\alpha_4=0$ and all nontrivial subalgebras are
	$\mathbb{F}e_2$ and $\mathbb{F}e_1.$
	
	The algebra $A_{11}=\left(
	\begin{array}{cccc}
	0 & 1 & 1 & 0 \\
	1 & 0 & 0 & -1
	\end{array}
	\right) $ has only three nontrivial subalgebras
	\[\mathbb{F}e_2,\ \ \mathbb{F}\left(e_1-
	\frac{e_2}{\sqrt{3}}
	\right)\ \mbox{and}\  \mathbb{F}\left(e_1+\frac{e_2}{\sqrt{3}}\right)
	.\]
	
	For $A_{12}=\left(
	\begin{array}{cccc}
	0 & 0 & 0 & 0 \\
	1 & 0 & 0 & 0
	\end{array}
	\right) $ we have
	$p_{A_{12}}(y)=-1 \neq 0$ and $A_{12}$ has only one nontrivial subalgebra, namely,
	$\mathbb{F}e_2.$
	\end{proof}	
	We present the corresponding results in $Char(\mathbb{F})= 2$ and $Char(\mathbb{F})= 3$ cases as follows.
	
	\begin{theorem} Let $Char(\mathbb{F})=2.$ Two-dimensional algebras over $\mathbb{F}$ with respect to the number of their subalgebras are given as follows.
\begin{itemize}
\item Any two-dimensional algebra over $\mathbb{F}$ with only one subalgebra is isomorphic to one of the following algebras
\begin{itemize}
\item[$\ast$] $A_{1,2}(\alpha_1,\alpha_2,\alpha_4,\beta_1)$, where $\alpha_4(\alpha_1+1)+(\alpha_2+1)^2=0,\  (\alpha_2+1)^3+\beta_1\alpha_4^2=0$ and $\alpha_4\neq 0$;
	\item[$\ast$]	$A_{1,2}(1,1,0,\beta_1)$, where $\beta_1\neq 0;$
	\item[$\ast$]	$A_{2,2}(\alpha_1,0,1)$;
	\item[$\ast$]	$A_{5,2}(\alpha_1)$;
\item[$\ast$] $A_{6,2}(1,0)$;
\item[$\ast$] $A_{9,2}$;
\item[$\ast$] $A_{12,2}.$
\end{itemize}
	\item	Any two-dimensional algebra over $\mathbb{F}$ with two subalgebras is isomorphic to one of the following algebras
\begin{itemize}
\item[$\ast$] $A_{1,2}(\alpha_1,\alpha_2,\alpha_4,\beta_1)$, where  $(\alpha_2+1)^3+\beta_1\alpha_4^2=0$ and  $\alpha_4(\alpha_4(\alpha_1+1)+(\alpha_2+1)^2)\neq 0$;
	\item[$\ast$]	$A_{1,2}(1,\alpha_2,0,\beta_1)$, where  $\alpha_2\neq 1$;
	\item[$\ast$]	$A_{1,2}(\alpha_1,1,0,\beta_1)$, where  $\alpha_1\neq 1$,
\item[$\ast$] $A_{2,2}(\alpha_1,0, \beta_2)$, where $\beta_2\neq 1$;
	\item[$\ast$]	$A_{3,2}(\alpha_1,1)$;
\item[$\ast$] $A_{4,2}(\alpha_1, \beta_2)$, where $\beta_2\neq 1$;
\item[$\ast$] $A_{6,2}(\alpha_1,0)$, where $\alpha_1\neq 1$;
\item[$\ast$] $A_{7,2}(1)$;
 \item[$\ast$] $A_{8,2}(\alpha_1)$, where $\alpha_1\neq 1$;
 \item[$\ast$] $A_{10,2}$.\end{itemize}
\item Any two-dimensional algebra over $\mathbb{F}$ with three subalgebras is isomorphic to one of the following algebras
\begin{itemize}
\item[$\ast$] 	$A_{1,2}(\alpha_1,\alpha_2,\alpha_4,\beta_1)$, where $\alpha_4(\alpha_4\beta_1+(\alpha_1+1)(\alpha_2+1))\neq 0$;
	\item[$\ast$]	$A_{1,2}(\alpha_1,\alpha_2,0,\beta_1)$, where $(\alpha_1+1)(\alpha_2+1)\neq 0$;
	\item[$\ast$]	$A_{2,2}(\alpha_1,\beta_1, \beta_2)$, where $\beta_1\neq 0$;
\item[$\ast$]	$A_{3,2}(\alpha_1, \beta_2)$, where $\beta_2\neq 1$;	
\item[$\ast$]	$A_{6,2}(\alpha_1,\beta_1)$, where $ \beta_1\neq 0;$
	\item[$\ast$]	$A_{7,2}(\alpha_1)$, where $\alpha_1\neq 1$;
\item[$\ast$] $A_{11,2}$.
\end{itemize}
\item Any two-dimensional algebra over $\mathbb{F}$ with infinitely many subalgebras is isomorphic to one of the following algebras
\begin{itemize}
\item[$\ast$] $A_{1,2}(1,1,0,0)$;
\item[$\ast$] $A_{4,2}(\alpha_1,1)$;
\item[$\ast$] $A_{8,2}(1).$
\end{itemize}
\end{itemize}	
	\end{theorem}
\begin{proof}
\begin{proof}
For $A_{1,2}=\left(
\begin{array}{cccc}
 \alpha _1 & \alpha _2 & \alpha _2+1 & \alpha _4 \\
 \beta _1 & -\alpha _1 & -\alpha _1+1 & -\alpha _2
\end{array}
\right)$ we have
\[p_{A_{1,2}}(y)=\alpha_4 y^3+(\alpha_2+1)y^2+(\alpha_1-1)y -\beta_1.\]
In $\alpha_4\neq 0$ case the algebra $A_{1,2}(\alpha_1,\alpha_2,\alpha_4,\beta_1)$ has only
\begin{itemize}
  \item one nontrivial subalgebra if and only if $\alpha_4(\alpha_1-1)-(\alpha_2+1)^2=0 \mbox{ and } (\alpha_2+1)^3-\beta_1\alpha_4^2=0.$
  \item two nontrivial subalgebras if and only if $\alpha_4\beta_1-(\alpha_1-1)(\alpha_2+1)=0$ and non of  $\alpha_4(\alpha_1-1)-(\alpha_2+1)^2 \ , \alpha_4\beta_1^2-(\alpha_1-1)^3$ is zero.
  \item three nontrivial subalgebras if and only $\alpha_4\beta_1-(\alpha_1-1)(\alpha_2+1)\neq 0.$
\end{itemize}

 In $\alpha_4= 0$. The subalgebras are
$\mathbb{F}(e_2)$ and $\mathbb{F}(e_1+y_0e_2)$, where $y_0$ is any root of the polynomial
$p_{A_{1,2}}(y)=y^2(\alpha_2+1)+y(\alpha_1-1)-\beta_1$.\\

 Therefore the algebra $A_{1,2}(\alpha_1,\alpha_2,0,\beta_1)$ has only
 \begin{itemize}
   \item one nontrivial subalgebra if and only if $\alpha_1-1=\alpha_2+1=0\mbox{ and } \beta_1\neq 0.$
   \item two nontrivial subalgebras if and only if $\alpha_2+1\neq 0 \ , (\alpha_1-1)=0,$ or $\alpha_2+1= 0 \mbox{ and } \alpha_1-1\neq 0.$
   \item three nontrivial subalgebras if and only if $\alpha_2+1\neq 0 \ , (\alpha_1-1)\neq 0.$
   \item infinitely many nontrivial $1$-dimensional subalgebras if and only if $\alpha_2+1=\alpha_1-1=\beta_1=0.$
 \end{itemize}

  For $A_{2,2}=\left(
\begin{array}{cccc}
 \alpha _1 & 0 & 0 & 1 \\
 \beta _1 & \beta _2 & -\alpha _1+1 & 0
\end{array}
\right)$ case $\alpha_4\neq 0$ then the nontrivial subalgebras are
  $\mathbb{F}(e_1+y_0e_2)$, where $y_0$ is the roots of the polynomial
  \[ p_{A_{2,2}}(y)=y^3+(\beta_2-1)y-\beta_1.\]

   Therefore the algebra $A_{2,2}(\alpha_1,\beta_1, \beta_2)$ has only
   \begin{itemize}
     \item one nontrivial subalgebra if and only if $\beta_2-1 = 0\mbox{ and } \beta_1= 0.$
     \item two nontrivial subalgebras if and only if $\beta_2-1\neq 0 \mbox{ and }\beta_1=0.$
     \item three nontrivial subalgebras if and only if  $\beta_1 \neq 0.$
   \end{itemize}

  For $A_{3,2}=\left(
\begin{array}{cccc}
 \alpha_1 & 1 & 1 & 0 \\
 0 & \beta _2 & 1-\alpha_1 & -1
\end{array}
\right)$ case $\alpha_4= 0.$\\
The subalgebras are $\mathbb{F}(e_2)$ and $\mathbb{F}(e_1+y_0e_2)$, where $y_0$ is the roots of the polynomial:
 \[y^2-y(1+\beta_2)=0.\]
Therefore $A_{3,2}(\alpha_1,\beta_2)$ has only
\begin{itemize}
  \item two nontrivial subalgebras if and only if $(\beta_2+1)=0.$
  \item three nontrivial subalgebras if and only if $(\beta_2+1) \neq 0.$
\end{itemize}

For $A_{4,2}=\left(
\begin{array}{cccc}
 \alpha _1 & 0 & 0 & 0 \\
 0 & \beta _2 & -\alpha _1+1 & 0
\end{array}
\right)$  case $\alpha_4= 0.$\\
The subalgebras are
$\mathbb{F}(e_2)$ and $\mathbb{F}(e_1+y_0e_2)$, where $y_0$ is the roots of the polynomial
 \[y(\beta_2-1)=0.\]
 Therefore $A_{4,2}(\alpha_1,\beta_2)$ has only
\begin{itemize}
  \item two nontrivial subalgebras if and only if $\beta_2-1 \neq 0.$
  \item infinitely many nontrivial subalgebras if and only if $\beta_2-1 = 0.$
\end{itemize}

For $A_{5,2}=\left(
\begin{array}{cccc}
 \alpha _1 & 0 & 0 & 0 \\
 1 & 1 & -\alpha _1+1 & 0
\end{array}
\right)$ case $\alpha_4= 0.$\\
$p_{A_{5,2}}(y)\neq 0$ then $\mathbb{F}(e_2)$ is only the nontrivial subalgebra of $ A_{5,2}(\alpha_1).$\\

For $A_{6,2}=\left(
\begin{array}{cccc}
 \alpha _1 & 0 & 0 & 1 \\
 \beta _1 & -\alpha _1+1 & -\alpha _1 & 0
\end{array}
\right),$ case $\alpha_4\neq 0.$\\
 The subalgebras are
 $\mathbb{F}(e_1+y_0e_2)$, where $y_0$ is the roots of the polynomial:
\[ y^3+y(\alpha_1-1)-\beta_1=0.\]
 Therefore $A_{6,2}(\alpha_1,\beta_1)$ has only
\begin{itemize}
\item one nontrivial subalgebra if and only if $\alpha_1-1 = 0 \mbox{ and } \beta_1=0.$
  \item two nontrivial subalgebras if and only if $\alpha_1-1 \neq 0 \mbox{ and } \beta_1=0.$
  \item three nontrivial subalgebras if and only if $\beta_1 \neq 0.$
\end{itemize}

For $A_{7,2}=\left(
\begin{array}{cccc}
 \alpha _1 & 1 & 1 & 0 \\
 0 & 1-\alpha _1 & -\alpha _1 & -1
\end{array}
\right)$ case.\\
The subalgebras are
$\mathbb{F}(e_2))$ and $\mathbb{F}(e_1+y_0e_2)$, where $y_0$ is the roots of the polynomial:\\
 \[y^2+y(\alpha_1-1)=0.\]
 Therefore $A_{7,2}(\beta_1)$ has only
\begin{itemize}
  \item two nontrivial subalgebras if and only if $(\alpha_1-1)=0.$
  \item three nontrivial subalgebras if and only if $(\alpha_1-1) \neq 0.$
\end{itemize}

For $A_{8,2}=\left(
\begin{array}{cccc}
 \alpha _1 & 0 & 0 & 0 \\
 0 & 1-\alpha _1 & -\alpha _1 & 0
\end{array}
\right)$ case.\\
 The subalgebras are
$\mathbb{F}(e_2)$ and $\mathbb{F}(e_1+y_0e_2)$, where $y_0$ is the roots of the polynomial:\\
 \[y(\alpha_1-1) =0.\]
 Therefore $A_{8,2}(\alpha_1)$ has only
\begin{itemize}
  \item two nontrivial subalgebras if and only if $\alpha_1-1 \neq 0.$
  \item infinitely nontrivial subalgebras if and only if $\alpha_1-1 = 0.$
\end{itemize}

For $A_{9,2}=\left(
\begin{array}{cccc}
 1 & 0 & 0 & 0 \\
 1 & 0 & 1 & 0
\end{array}
\right),$ case.\\
$p_{A_{9,2}} (y)\neq 0$ then $A_{9,2}$ has only one nontrivial subalgebra
$\mathbb{F}(e_2).$\\

For $A_{10,2}=\left(
\begin{array}{cccc}
 0 & 1 & 1 & 0 \\
 0 & 0 & 0 & -1
\end{array}
\right),$ case.\\
The subalgebras are
$\mathbb{F}(e_2)$ and $\mathbb{F}(e_1).$\\

For $A_{11,2}=\left(
\begin{array}{cccc}
 1 & 1 & 1 & 0 \\
 0 & -1 & -1 & -1
\end{array}
\right),$ case.\\
The subalgebras are
$\mathbb{F}(e_2)$ and $\mathbb{F}(e_1+y_0e_2)$, where $y_0$ is the roots of the polynomial:\\
 \[y^2+y=0.\]
 Then $A_{11,2}$ has three nontrivial subalgebras
 \[\mathbb{F}(e_2)\ , \mathbb{F}(e_1)\ , \mathbb{F}(e_1+y_0e_2)\]

For $A_{12,2}=\left(
\begin{array}{cccc}
 0 & 0 & 0 & 0 \\
 1 & 0 & 0 & 0
\end{array}
\right)$ case.\\
$p_{A_{12,2}}(y) \neq 0$ then $A_{12,2}$ has only one nontrivial subalgebra, namely,
$\mathbb{F}(e_2).$
\end{proof}
		
	\begin{theorem} Let $Char(\mathbb{F})= 3.$ Two-dimensional algebras over $\mathbb{F}$ are described as follows.
\begin{itemize}
\item Any two-dimensional algebra over $\mathbb{F}$ with only one subalgebra is isomorphic to one of the following algebras
\begin{itemize}
\item[$\ast$] $A_{2,3}(\alpha_1,\beta_1,2\alpha_1-1)$;
	\item[$\ast$]	$A_{3,3}(\beta_1,-1)$, where $\beta_1\neq 0$;
	\item[$\ast$]	$A_{5,3}(\alpha_1)$;
\item[$\ast$] $A_{9,3}$;
\item[$\ast$] $A_{11,3}$;
\item[$\ast$] $A_{12,3}.$
\end{itemize}
\item Any two-dimensional algebra over $\mathbb{F}$ with two subalgebras is isomorphic to one of the following algebras
\begin{itemize}
\item[$\ast$]	$A_{1,3}(\alpha_1,\alpha_2,\alpha_4,-\alpha_4-1)$;
	\item[$\ast$]	$A_{3,3}(\beta_1,\beta_2)$, where $\beta_2\neq -1$,
\item[$\ast$] $A_{4,3}(\alpha_1, \beta_2)$, where $2\alpha_1-\beta_2-1\neq 0$,
\item[$\ast$] $A_{7,3}(\beta_1)$;
\item[$\ast$] $A_{8,3}(\alpha_1)$.
\end{itemize}
\item Any two-dimensional algebra over $\mathbb{F}$ with three subalgebras is isomorphic to one of the following algebras
\begin{itemize}
\item[$\ast$]	$A_{1,3}(\alpha_1,\alpha_2,\alpha_4,\beta_1)$, where $\alpha_4(\beta_1+\alpha_4+1)\neq 0$;
	\item[$\ast$]	$A_{1,3}(\alpha_1,\alpha_2,0,\beta_1)$, where $\beta_1\neq -1$;
	\item[$\ast$]	$A_{2,3}(\alpha_1,\beta_1, \beta_2)$, where $\beta_1\neq 2\alpha_1-1$;
	\item[$\ast$]	$A_{6,3}(\alpha_1,\beta_1).$
\end{itemize}
\item Any two-dimensional algebra over $\mathbb{F}$ with infinitely many subalgebras is isomorphic to one of the following algebras
\begin{itemize}
\item[$\ast$]	$A_{3,3}(0,-1)$;
\item[$\ast$] $A_{4,3}(\alpha_1,2\alpha_1-1)$;
\item[$\ast$] $A_{10,3}.$
\end{itemize}
\end{itemize}
	\end{theorem}
\begin{proof}
  For $A_{1,3}(\alpha_1,\alpha_2,\alpha_4,\beta_1)=\left(
\begin{array}{cccc}
 \alpha _1 & \alpha _2 & \alpha _2+1 & \alpha _4 \\
 \beta _1 & -\alpha _1 & -\alpha _1+1 & -\alpha _2
\end{array}
\right)$ case
\[p_{A_{1,3}}(y)=\alpha_4 y^3+y^2-y-\beta_1.\]
If  $\alpha_4\neq 0,$ the algebra $A_{1,3}(\alpha_1,\alpha_2,\alpha_4,\beta_1)$ has only
\begin{itemize}
  \item two nontrivial subalgebras if and only if $\alpha_4+\beta_1+1=0.$
  \item three nontrivial subalgebras if and only if $\alpha_4+\beta_1+1 \neq 0.$
\end{itemize}
If $\alpha_4= 0,$ the subalgebras are
$\mathbb{F}(e_2)$ and $\mathbb{F}(e_1+y_0e_2)$, where $y_0$ is a root of the polynomial
$p_{A_{1,3}}(y)=y^2-y-\beta_1$.\\

 Therefore the algebra $A_{1,3}(\alpha_1,\alpha_2,0,\beta_1)$ has only
 \begin{itemize}
   \item two nontrivial subalgebras if and only if $\beta_1 +1= 0.$
   \item three nontrivial subalgebras if and only if $\beta_1+1 \neq 0.$
 \end{itemize}

  For $A_{2,3}(\alpha_1,\beta_1,\beta_2)=\left(
  \begin{array}{cccc}
  \alpha _1 & 0 & 0 & 1 \\
  \beta _1 & \beta _2 & -\alpha _1+1 & 0
  \end{array}
  \right)$ case
  \[ p_{A_{2,3}}(y)=y^3+(2\alpha_1-\beta_2-1)y-\beta_1.\]

   Therefore the algebra $A_{2,3}(\alpha_1,\beta_1, \beta_2)$ has only
   \begin{itemize}
     \item one nontrivial subalgebra if and only if $2\alpha_1-\beta_2-1=0.$
     \item three nontrivial subalgebras $ 2\alpha_1-\beta_2-1 \neq 0.$
   \end{itemize}

  For $A_{3,3}(\beta_1,\beta_2)=\left(
\begin{array}{cccc}
 0 & 1 & 1 & 0 \\
 \beta _1 & \beta _2 & 1 & -1
\end{array}
\right)$ case $\alpha_4= 0.$\\
The subalgebras are $\mathbb{F}(e_2)$ and $\mathbb{F}(e_1+y_0e_2)$, where $y_0$ is the roots of the polynomial:
 \[-y(1+\beta_2)-\beta_1=0.\]
Therefore $A_{3,3}(\beta_1,\beta_2)$ has only
\begin{itemize}
  \item one nontrivial subalgebra if and only if $\beta_2+1 =0, \ \beta_1 \neq 0.$
  \item two  nontrivial subalgebras if and only if $\beta_2+1 \neq 0.$
  \item infinitely nontrivial subalgebras if and only if $\beta_2+1= \beta_1= 0.$
\end{itemize}

For $A_{4,3}(\alpha_1,\beta_2)=\left(
\begin{array}{cccc}
 \alpha _1 & 0 & 0 & 0 \\
 0 & \beta _2 & -\alpha _1+1 & 0
\end{array}
\right)$  case $\alpha_4= 0.$\\
The subalgebras are
$\mathbb{F}(e_2)$ and $\mathbb{F}(e_1+y_0e_2)$, where $y_0$ is the roots of the polynomial
 \[y(2\alpha_1-\beta_2-1)=0.\]
 Therefore $A_{4,3}(\alpha_1,\beta_2)$ has only
\begin{itemize}
  \item two nontrivial subalgebras if and only if $2\alpha_1-\beta_2-1 \neq 0.$
  \item infinitely nontrivial subalgebras if and only if $2\alpha_1-\beta_2-1 = 0.$
\end{itemize}

For $ A_{5,3}(\alpha_1)=\left(
\begin{array}{cccc}
 \alpha _1 & 0 & 0 & 0 \\
 1 & 2\alpha _1-1 & -\alpha _1+1 & 0
\end{array}
\right)$ case $\alpha_4= 0.$\\
$p_{A_{5,3}}(y)\neq 0$ then $\mathbb{F}(e_2)$ is only the nontrivial subalgebra of $ A_{5,3}(\alpha_1).$\\

For $ A_{6,3}(\alpha_1,\beta_1)=\left(
\begin{array}{cccc}
 \alpha _1 & 0 & 0 & 1 \\
 \beta _1 & -\alpha _1+1 & -\alpha _1 & 0
\end{array}
\right)$ case $\alpha_4\neq 0.$\\
 The subalgebras are
 $\mathbb{F}(e_1+y_0e_2)$, where $y_0$ is the roots of the polynomial:\\
\[ y^3-y-\beta_1=0.\]
 Therefore $A_{6,3}(\alpha_1,\beta_1)$ has only three nontrivial subalgebras.\\

For $ A_{7,3}(\beta_1)=\left(
\begin{array}{cccc}
 0 & 1 & 1 & 0 \\
 \beta _1 & 1 & 0 & -1
\end{array}
\right)$ case.\\
The subalgebras are
$\mathbb{F}(e_2)$ and $\mathbb{F}(e_1+y_0e_2)$, where $y_0$ is the roots of the polynomial:
 \[-y-\beta_1=0.\]
 Therefore $A_{7,3}(\beta_1)$ has only two nontrivial subalgebras $\mathbb{F}(e_2)$ and $\mathbb{F}(e_1+y_0e_2).$\\

For $ A_{8,3}(\alpha_1)=\left(
\begin{array}{cccc}
 \alpha _1 & 0 & 0 & 0 \\
 0 & 1-\alpha _1 & -\alpha _1 & 0
\end{array}
\right)$ case.\\
 The subalgebras are only
$\mathbb{F}(e_2)$ and $\mathbb{F}(e_1).$\\

For $A_{9,3}=\left(
\begin{array}{cccc}
 0 & 1 & 1 & 0 \\
 1 & 0 & 0 & -1
\end{array}
\right) $ case.\\
$p_{A_{9,3}} (y)\neq 0$ then $A_{9,3}$ has only one nontrivial subalgebra
$\mathbb{F}(e_2).$\\

For $A_{10,3}=\left(
\begin{array}{cccc}
 0 & 1 & 1 & 0 \\
 0 & 0 & 0 & -1
\end{array}
\right) $ case.\\
$p_{A_{10,3}} (y)= 0$ for any $y$ then $A_{10,3}$ has infinitely many nontrivial subalgebras.\\

For $A_{11,3}=\left(
\begin{array}{cccc}
 1 & 0 & 0 & 0 \\
 1 & -1 & -1 & 0
\end{array}
\right) $ case.\\
$p_{A_{11,3}} (y)\neq 0$ then $A_{11,3}$ has only one nontrivial subalgebra $\mathbb{F}(e_2).$\\

For $A_{12,3}=\left(
\begin{array}{cccc}
 0 & 0 & 0 & 0 \\
 1 & 0 & 0 & 0
\end{array}
\right) $ case.\\
$p_{A_{12,3}}(y) \neq 0$ then $A_{12,3}$ has only one nontrivial subalgebra
$\mathbb{F}(e_2).$
\end{proof}
\end{proof}

\section{The idempotents}
In this section we describe the idempotent elements of two-dimensional algebras over $\mathbb{F}$ (for assumptions on the field $\mathbb{F}$ see Section 1). Recall that a nonzero element $\mathbf{v}$ of an algebra $\mathbb{A}$ is an idempotent if
\begin{equation} \label{Idem}
\mathbf{v}^2=\mathbf{v}.
\end{equation}
  Let $\mathbb{A}$ be a two-dimensional algebra. If $\mathbf{v}=ye_2$ then (\ref{Idem}) is equivalent to $\mathbf{u}^2=\frac{1}{y}\mathbf{u},$ where  $\mathbf{u}=e_2$ and if $\mathbf{v}=xe_1+ye_2$,\ $x\neq 0$ then (\ref{Idem}) is equivalent to $\mathbf{u}^2=\frac{1}{x}\mathbf{u},$ where $\mathbf{u}=e_1+\frac{y}{x}e_2$. Due to the discussion made at the begining of Section \ref{Sec3} here we are interested in nonzero $\mathbf{u} \in \mathbb{A}$ with  $\mathbf{u}^2=\lambda\mathbf{u},$ $\lambda\neq 0$.

Let  $Id(A)$ stand for the set of all idempotents of the algebra $\mathbb{A}$ with MSC $A=\left(
\begin{array}{cccc}
\alpha _1 & \alpha _2 & \alpha _3 & \alpha _4 \\
\beta _1 & \beta _2 & \beta _3 & \beta _4
\end{array}
\right)$ and $\lambda_A(y)=\alpha_4y^2+(\alpha_2+\alpha_3)y+\alpha_1.$
\begin{proposition} \emph{} \label{Pr41}

\begin{itemize}
\item If  $\alpha_4\neq 0$ or $\alpha_4= 0$ and $\beta_4=0$ then  \[Id(A)=\left\{\mathbf{u}=\frac{1}{\lambda_A(y_0)}(e_1+y_0e_2):\ y_0\in \mathbb{F}, \ p_A(y_0)=0,\ \lambda_A(y_0)\neq 0\right\}.\]
\item If $\alpha_4= 0$ and $\beta_4\neq 0$ then  \[Id(A)=\left\{\mathbf{u}=\frac{1}{\lambda_A(y_0)}(e_1+y_0e_2):\ y_0\in \mathbb{F},\ p_A(y_0)=0,\ \lambda_A(y_0)\neq 0\right\}\cup \left\{\frac{1}{\beta_4}e_2\right\}.\]
\end{itemize}
\end{proposition}
\begin{proof}
Indeed, if $\alpha _4\neq 0$ then the first equation of (\ref{SEs}) implies $\lambda=\alpha_4y_0^2+(\alpha_2+\alpha_3)y_0+\alpha_1,$ where $y_0$ is a root of \[p_A(y)=y^3 \alpha _4+y^2(\alpha _2+\alpha _3-\beta _4) +y(\alpha _1-\beta _2-\beta _3)-\beta _1.\]
Therefore, in this case we have
\[Id(A)=\left\{\mathbf{u}=\frac{1}{\lambda_A(y_0)}(e_1+y_0e_2): y_0\in \mathbb{F}, p_A(y_0)=0,\lambda_A(y_0)\neq 0\right\}.\]

Let us consider $\alpha_4= 0.$ Then one has $e_2^2=\beta_4e_2$ and the element $\frac{1}{\beta_4}e_2$ is idempotent provided that $\beta_4\neq 0$. Moreover, in this case each root $y_0$ of $p_A(y)$ generates an idempotent $\frac{1}{\lambda_A(y_0)}(e_1+y_0e_2)$, whenever $\lambda_A(y_0)\neq 0$.
\end{proof}
Due to Proposition \ref{Pr41} and space reasons we abstain from writing the conditions for $A_1, A_{2}$, $A_6$ in $Char(\mathbb{F})\neq 2,3$ ($A_{1,2}, A_{2,2}$, $A_{6,2}$ in $Char(\mathbb{F})=2$ and  $A_{1,3}, A_{2,3}$, $A_{6,3}$ in $Char(\mathbb{F})=3$) case to have no idempotent, one idempotent, two, three and infinitely many idempotents. For the rest the final result we give as follows.
\begin{theorem}\emph{}

\begin{itemize}
\item Let $Char(\mathbb{F})\neq 2,3$. Then
\begin{itemize}
\item[$\ast$] $Id(A_4(0,\beta_2))=Id(A_5(\alpha_1))=Id(A_8(0))=Id(A_9)=Id(A_{12})=\emptyset;$
\item[$\ast$] $Id(A_3(0,-1))=\{-e_2\};$
\item[$\ast$] $Id(A_3(0,\beta_2))=\left\{-e_2, \frac{3}{2(\beta_2+1)}e_1+\frac{1}{2}e_2\right\},$ where $\beta_2\neq -1$
\item[$\ast$] $Id(A_3(\beta_1,\beta_2))=\left\{-e_2, \frac{3}{\beta_2+1-\sqrt{(\beta_2+1)^2+2\beta_1}}e_1+\frac{1}{2}e_2, \frac{3}{\beta_2+1+\sqrt{(\beta_2+1)^2+2\beta_1}}e_1+\frac{1}{2}e_2\right\},$

    \hfill where $\beta_1\neq 0$
\item[$\ast$] $Id(A_4(\alpha_1,\beta_2))=\left\{\frac{1}{\alpha_1}e_1\right\},$ where $\alpha_1\neq 0,\beta_2\neq 2\alpha_1-1$
\item[$\ast$] $Id(A_4(\alpha_1,2\alpha_1-1))=\left\{ \frac{1}{\alpha_1}e_1+te_2:\ t\in \mathbb{F}\right\},$ where $\alpha_1\neq 0$
\item[$\ast$] $Id(A_7(0))=\left\{-e_2, \frac{3}{2}e_1+\frac{1}{2}e_2\right\},$
\item[$\ast$] $Id(A_7(\beta_1))=\left\{-e_2, \frac{3}{1-\sqrt{12\beta_1+1}}e_1+\frac{1}{2}e_2, \frac{3}{1+\sqrt{12\beta_1+1}}e_1+\frac{1}{2}e_2\right\},$ where $\beta_1\neq 0$
\item[$\ast$] $Id(A_8(\alpha_1))=\left\{ \frac{1}{\alpha_1}e_1\right\},$ where $\alpha_1(3\alpha_1-1)\neq 0$
\item[$\ast$] $Id(A_8\left(\frac{1}{3}\right))=\left\{ 3e_1+te_2:\ t\in \mathbb{F}\right\},$
\item[$\ast$] $Id(A_{10})=\left\{-e_2\right\},$
\item[$\ast$] $Id(A_{11})=\left\{-e_2, \frac{\sqrt{3}}{2}e_1+\frac{1}{2}e_2,\ -\frac{\sqrt{3}}{2}e_1+\frac{1}{2}e_2\right\}.$
\end{itemize}
\item Let $Char(\mathbb{F})= 2.$ Then
\begin{itemize}
\item[$\ast$] $Id(A_{4,2}(0,\beta_2))=Id(A_{5,2}(\alpha))=Id(A_{8,2}(0))
    =Id(A_{9,2})=Id(A_{12,2})=\emptyset;$
\item[$\ast$] $Id(A_{3,2}(0,\beta_2))=\{e_2\};$
\item[$\ast$] $Id(A_{3,2}(\alpha_1,1))=\{e_2, \frac{1}{\alpha_1}e_1\}$ where $\alpha_1 \neq 0;$
\item[$\ast$] $Id(A_{3,2}(\alpha_1,\beta_2))=\{e_2, \frac{1}{\alpha_1}e_1,\frac{1}{\alpha_1}(e_1+(\beta_2+1)e_2) \}$ where $\beta_2\neq 1, \alpha_1 \neq 0;$
\item[$\ast$] $Id(A_{4,2}(\alpha_1,\beta_2))=\left\{\frac{1}{\alpha_1}e_1\right\}\mbox{if}\ \alpha_1\neq 0,\beta_2\neq 1;$
\item[$\ast$] $Id(A_{4,2}(\alpha_1,1))=\left\{ \frac{1}{\alpha_1}e_1+te_2:\ t\in \mathbb{F}\right\},$ where $\alpha_1\neq 0;$
\item[$\ast$] $Id(A_{7,2}(0))=\{-e_2\};$
\item[$\ast$] $Id(A_{7,2}(1))=\{e_1,-e_2\};$
\item[$\ast$] $Id(A_{7,2}(\alpha_1))=\left\{\frac{1}{\alpha_1}e_1, \frac{1}{\alpha_1}e_1+\frac{1-\alpha_1}{\alpha_1}e_2\right\},$ where $\alpha_1\neq 0, 1;$
\item[$\ast$] $Id(A_{8,2}(\alpha_1))=\left\{ \frac{1}{\alpha_1}e_1\right\},$ where $\alpha_1\neq 0,1;$
\item[$\ast$] $Id(A_{8,2}(1))=\{ e_1+te_2:\ t\in \mathbb{F}\};$
\item[$\ast$] $Id(A_{10,2})=\{e_2\};$
\item[$\ast$] $Id(A_{11,2})=\{e_1, -e_1+e_2, -e_2\}.$
\end{itemize}
\item Let $Char(\mathbb{F})= 3.$ Then
\begin{itemize}
\item[$\ast$] $Id(A_{4,3}(0,\beta_2))=Id(A_{5,3}(\alpha_1))=
Id(A_{8,3}(0))=Id(A_{11,3})=Id(A_{12,3})=\emptyset;$
\item[$\ast$] $Id(A_{3,3}(0,\beta_2))=\{-e_2\}$, where $\beta_2\neq -1;$
\item[$\ast$] $Id(A_{3,3}(\beta_1,\beta_2))=\left\{\frac{1+\beta_2}{\beta_1}e_1-e_2,-e_2 \right\}$, where $\beta_1\neq 0;$
 \item[$\ast$] $Id(A_{3,3}(0,-1))=\left\{te_1+\frac{1}{2}e_2,-e_2: t\in \mathbb{F}\right\};$
\item[$\ast$] $Id(A_{4,3}(\alpha_1,\beta_2))=\left\{\frac{1}{\alpha_1}e_1\right\},$ where $\alpha_1\neq 0, \  \beta_2\neq 2\alpha_1-1;$
\item[$\ast$] $Id(A_{4,3}(\alpha_1,2\alpha_1-1))=\left\{ \frac{1}{\alpha_1}e_1+te_2:\ t\in \mathbb{F}\right\},$ where $\alpha_1\neq 0;$
\item[$\ast$] $Id(A_{7,3}(0))=\{-e_2\};$
\item[$\ast$] $Id(A_{7,3}(\beta_1))=\left\{\frac{1}{\beta_1}e_1-e_2,-e_2\right\},$  where $\beta_1\neq 0;$
\item[$\ast$] $Id(A_{8,3}(\alpha_1))=\left\{\frac{1}{\alpha_1}e_1\right\},$ where $\alpha_1\neq0;$
\item[$\ast$] $Id(A_{9,3})=\{-e_2\};$
\item[$\ast$] $Id(A_{10,3})=\{ te_1-e_2:\ t\in \mathbb{F}\}.$
\end{itemize}
\end{itemize}
\end{theorem}
\section{Left Ideals}
A one-dimensional subspace $\mathbb{F}\mathbf{u}$ of an algebra $\mathbb{A}$ is a left ideal if \begin{equation}\label{41} \mathbf{v}\cdot \mathbf{u}= \lambda_l(\mathbf{v}) \mathbf{u} \end{equation} for all $\mathbf{v} \in \mathbb{A}$, where $\lambda_l(\mathbf{v})\in \mathbb{F}$.
In terms of the matrix of structure constants $A$ of a two-dimensional algebra $\mathbb{A}$ due to $\mathbf{v}\cdot \mathbf{u}=\mathrm{e}A(v\otimes u),$ where  $\mathbf{u}=\mathrm{e}u,\mathbf{v}=\mathrm{e}v$ the equation (\ref{41}) is equivalent to
\begin{equation}\label{4n}
A(v\otimes u)- \lambda_l(v) u=0.
\end{equation}
Let $\mathbb{A}$ be a two-dimensional algebra with MSC $A=\left(
\begin{array}{cccc}
\alpha _1 & \alpha _2 & \alpha _3 & \alpha _4 \\
\beta _1 & \beta _2 & \beta _3 & \beta _4
\end{array}
\right)$, $u=\left(
\begin{array}{c}
x \\
y
\end{array}
\right)$, and $v=\left(
\begin{array}{c}
s \\
t
\end{array}
\right)$. Then applying (\ref{4n}) we obtain the system of equations
\begin{equation}\label{SErin}
\begin{array}{rr}
-x \lambda_l(s,t) +s x \alpha _1+s y \alpha _2+t x \alpha _3+t y \alpha _4& =0, \\
-y \lambda_l(s,t) +s x \beta _1+s y \beta _2+t x \beta _3+t y \beta _4&=0.
\end{array}
\end{equation}
If $x=0$ then due to $y\neq 0$ the system (\ref{SErin}) becomes
$$
\begin{array}{rr}
s \alpha _2+t \alpha _4& =0, \\
- \lambda_l(s,t) +s  \beta _2+t  \beta _4&=0.
\end{array}
$$
Evidently, in $\alpha _2= \alpha _4=0$ case $\mathbb{F}e_2$ is a left ideal of $\mathbb{A}$.
If $x\neq 0$ then one can assume that $x=1$. The first equation of
(\ref{SErin}) implies
$\lambda_l(s,t)=s(\alpha_1+y\alpha_2)+t(\alpha_3+y\alpha_4)$ and therefore,  its second equation yields $s(\beta_1+y(\beta_2-\alpha_1)-y^2\alpha_3)+t(\beta_3+y(\beta_4-\alpha_3)-y^2\alpha_4)=0$. So we get the system of equations
\begin{equation}\label{S}
\begin{array}{rr}
l_{1A}(y)=y^2\alpha_4+y(\alpha_3-\beta_4)-\beta_3=0\\
l_{2A}(y)=y^2\alpha_2+y(\alpha_1-\beta_2)-\beta_1=0
\end{array}\end{equation}
and any solution $y_0$ of (\ref{S}) provides the left ideal $\mathbb{F}(e_1+y_0e_2)$. Hence, we need to study the existence of the solutions of (\ref{S}).

\begin{proposition}  Let $Char(\mathbb{F})\neq 2$.

\begin{itemize}
\item  The system $(\ref{S})$ has no solution if and only if one of the following holds
\begin{itemize}
\item[*]	$\alpha_4=\alpha_3-\beta_4=0$ and $\beta_3\neq 0;$
\item[*]	$\alpha_4=0$, $\alpha_3-\beta_4\neq 0$ and $\alpha_2\beta_3^2+\beta_3(\alpha_1-\beta_2)(\alpha_3-\beta_4)-\beta_1(\alpha_3-\beta_4)^2\neq 0;$
\item[*]	$\alpha_4\neq 0$, $(\alpha_3-\beta_4)^2+4\alpha_4\beta_3= 0$ and
	$\alpha_2(\alpha_3-\beta_4)^2-2\alpha_4(\alpha_1-\beta_2)(\alpha_3-\beta_4)-4\beta_1\alpha_4^2\neq 0;$
\item[*] $\alpha_4\neq 0$, $(\alpha_3-\beta_4)^2+4\alpha_4\beta_3\neq 0$, $\alpha_4(\alpha_1-\beta_2)-\alpha_2(\alpha_3-\beta_4)\neq 0$ and

\hfill $l_{1A}\left(\frac{\alpha_4\beta_1-\alpha_2\beta_3}{\alpha_4(\alpha_1-\beta_2)-\alpha_2(\alpha_3-\beta_4)}\right)\neq 0;$
\item[*] $\alpha_4\neq 0$, $(\alpha_3-\beta_4)^2+4\alpha_4\beta_3\neq 0$, $\alpha_4(\alpha_1-\beta_2)-\alpha_2(\alpha_3-\beta_4)= 0$ and
	$\alpha_4\beta_1-\alpha_2\beta_3\neq 0;$
\item[*] $\alpha_4=\alpha_3-\beta_4=\beta_3=\alpha_2=\alpha_1-\beta_2=0$ and $\beta_1\neq 0$.
\end{itemize}
\item	The system $(\ref{S})$ has an unique solution if and only if one of the following holds
\begin{itemize}
\item[*] $\alpha_4=0$, $\alpha_3-\beta_4\neq 0$ and $\alpha_2\beta_3^2+\beta_3(\alpha_1-\beta_2)(\alpha_3-\beta_4)-\beta_1(\alpha_3-\beta_4)^2=0;$
    \item[*] $\alpha_4\neq 0$, $(\alpha_3-\beta_4)^2+4\alpha_4\beta_3= 0$ and
$\alpha_2(\alpha_3-\beta_4)^2-2\alpha_4(\alpha_1-\beta_2)(\alpha_3-\beta_4)-4\beta_1\alpha_4^2=0;$
\item[*]
$\alpha_4\neq 0$, $(\alpha_3-\beta_4)^2+4\alpha_4\beta_3\neq 0$,
	$\alpha_4(\alpha_1-\beta_2)-\alpha_2(\alpha_3-\beta_4)\neq 0$ and

\hfill $l_{1A}\left(\frac{\alpha_4\beta_1-\alpha_2\beta_3}{\alpha_4(\alpha_1-\beta_2)-\alpha_2(\alpha_3-\beta_4)}\right)=0;$
\item[*]
$\alpha_3-\beta_4=\alpha_4=\beta_3=(\alpha_1-\beta_2)^2+4\alpha_2\beta_1= 0$ and $\alpha_2\neq 0;$
\item[*]
$\alpha_3-\beta_4=\alpha_4=\beta_3=\alpha_2= 0$ and $\alpha_1-\beta_2\neq 0.$
\end{itemize}
\item	The system $(\ref{S})$ has two different solutions if and only if one of the following holds
\begin{itemize}
\item[*]	$\alpha_4\neq 0$, $(\alpha_3-\beta_4)^2+4\alpha_4\beta_3\neq  0$ and $\alpha_4(\alpha_1-\beta_2)-\alpha_2(\alpha_3-\beta_4)=\alpha_4\beta_1-\alpha_2\beta_3= 0;$
    \item[*]
	$\alpha_3-\beta_4=\alpha_4=\beta_3=0$
	and $\alpha_2\neq 0$, $(\alpha_1-\beta_2)^2+4\alpha_2\beta_1\neq 0$.
\end{itemize}	
\item	The system $(\ref{S})$ has infinitely many solutions if and only if
\begin{itemize}
\item[*]    $\alpha_1-\beta_2=\alpha_2=\beta_1=\alpha_3-\beta_4=\alpha_4=\beta_3=0$.
\end{itemize}
\end{itemize}   	
\end{proposition}

\begin{proof} The systems of equations (\ref{S}) has no solution, one solution, two solutions and infinitely many solutions is the same for the polynomials $l_{1A}(y)$ and $l_{2A}(y)$ to have no common root, one common root, two common roots and infinitely many common roots. In this study there is a ``symmetricity'' with respect to $l_{1A}(y)$ and $l_{2A}(y)$. To avoid repetitions we give the conditions with respect to one of them (in our case it is $l_{1A}(y)$) the other conditions will come out verifying for $l_{1A}(y)$ and $l_{2A}(y)$ to have no common root, one common root, two common roots and infinitely many common roots. By the way, we keep this strategy further in Propositions \ref{5.2}, \ref{6.1} and \ref{6.2} as well (the later two cases the conditions are given with respect to $r_{1A}(y)$).

The polynomial $l_{1A}(y)$ has no root if and only if $\alpha_4=\alpha_3-\beta_4=0$ and $\beta_3\neq 0$. Therefore,  in this case the system of equations (\ref{S}) is inconsistent.

Let $l_{1A}(y)$ have only one root. The polynomial $l_{1A}(y)$ has only one root $y_0$ if and only if $\alpha_4=0$ and $\alpha_3-\beta_4\neq 0$ or $\alpha_4\neq 0$ and $(\alpha_3-\beta_4)^2+4\alpha_4\beta_3= 0$. In the first case we have $y_0=\frac{\beta_3}{\alpha_3-\beta_4}$. Then $y_0$ satisfies (\ref{S}) if
$\alpha_2\beta_3^2+\beta_3(\alpha_1-\beta_2)(\alpha_3-\beta_4)-\beta_1(\alpha_3-\beta_4)^2=0$ and in the case $\alpha_2\beta_3^2+\beta_3(\alpha_1-\beta_2)(\alpha_3-\beta_4)-\beta_1(\alpha_3-\beta_4)^2\neq 0$ the system has no solution. In the second case
$y_0=-\frac{\alpha_3-\beta_4}{2\alpha_4}$, it satisfies (\ref{S}) if   \[\alpha_2(\alpha_3-\beta_4)^2-2\alpha_4(\alpha_1-\beta_2)(\alpha_3-\beta_4)-4\beta_1\alpha_4^2=0\]  and in the case $\alpha_2(\alpha_3-\beta_4)^2-2\alpha_4(\alpha_1-\beta_2)(\alpha_3-\beta_4)-4\beta_1\alpha_4^2\neq 0$ the system has no solution.

Let $l_{1A}(y)$ have two different roots. The polynomial $l_{1A}(y)$ has two roots if and only if $\alpha_4\neq 0$ and $(\alpha_3-\beta_4)^2+4\alpha_4\beta_3\neq 0$.  Therefore, due to \[l_{2A}(y)=\frac{\alpha_2}{\alpha_4}l_{1A}(y)+\frac{\alpha_4(\alpha_1-\beta_2)-\alpha_2(\alpha_3-\beta_4)}{\alpha_4}y-\frac{\alpha_4\beta_1-\alpha_2\beta_3}{\alpha_4}\] one has the options: under $\alpha_4(\alpha_1-\beta_2)-\alpha_2(\alpha_3-\beta_4)\neq 0$ the system of equations (\ref{S}) has only one solution, namely $y_0=\frac{\alpha_4\beta_1-\alpha_2\beta_3}{\alpha_4(\alpha_1-\beta_2)-\alpha_2(\alpha_3-\beta_4)}$, if  $l_{1A}(\frac{\alpha_4\beta_1-\alpha_2\beta_3}{\alpha_4(\alpha_1-\beta_2)-\alpha_2(\alpha_3-\beta_4)})=0$ and the system of equations (\ref{S}) has no solution if $l_{1A}(\frac{\alpha_4\beta_1-\alpha_2\beta_3}{\alpha_4(\alpha_1-\beta_2)-\alpha_2(\alpha_3-\beta_4)})\neq 0$. If $\alpha_4(\alpha_1-\beta_2)-\alpha_2(\alpha_3-\beta_4)= 0$ then (\ref{S})
has no solution in the case $\alpha_4\beta_1-\alpha_2\beta_3\neq 0$  and the system has two different solutions if $\alpha_4\beta_1-\alpha_2\beta_3= 0$.

Let $l_{1A}(y)$ have infinitely many roots. The polynomial $l_1(y)$ has infinitely many roots if and only if
$\alpha_3-\beta_4=\alpha_4=\beta_3=0$ and therefore, the system has no solution in $\alpha_2=\alpha_1-\beta_2=0$ and $\beta_1\neq 0$ case, it has only one solution in  $\alpha_2=0$, $\alpha_1-\beta_2\neq 0$ and $\alpha_2\neq 0$, $(\alpha_1-\beta_2)^2+4\alpha_2\beta_1= 0$ cases, it has two different solutions in $\alpha_2\neq 0$, $(\alpha_1-\beta_2)^2+4\alpha_2\beta_1\neq 0$ case, it has infinitely many solutions in
$\alpha_3-\beta_4=\alpha_4=\beta_3=\alpha_1-\beta_2=\alpha_2=\beta_1=0$ case.
\end{proof}
In $Char(\mathbb{F})=2$ we make use the following proposition.
\begin{proposition} \label{5.2} \emph{}
\begin{itemize}
\item The system $(\ref{S})$ has no solution if and only if one of the following holds
\begin{itemize}
\item[$\ast$] $\alpha_4=\alpha_3-\beta_4=0$ and $\beta_3\neq 0;$
\item[$\ast$]	$\alpha_4=0$, $\alpha_3-\beta_4\neq 0$ and $\alpha_2\beta_3^2+\beta_3(\alpha_1-\beta_2)(\alpha_3-\beta_4)
    -\beta_1(\alpha_3-\beta_4)^2\neq 0;$
\item[$\ast$]
$\alpha_4\neq 0$, $\alpha_3-\beta_4= 0$ and
	$\alpha_2^2\beta_3^2+\beta_3\alpha_4(\alpha_1-\beta_2)^2-\beta_1^2\alpha_4^2\neq 0;$
\item[$\ast$]
$\alpha_4\neq 0$, $\alpha_3-\beta_4\neq 0$, $\alpha_4(\alpha_1-\beta_2)-\alpha_2(\alpha_3-\beta_4)\neq 0$ and
$l_{1A}\left(\frac{\alpha_4\beta_1-\alpha_2\beta_3}{\alpha_4(\alpha_1-\beta_2)
-\alpha_2(\alpha_3-\beta_4)}\right)\neq 0;$
\item[$\ast$]
$\alpha_4\neq 0$, $\alpha_3-\beta_4\neq 0$, $\alpha_4(\alpha_1-\beta_2)-\alpha_2(\alpha_3-\beta_4)= 0$ and
	$\alpha_4\beta_1-\alpha_2\beta_3\neq 0;$
\item[$\ast$]
	$\alpha_4=\alpha_3-\beta_4=\beta_3=\alpha_2=\alpha_1-\beta_2=0$ and $\beta_1\neq 0.$
\end{itemize}
\item	The system $(\ref{S})$ has unique solution if and only if one of the following holds
\begin{itemize}
\item[$\ast$]	$\alpha_4=0$, $\alpha_3-\beta_4\neq 0$ and $\alpha_2\beta_3^2+\beta_3(\alpha_1-\beta_2)(\alpha_3-\beta_4)-\beta_1(\alpha_3-\beta_4)^2=0$
\item[$\ast$]
	$\alpha_4\neq 0$, $\alpha_3-\beta_4= 0$ and
	$\alpha_2^2\beta_3^2+\beta_3\alpha_4(\alpha_1-\beta_2)^2-\beta_1^2\alpha_4^2=0$
\item[$\ast$]
	$\alpha_4\neq 0$, $\alpha_3-\beta_4\neq 0$,
	$\alpha_4(\alpha_1-\beta_2)-\alpha_2(\alpha_3-\beta_4)\neq 0$ and $l_{1A}\left(\frac{\alpha_4\beta_1-\alpha_2\beta_3}{\alpha_4(\alpha_1-\beta_2)
-\alpha_2(\alpha_3-\beta_4)}\right))=0$
\item[$\ast$]
	$\alpha_3-\beta_4=\alpha_4=\beta_3=\alpha_1-\beta_2= 0$ and $\alpha_2\neq 0$
\item[$\ast$]
	$\alpha_3-\beta_4=\alpha_4=\beta_3=\alpha_2= 0$ and $\alpha_1-\beta_2\neq 0.$
\end{itemize}	
\item	The system $(\ref{S})$ has two different solutions if and only if one of the following holds
	\begin{itemize}
\item[$\ast$] $\alpha_4\neq 0$, $\alpha_3-\beta_4\neq  0$ and $\alpha_4(\alpha_1-\beta_2)-\alpha_2(\alpha_3-\beta_4)=\alpha_4\beta_1-\alpha_2\beta_3= 0$
\item[$\ast$]
	$\alpha_3-\beta_4=\alpha_4=\beta_3=0$
	and $\alpha_2\neq 0$, $\alpha_1-\beta_2\neq 0$.
\end{itemize}
\item	The system $(\ref{S})$ has infinitely many solutions if and only if
 \begin{itemize}
\item[$\ast$]   $\alpha_1-\beta_2=\alpha_2=\beta_1=\alpha_3-\beta_4=\alpha_4=\beta_3=0$.
\end{itemize}
\end{itemize}	
\end{proposition}
\begin{proof}
   \underline{$l_1(y)$ has no root case.} The polynomial $l_1(y)$ has no root if and only if 	\\
$\alpha_4=\alpha_3-\beta_4=0$ and $\beta_3\neq 0$. Therefore (\ref{S}) has no root in this case.

\underline{ $l_1(y)$ has only one root case.} The polynomial $l_1(y)$ has only one root $y_0$ if and only if $\alpha_4=0$ and $\alpha_3-\beta_4\neq 0$ or $\alpha_4\neq 0$ and $\alpha_3-\beta_4= 0$. In the first case $y_0=\frac{\beta_3}{\alpha_3-\beta_4}$, it satisfies (\ref{S}) if
$\alpha_2\beta_3^2+\beta_3(\alpha_1-\beta_2)(\alpha_3-\beta_4)-\beta_1(\alpha_3-\beta_4)^2=0$ and in $\alpha_2\beta_3^2+\beta_3(\alpha_1-\beta_2)(\alpha_3-\beta_4)-\beta_1(\alpha_3-\beta_4)^2\neq 0$ case the system has no solution. In the second case
$y_0=\sqrt{\frac{\beta_3}{\alpha_4}}$, it satisfies (\ref{S}) if   \[\alpha_2^2\beta_3^2+\beta_3\alpha_4(\alpha_1-\beta_2)^2-\beta_1^2\alpha_4^2=0\]  and in \[\alpha_2^2\beta_3^2+\beta_3\alpha_4(\alpha_1-\beta_2)^2-\beta_1^2\alpha_4^2\neq 0\] case the system has no solution.

\underline{  $l_1(y)$ has two different roots case.} The polynomial $l_1(y)$ has two roots if and only if $\alpha_4\neq 0$ and $\alpha_3-\beta_4\neq 0$  Therefore due to \[l_2(y)=\frac{\alpha_2}{\alpha_4}l_1(y)+\frac{\alpha_4(\alpha_1-\beta_2)-\alpha_2(\alpha_3-\beta_4)}{\alpha_4}y-\frac{\alpha_4\beta_1-\alpha_2\beta_3}{\alpha_4}\] if $\alpha_4(\alpha_1-\beta_2)-\alpha_2(\alpha_3-\beta_4)\neq 0$ the system (\ref{S}) has only one solution, namely $y_0=\frac{\alpha_4\beta_1-\alpha_2\beta_3}{\alpha_4(\alpha_1-\beta_2)-\alpha_2(\alpha_3-\beta_4)}$, if  $l_1(\frac{\alpha_4\beta_1-\alpha_2\beta_3}{\alpha_4(\alpha_1-\beta_2)-\alpha_2(\alpha_3-\beta_4)})=0$, if $l_1(\frac{\alpha_4\beta_1-\alpha_2\beta_3}{\alpha_4(\alpha_1-\beta_2)-\alpha_2(\alpha_3-\beta_4)})\neq 0$ then the system has no solution. If $\alpha_4(\alpha_1-\beta_2)-\alpha_2(\alpha_3-\beta_4)= 0$ then (\ref{S})
has no solution in $\alpha_4\beta_1-\alpha_2\beta_3\neq 0$ case, and in $\alpha_4\beta_1-\alpha_2\beta_3= 0$ case the system has two different solutions.

\underline{$l_1(y)$ has infinitely roots.} The polynomial $l_1(y)$ has infinitely many roots if and only if
$\alpha_3-\beta_4=\alpha_4=\beta_3=0$ and therefore the system has no solution in $\alpha_2=\alpha_1-\beta_2=0$ and $\beta_1\neq 0$ case, it has only one solution in  $\alpha_2=0$, $\alpha_1-\beta_2\neq 0$ and $\alpha_2\neq 0$, $\alpha_1-\beta_2= 0$ cases, it has two different solutions in $\alpha_2\neq 0$, $\alpha_1-\beta_2\neq 0$ case, it has infinitely many solutions in
\[\alpha_3-\beta_4=\alpha_4=\beta_3=\alpha_1-\beta_2=\alpha_2=\beta_1=0\] case.
\end{proof}
We summarize the results on left ideals in Tables \ref{T1} and \ref{T2} (See Section APPENDIX).

\section{Right Ideals}
Recall that a one-dimensional subspace $\mathbb{F}\mathbf{u}$, where $0\neq \mathbf{u}\in\mathbb{A}$, is a right ideal of an algebra $\mathbb{A}$ if
\begin{equation} \label{61} \mathbf{u}\cdot\mathbf{v}= \lambda_r(\mathbf{v}) \mathbf{u}\end{equation} for all $\mathbf{v} \in \mathbb{A},$ where $\lambda_r(\mathbf{v})\in \mathbb{F}.$
 In two-dimensional case due to $\mathbf{u}\cdot\mathbf{v}=\mathrm{e}A(u\otimes v),$ where  $\mathbf{u}=\mathrm{e}u,\mathbf{v}=\mathrm{e}v$ we get
\begin{equation}\label{4m}
A(v\otimes u)- \lambda_r(v) u=0.
\end{equation}
Let $A=\left(
\begin{array}{cccc}
\alpha _1 & \alpha _2 & \alpha _3 & \alpha _4 \\
\beta _1 & \beta _2 & \beta _3 & \beta _4
\end{array}
\right)$, $u=\left(
\begin{array}{c}
x \\
y
\end{array}
\right)$, and $v=\left(
\begin{array}{c}
s \\
t
\end{array}
\right)$. Then the equation (\ref{4m}) is written as a system of equations as follows
\begin{equation}\label{SErim}
\begin{array}{rr}
-x \lambda_r(s,t) +s x \alpha _1+xt\alpha _2+ys\alpha _3+t y \alpha _4& =0, \\
-y \lambda_r(s,t) +s x \beta _1+xt\beta _2+ys\beta _3+t y \beta _4&=0.
\end{array}
\end{equation}
Let us consider two cases again. Let $x=0.$ Then due to $y\neq 0$ the system (\ref{SErim}) has the form
 $$
\begin{array}{rr}
s \alpha _3+t \alpha _4& =0, \\
- \lambda_r(s,t) +s  \beta _3+t  \beta _4&=0.
\end{array}
$$
Therefore, only in $\alpha _3= \alpha _4=0$ case $\mathbb{F}e_2$ is a right ideal of $\mathbb{A}$.
If $x\neq 0$ then one can assume that $x=1$. The first equation of
(\ref{SErim}) implies
$\lambda_r(s,t)=s(\alpha_1+y\alpha_3)+t(\alpha_2+y\alpha_4)$ and therefore,  its second equation gives $s(\beta_1+y(\beta_3-\alpha_1)-y^2\alpha_3)+t(\beta_3+y(\beta_2-\alpha_2)-y^2\alpha_4)=0$. Hence, we get the system of equations
\begin{equation}\label{Srim}
\begin{array}{rr}
r_{1A}(y)=y^2\alpha_4+y(\alpha_2-\beta_4)-\beta_2=0\\
r_{2A}(y)=y^2\alpha_3+y(\alpha_1-\beta_3)-\beta_1=0
\end{array}\end{equation}
and any its solution $y_0$ generates a right ideal in the form $\mathbb{F}(e_1+y_0e_2)$.

Further we make use the following two propositions
\begin{proposition}\label{6.1} Let $Char(\mathbb{F})\neq 2$.

\begin{itemize}
\item The system $(\ref{Srim})$ has no solution if and only if
\begin{itemize}
\item[$\ast$]
	$\alpha_4=\alpha_2-\beta_4=0$ and $\beta_2\neq 0$
\item[$\ast$]
	$\alpha_4=0$, $\alpha_2-\beta_4\neq 0$ and $\alpha_3\beta_2^2+\beta_2(\alpha_1-\beta_3)(\alpha_2-\beta_4)-\beta_1(\alpha_2-\beta_4)^2\neq 0$
\item[$\ast$]
	$\alpha_4\neq 0$, $(\alpha_2-\beta_4)^2+4\alpha_4\beta_2= 0$ and
	$\alpha_3(\alpha_2-\beta_4)^2-2\alpha_4(\alpha_1-\beta_3)(\alpha_2-\beta_4)-4\beta_1\alpha_4^2\neq 0$
\item[$\ast$]
	$\alpha_4\neq 0$, $(\alpha_2-\beta_4)^2+4\alpha_4\beta_2\neq 0$, $\alpha_4(\alpha_1-\beta_3)-\alpha_3(\alpha_2-\beta_4)\neq 0$ and

	\hfill $r_{1A}\left(\frac{\alpha_4\beta_1-\alpha_3\beta_2}{\alpha_4(\alpha_1-\beta_3)-\alpha_3(\alpha_2-\beta_4)}\right)\neq 0$
\item[$\ast$]
	$\alpha_4\neq 0$, $(\alpha_2-\beta_4)^2+4\alpha_4\beta_2\neq 0$, $\alpha_4(\alpha_1-\beta_3)-\alpha_3(\alpha_2-\beta_4)= 0$ and
	$\alpha_4\beta_1-\alpha_3\beta_2\neq 0$
\item[$\ast$]
	$\alpha_4=\alpha_2-\beta_4=\beta_2=\alpha_3=\alpha_1-\beta_3=0$ and $\beta_1\neq 0$.
\end{itemize}	
\item	The system $(\ref{Srim})$ has unique solution if and only if\begin{itemize}
\item[$\ast$] $\alpha_4=0$, $\alpha_2-\beta_4\neq 0$ and $\alpha_3\beta_2^2+\beta_2(\alpha_1-\beta_3)(\alpha_2-\beta_4)-\beta_1(\alpha_2-\beta_4)^2=0$
\item[$\ast$]
	$\alpha_4\neq 0$, $(\alpha_2-\beta_4)^2+4\alpha_4\beta_2= 0$ and
	$\alpha_3(\alpha_2-\beta_4)^2-2\alpha_4(\alpha_1-\beta_3)(\alpha_2-\beta_4)-4\beta_1\alpha_4^2=0$
\item[$\ast$]
	$\alpha_4\neq 0$, $(\alpha_2-\beta_4)^2+4\alpha_4\beta_2\neq 0$,
	$\alpha_4(\alpha_1-\beta_3)-\alpha_3(\alpha_2-\beta_4)\neq 0$ and

\hfill $r_{1A}\left(\frac{\alpha_4\beta_1-\alpha_3\beta_2}{\alpha_4(\alpha_1-\beta_3)-\alpha_3(\alpha_2-\beta_4)}\right)=0$
\item[$\ast$]
	$\alpha_2-\beta_4=\alpha_4=\beta_2=(\alpha_1-\beta_3)^2+4\alpha_3\beta_1= 0$ and $\alpha_3\neq 0$
\item[$\ast$]
	$\alpha_2-\beta_4=\alpha_4=\beta_2=\alpha_3= 0$ and $\alpha_1-\beta_3\neq 0$
\end{itemize}	
	
\item	The system $(\ref{Srim})$ has two different solutions if and only if
\begin{itemize}
\item[$\ast$]	$\alpha_4\neq 0$, $(\alpha_2-\beta_4)^2+4\alpha_4\beta_2\neq  0$ and $\alpha_4(\alpha_1-\beta_3)-\alpha_3(\alpha_2-\beta_4)=\alpha_4\beta_1-\alpha_3\beta_2= 0$
    \item[$\ast$]
	$\alpha_2-\beta_4=\alpha_4=\beta_2=0$
	and $\alpha_3\neq 0$, $(\alpha_1-\beta_3)^2+4\alpha_3\beta_1\neq 0$.
\end{itemize}	
\item	The system $(\ref{Srim})$ has infinitely many solutions if and only if
\begin{itemize}
\item[$\ast$] $\alpha_1-\beta_3=\alpha_3=\beta_1=\alpha_2-\beta_4=\alpha_4=\beta_2=0$.
\end{itemize}
\end{itemize}	
\end{proposition}

\begin{proof} Let $r_{1A}(y)$ have no root. The polynomial $r_{1A}(y)$ has no root if and only if $\alpha_4=\alpha_2-\beta_4=0$ and $\beta_2\neq 0$. Therefore, (\ref{Srim}) has no root in this case.

Let $r_{1A}(y)$ have only one root. The polynomial $r_{1A}(y)$ has only one root $y_0$ if and only if $\alpha_4=0$ and $\alpha_2-\beta_4\neq 0$ or $\alpha_4\neq 0$ and $(\alpha_2-\beta_4)^2+4\alpha_4\beta_2= 0$. In the first case $y_0=\frac{\beta_2}{\alpha_2-\beta_4}$, it satisfies (\ref{Srim}) if
$\alpha_3\beta_2^2+\beta_2(\alpha_1-\beta_3)(\alpha_2-\beta_4)-\beta_1(\alpha_2-\beta_4)^2=0$. 
In the second case
$y_0=-\frac{\alpha_2-\beta_4}{2\alpha_4}$ and it satisfies (\ref{Srim}) if   \[\alpha_3(\alpha_2-\beta_4)^2-2\alpha_4(\alpha_1-\beta_3)(\alpha_2-\beta_4)-4\beta_1\alpha_4^2=0\].

Let  $r_{1A}(y)$ have two different roots. The polynomial $r_{1A}(y)$ has two roots if and only if $\alpha_4\neq 0$ and $(\alpha_2-\beta_4)^2+4\alpha_4\beta_2\neq 0$.  Therefore, due to \[r_{2A}(y)=\frac{\alpha_3}{\alpha_4}r_{1A}(y)+\frac{\alpha_4(\alpha_1-\beta_3)-\alpha_3(\alpha_2-\beta_4)}{\alpha_4}y-\frac{\alpha_4\beta_1-\alpha_3\beta_2}{\alpha_4}\] if $\alpha_4(\alpha_1-\beta_3)-\alpha_3(\alpha_2-\beta_4)\neq 0$ the system (\ref{Srim}) has only one solution, namely $y_0=\frac{\alpha_4\beta_1-\alpha_3\beta_2}{\alpha_4(\alpha_1-\beta_3)-\alpha_3(\alpha_2-\beta_4)}$, if  $r_{1A}\left(\frac{\alpha_4\beta_1-\alpha_3\beta_2}{\alpha_4(\alpha_1-\beta_3)-\alpha_3(\alpha_2-\beta_4)}\right)=0$.
 If $\alpha_4(\alpha_1-\beta_3)-\alpha_3(\alpha_2-\beta_4)= 0$ then (\ref{Srim})
has no solution in $\alpha_4\beta_1-\alpha_3\beta_2\neq 0$ case, and in $\alpha_4\beta_1-\alpha_3\beta_2= 0$ case the system has two different solutions.

Let $r_{1A}(y)$ have infinitely many roots. The polynomial $r_{1A}(y)$ has infinitely many roots if and only if
$\alpha_2-\beta_4=\alpha_4=\beta_2=0$ and therefore, the system has no solution in $\alpha_3=\alpha_1-\beta_3=0$ and $\beta_1\neq 0$ case, it has only one solution in  $\alpha_3=0$, $\alpha_1-\beta_3\neq 0$ and $\alpha_3\neq 0$, $(\alpha_1-\beta_3)^2+4\alpha_3\beta_1= 0$ cases, it has two different solutions in $\alpha_3\neq 0$, $(\alpha_1-\beta_3)^2+4\alpha_3\beta_1\neq 0$ case, it has infinitely many solutions in
\[\alpha_2-\beta_4=\alpha_4=\beta_2=\alpha_1-\beta_3=\alpha_3=\beta_1=0\] case.
\end{proof}

\begin{proposition} \label{6.2} Let $Char(\mathbb{F})=2.$
\begin{itemize}
\item The system $(\ref{Srim})$ has no solution if and only if one of the following holds
\begin{itemize}
\item[*] $\alpha_4=\alpha_2-\beta_4=0$ and $\beta_2\neq 0;$
\item[*]	$\alpha_4=0$, $\alpha_2-\beta_4\neq 0$ and $\alpha_3\beta_2^2+\beta_2(\alpha_1-\beta_3)(\alpha_2-\beta_4)
    -\beta_1(\alpha_2-\beta_4)^2\neq 0;$
\item[*]
$\alpha_4\neq 0$, $\alpha_2-\beta_4= 0$ and
	$\alpha_3^2\beta_2^2+\beta_2\alpha_4(\alpha_1-\beta_3)^2-\beta_1^2\alpha_4^2\neq 0;$
\item[*]
$\alpha_4\neq 0$, $\alpha_2-\beta_4\neq 0$, $\alpha_4(\alpha_1-\beta_3)-\alpha_3(\alpha_2-\beta_4)\neq 0$ and
$r_{1A}\left(\frac{\alpha_4\beta_1-\alpha_3\beta_2}{\alpha_4(\alpha_1-\beta_3)
-\alpha_3(\alpha_2-\beta_4)}\right)\neq 0;$
\item[*]
$\alpha_4\neq 0$, $\alpha_2-\beta_4\neq 0$, $\alpha_4(\alpha_1-\beta_3)-\alpha_3(\alpha_2-\beta_4)= 0$ and
	$\alpha_4\beta_1-\alpha_3\beta_2\neq 0;$
\item[*]
	$\alpha_4=\alpha_2-\beta_4=\beta_2=\alpha_3=\alpha_1-\beta_3=0$ and $\beta_1\neq 0.$
\end{itemize}
\item	The system $(\ref{Srim})$ has unique solution if and only if one of the following holds
\begin{itemize}
\item[*]	$\alpha_4=0$, $\alpha_2-\beta_4\neq 0$ and $\alpha_3\beta_2^2+\beta_2(\alpha_1-\beta_3)(\alpha_2-\beta_4)-\beta_1(\alpha_2-\beta_4)^2=0$
\item[*]
	$\alpha_4\neq 0$, $\alpha_2-\beta_4= 0$ and
	$\alpha_3^2\beta_2^2+\beta_2\alpha_4(\alpha_1-\beta_3)^2-\beta_1^2\alpha_4^2=0$
\item[*]
	$\alpha_4\neq 0$, $\alpha_2-\beta_4\neq 0$,
	$\alpha_4(\alpha_1-\beta_3)-\alpha_3(\alpha_2-\beta_4)\neq 0$ and $r_{1A}\left(\frac{\alpha_4\beta_1-\alpha_3\beta_2}{\alpha_4(\alpha_1-\beta_3)
-\alpha_3(\alpha_2-\beta_4)}\right))=0$
\item[*]
	$\alpha_2-\beta_4=\alpha_4=\beta_2=\alpha_1-\beta_3= 0$ and $\alpha_3\neq 0$
\item[*]
	$\alpha_2-\beta_4=\alpha_4=\beta_2=\alpha_3= 0$ and $\alpha_1-\beta_3\neq 0.$
\end{itemize}	
\item	The system $(\ref{Srim})$ has two different solutions if and only if one of the following holds
	\begin{itemize}
\item[*] $\alpha_4\neq 0$, $\alpha_2-\beta_4\neq  0$ and $\alpha_4(\alpha_1-\beta_3)-\alpha_3(\alpha_2-\beta_4)=\alpha_4\beta_1-\alpha_3\beta_2= 0$
\item[*]
	$\alpha_2-\beta_4=\alpha_4=\beta_2=0$
	and $\alpha_3\neq 0$, $\alpha_1-\beta_3\neq 0$.
\end{itemize}
\item	The system $(\ref{Srim})$ has infinitely many solutions if and only if
 \begin{itemize}
\item[*]   $\alpha_1-\beta_3=\alpha_3=\beta_1=\alpha_2-\beta_4=\alpha_4=\beta_2=0$.
\end{itemize}
\end{itemize}	
\end{proposition}
\begin{proof}
  \underline{$l_1(y)$ has no root case.} The polynomial $l_1(y)$ has no root if and only if 	\\
$\alpha_4=\alpha_2-\beta_4=0$ and $\beta_2\neq 0$. Therefore (\ref{Srim}) has no root in this case.

\underline{ $l_1(y)$ has only one root case.} The polynomial $l_1(y)$ has only one root $y_0$ if and only if $\alpha_4=0$ and $\alpha_2-\beta_4\neq 0$ or $\alpha_4\neq 0$ and $\alpha_2-\beta_4= 0$. In the first case $y_0=\frac{\beta_2}{\alpha_2-\beta_4}$, it satisfies (\ref{Srim}) if
$\alpha_3\beta_2^2+\beta_2(\alpha_1-\beta_3)(\alpha_2-\beta_4)-\beta_1(\alpha_2-\beta_4)^2=0$ and in $\alpha_3\beta_2^2+\beta_2(\alpha_1-\beta_3)(\alpha_2-\beta_4)-\beta_1(\alpha_2-\beta_4)^2\neq 0$ case the system has no solution. In the second case
$y_0=\sqrt{\frac{\beta_2}{\alpha_4}}$, it satisfies (\ref{Srim}) if   \[\alpha_3^2\beta_2^2+\beta_2\alpha_4(\alpha_1-\beta_3)^2-\beta_1^2\alpha_4^2=0\]  and in \[\alpha_3^2\beta_2^2+\beta_2\alpha_4(\alpha_1-\beta_3)^2-\beta_1^2\alpha_4^2\neq 0\] case the system has no solution.

\underline{  $l_1(y)$ has two different roots case.} The polynomial $l_1(y)$ has two roots if and only if $\alpha_4\neq 0$ and $\alpha_2-\beta_4\neq 0$  Therefore due to \[l_2(y)=\frac{\alpha_3}{\alpha_4}l_1(y)+\frac{\alpha_4(\alpha_1-\beta_3)-\alpha_3(\alpha_2-\beta_4)}{\alpha_4}y-\frac{\alpha_4\beta_1-\alpha_3\beta_2}{\alpha_4}\] if $\alpha_4(\alpha_1-\beta_3)-\alpha_3(\alpha_2-\beta_4)\neq 0$ the system (\ref{Srim}) has only one solution, namely $y_0=\frac{\alpha_4\beta_1-\alpha_3\beta_2}{\alpha_4(\alpha_1-\beta_3)-\alpha_3(\alpha_2-\beta_4)}$, if  $l_1(\frac{\alpha_4\beta_1-\alpha_3\beta_2}{\alpha_4(\alpha_1-\beta_3)-\alpha_3(\alpha_2-\beta_4)})=0$, if $l_1(\frac{\alpha_4\beta_1-\alpha_3\beta_2}{\alpha_4(\alpha_1-\beta_3)-\alpha_3(\alpha_2-\beta_4)})\neq 0$ then the system has no solution. If $\alpha_4(\alpha_1-\beta_3)-\alpha_3(\alpha_2-\beta_4)= 0$ then (\ref{Srim})
has no solution in $\alpha_4\beta_1-\alpha_3\beta_2\neq 0$ case, and in $\alpha_4\beta_1-\alpha_3\beta_2= 0$ case the system has two different solutions.

\underline{$l_1(y)$ has infinitely roots.} The polynomial $l_1(y)$ has infinitely many roots if and only if
$\alpha_2-\beta_4=\alpha_4=\beta_2=0$ and therefore the system has no solution in $\alpha_3=\alpha_1-\beta_3=0$ and $\beta_1\neq 0$ case, it has only one solution in  $\alpha_3=0$, $\alpha_1-\beta_3\neq 0$ and $\alpha_3\neq 0$, $\alpha_1-\beta_3= 0$ cases, it has two different solutions in $\alpha_3\neq 0$, $\alpha_1-\beta_3\neq 0$ case, it has infinitely many solutions in
\[\alpha_2-\beta_4=\alpha_4=\beta_2=\alpha_1-\beta_3=\alpha_3=\beta_1=0\] case.
\end{proof}
The results on right ideals are given in Tables \ref{T3} and \ref{T4} (see APPENDIX).

\section{Two-sided Ideals}
A one-dimensional subspace $\mathbb{F}\mathbf{u}$, where $0\neq \mathbf{u}\in\mathbb{A}$, is a two-sided ideal of an algebra $\mathbb{A}$ if
\[\mathbf{v}\cdot \mathbf{u}= \lambda_l(\mathbf{v}) \mathbf{u},\ \mathbf{u}\cdot\mathbf{v}= \lambda_r(\mathbf{v}) \mathbf{u}\ \mbox{for all}\ \mathbf{v} \in \mathbb{A}.\]
If $\mathbb{A}$ a two-dimensional with a basis $\mathrm{e}=\{e_1,e_2\}$ then by using \[\mathbf{u}\cdot\mathbf{v}=\mathrm{e}A(u\otimes v),\] where  $\mathbf{u}=\mathrm{e}u,\mathbf{v}=\mathrm{e}v$, we come to
\[
A(v\otimes u)- \lambda_l(v) u=0,\ A(u\otimes v)- \lambda_r(v) u=0
\] and it yields the system of equations
\begin{equation}\label{SEri}
\begin{array}{rrl}
-x \lambda_l(s,t) +s x \alpha _1+s y \alpha _2+t x \alpha _3+t y \alpha _4 & =&0 \\
-y \lambda_l(s,t) +s x \beta _1+s y \beta _2+t x \beta _3+t y \beta _4 & =&0\\
-x \lambda_r(s,t) +s x \alpha _1+xt\alpha _2+ys\alpha _3+t y \alpha _4 & =&0 \\
-y \lambda_r(s,t) +s x \beta _1+xt\beta _2+ys\beta _3+t y \beta _4 & =&0,
\end{array}
\end{equation} where $u=\left(
\begin{array}{c}
x \\
y
\end{array}
\right)$, $v=\left(
\begin{array}{c}
s \\
t
\end{array}
\right).$
If $x=0$ then due to $y\neq 0$ the subspace $\mathbb{F}e_2$ is a two-sided ideal if and only if $\alpha _2=\alpha _3=\alpha _4=0$. If $x\neq 0$ then one can assume that $x=1$ and the system (\ref{SEri}) is equivalent to

\[\begin{array}{rrl}
l_{1A}(y)=y^2\alpha_4+y(\alpha_3-\beta_4)-\beta_3& =&0\\
l_{2A}(y)=y^2\alpha_2+y(\alpha_1-\beta_2)-\beta_1& =&0\\
r_{1A}(y)=y^2\alpha_4+y(\alpha_2-\beta_4)-\beta_2& =&0\\
r_{2A}(y)=y^2\alpha_3+y(\alpha_1-\beta_3)-\beta_1& =&0.
\end{array}\]

The subtraction the first equation from the third yields $y(\alpha_3-\alpha_2)-\beta_3+\beta_2=0$ and therefore, one easily gets  the following  result.

\begin{proposition}\label{pri}  The algebra $\mathbb{A}$ is simple, i.e., it  has no nontrivial two-sided ideal, only in the following cases:
\begin{itemize}	
	 \item[$1.$] $\beta_2\neq\beta_3$ and $\alpha_2=\alpha_3\neq 0$.
	 \item[$2.$]  $\beta_2\neq\beta_3$ and $\alpha_2=\alpha_3=0,\ \alpha_4\neq 0$.
	\item[$3.$] $\alpha_2\neq\alpha_3$ and at least one of $$(\beta_3-\beta_2)^2\alpha_2+(\beta_3-\beta_2)(\alpha_3-\alpha_2)(\alpha_1-\beta_2)-(\alpha_3-\alpha_2)^2\beta_1,$$ $$(\beta_3-\beta_2)^2\alpha_4+(\beta_3-\beta_2)(\alpha_3-\alpha_2)(\alpha_3-\beta_4)-(\alpha_3-\alpha_2)^2\beta_3$$ is not zero.

Moreover,
\begin{itemize}
\item[$\ast$] in $\alpha_2=\alpha_3=\alpha_4=0$ case $\mathbb{F}e_2$ is a two-sided ideal of $\mathbb{A};$
\item[$\ast$] in $\alpha_2\neq \alpha_3$ case $\mathbb{F}(e_1+\frac{\beta_3-\beta_2}{\alpha_3-\alpha_2}e_2)$ is a two-sided ideal if and only if

    $(\beta_3-\beta_2)^2\alpha_2+(\beta_3-\beta_2)(\alpha_3-\alpha_2)(\alpha_1-\beta_2)-(\alpha_3-\alpha_2)^2\beta_1$

     \hfill $=(\beta_3-\beta_2)^2\alpha_4+(\beta_3-\beta_2)(\alpha_3-\alpha_2)(\alpha_3-\beta_4)-(\alpha_3-\alpha_2)^2\beta_3=0$.

\item[$\ast$] If $\beta_3=\beta_2,\ \alpha_3=\alpha_2$ then $\mathbb{A}$ is a commutative algebra and every left (right) ideal is two-sided.
\end{itemize}
\end{itemize}
\end{proposition}

This proposition implies that for any nontrivial two-dimensional algebra only the following options occur:
\begin{itemize}
\item no two-sided ideal;
\item only one two-sided ideal;
\item two two-sided ideals.
\end{itemize}
Now we review the algebras from the list of isomorphism classes in Section 2 for the number of two-sided ideals.
\begin{itemize}
\item $Char(\mathbb{F})\neq 2,3$.
	
 Let us consider $A_1$. The case  $\alpha_2=\alpha_3=\alpha_4 =0$ is impossible. Note that for $A_1$ one has $\alpha_3-\alpha_2=\beta_3-\beta_2=1 $ and therefore, only in  $\alpha_4=-2\alpha_2-\alpha_1, \beta_1=\alpha_2+2\alpha_1$ case  $A_1$ has only one two-sided ideal $\mathbb{F}(e_1+e_2)$.

In $A_2$ case we have $\alpha_3-\alpha_2=0$ and therefore, if $\beta_2= 1-\alpha_1$ then in $A_2$ every left ideal is two-sided. In other cases $A_2$ has no two-sided ideal.

If $\beta_2= 1$ then every left ideal is two-sided of $A_3$. Otherwise,  it has no two-sided ideal.

If $\beta_2= 1-\alpha_1$ then every left ideal is two-sided for $A_4$. Otherwise, $A_4$ has only one two-sided ideal.

If $\alpha_1=\frac{2}{3}$ then every left ideal is two-sided for $A_5$. Otherwise, $A_5$ has only one two-sided ideal. As for $A_6$, $A_7$ and $A_{11}$ they do not have two-sided ideals. The algebras $A_8$, $A_9$, $A_{10}$ and $A_{12}$ have only one two-sided ideal.

\item $Char(\mathbb{F})=2.$

 In $A_{1,2}$ case  $\alpha_2=\alpha_3=\alpha_4 =0$ is impossible. Here we have $\alpha_3-\alpha_2=\beta_3-\beta_2=1 $ and therefore, only in  $\alpha_4=\alpha_1, \beta_1=\alpha_2$ case  $A_{1,2}$ has one two-sided ideal $\mathbb{F}(e_1+e_2)$.

In $A_{2,2}$ case one has $\alpha_3-\alpha_2=0$ and therefore, if $\beta_2= 1-\alpha_1$ then $A_{2,2}$ every left ideal is two-sided. In the other cases $A_{2,2}$ has no two-sided ideal.

If $ \alpha_1 = 1 , \ \beta_2=0$ then every left ideal of $A_{3,2}$ is two-sided. Otherwise it has no two-sided ideal.

If $\beta_2= 1-\alpha_1$ then every left ideal is two-sided for $A_{4,2}$. Otherwise has only one two-sided ideal.

If $\alpha_1=0$ then every left ideal is two-sided for $A_{5,2}$. Otherwise $A_{5,2}$ has only one two-sided ideal. The algebras $A_{6,2}$, $A_{7,2}$ $A_{11,2}$ and $A_{12,2}$ have no two-sided ideal. The algebras $A_{8,2}$, $A_{9,2}$ and $A_{10,2}$ have one two-sided ideal.

\item $Char(\mathbb{F})=3.$

$A_{1,3}$: In this case  $\alpha_2=\alpha_3=\alpha_4 =0$ is impossible and $\alpha_3-\alpha_2=\beta_3-\beta_2=1 $ and therefore only in  $\alpha_4=-\alpha_1-2\alpha_2, \beta_1=\alpha_2+2\alpha_2$ case  $A_{1,2}$ has only one two-sided ideal $\mathbb{F}(e_1+e_2)$.

In $A_{2,2}$ case $\alpha_3-\alpha_2=0$ and therefore if $\beta_2= 1-\alpha_1$ then $A_{2,3}$ every left ideal is two-sided. In other cases $A_{2,3}$ has no two-sided ideal.

If $\beta_2= 1$ then every left ideal is two-sided for $A_{3,3}$. Otherwise  it has no two-sided ideal.

If $\beta_2= 1-\alpha_1$ then every left ideal is two-sided for $A_{4,3}$. Otherwise has only one two-sided ideal.

 $A_{5,3}$ has only one two-sided ideal. $A_{6,3}$ has no two-sided ideal. $A_{7,3}$ has no two-sided ideal. $A_{8,3}$ has one two-sided ideal. $A_{9,3}$ has one two-sided ideal. $A_{10,3}$ has one two-sided ideal. $A_{11,2}$ has one two-sided ideal  and $A_{12,2}$ has one two-sided ideal.
\end{itemize}
 The results on two-sided ideals are given in Tables \ref{T5} and \ref{T6} (see APPENDIX).

 Note that the second (together with the first) columns of these tables  present classification of all simple two-dimensional algebras.
\section{Left quasiunits}
This section deals with left quasiunits of all two-dimensional algebras (for the motivation and the importance of the concept see \cite{Kantor} and \cite{K1} with occasion of the study of so-called \emph{Terminal algebras}).
\begin{definition} An element $\mathbf{e}_q$ of an algebra $\mathbb{A}$ is said to be a left quasiunit if
$$\mathbf{e}_q(\mathbf{uv})=(\mathbf{e}_q\mathbf{u})\mathbf{v}+\mathbf{u}(\mathbf{e}_q\mathbf{v})-\mathbf{uv}\ \mbox{ for all}\ \mathbf{u, v }\in \mathbb{A}.$$
\end{definition}
In terms of MSC of $\mathbb{A}$ it is nothing but to write
\[A(e_q\otimes A(u\otimes v))=A(A(e_q\otimes u)\otimes v)+A(u\otimes A(e_q\otimes v))-A(u\otimes v)\] and therefore assuming that
$e_q=\left(
\begin{array}{c}
x_0 \\
y_0
\end{array}
\right)$, $u=\left(
\begin{array}{c}
x\\
y
\end{array}
\right)$ and $v=\left(
\begin{array}{c}
s \\
t
\end{array}
\right)$ one gets the following system of equations\\
$s x (\alpha _1+\alpha _4 \beta _1 y_0-\alpha _2 \beta _3 y_0-\alpha _3 \beta _3 y_0-\alpha _3 \alpha _1 y_0-\alpha _3 \beta _1 x_0-\alpha _1^2 x_0 )+t x (\alpha _2-\alpha _2 \beta _4 y_0+\alpha _4 \beta _2 y_0-\alpha _4 \beta _3 y_0 -\alpha _1 \alpha _4 y_0-\alpha _4 \beta _1 x_0- \alpha _1 \alpha _2 x_0 )+s y (\alpha _3-\alpha _3 \beta _4 y_0-\alpha _1 \alpha _4 y_0+\alpha _2 \beta _3 x_0-\alpha _4 \beta _1 x_0-\alpha _3 \beta _2 x_0-\alpha _1 \alpha _2 x_0)$

\hfill$+t y (\alpha _4-\alpha _4 \beta _4 y_0-\alpha _4 \alpha _2 y_0+\alpha _2 \beta _4 x_0-2 \alpha _4 \beta _2 x_0-\alpha _2^2 x_0 -\alpha _3 \alpha _2 x_0+\alpha _1 \alpha _4 x_0)=0$,\\
$s x (\beta _1+\alpha _1 \beta _3 y_0-2 \alpha _3 \beta _1 y_0-\beta _3^2 y_0-\beta _2 \beta _3 y_0+\beta _1 \beta _4 y_0-\alpha _1 \beta _1 x_0-\beta _1 \beta _3 x_0)+t x (\beta _2-\alpha _4 \beta _1 y_0-\alpha _3 \beta _2 y_0+\alpha _2 \beta _3 y_0-\beta _3 \beta _4 y_0-\alpha _1 \beta _2 x_0-\beta _4 \beta _1 x_0)+s y (\beta _3-\alpha _4 \beta _1 y_0-\beta _3 \beta _4 y_0-\alpha _2 \beta _1 x_0+\alpha _3 \beta _1 x_0-\alpha _1 \beta _3 x_0-\beta _4 \beta _1 x_0 )$

\hfill $+t y (\beta _4-\alpha _4 \beta _2 y_0-\beta _4^2 y_0+\alpha _4 \beta _1 x_0-\alpha _2 \beta _2 x_0-\alpha _2 \beta _3 x_0-\beta _2 \beta _4 x_0)=0,$\\ which is equivalent to
\begin{equation}\label{SEQ}\begin{array}{rr}
-x_0 \alpha _1^2-y_0 \alpha _3 \alpha _1+\alpha _1-x_0 \alpha _3 \beta _1+y_0 \alpha _4 \beta _1-y_0 \alpha _2 \beta _3-y_0 \alpha _3 \beta _3 &=0,\\
-x_0 \alpha _1 \alpha _2-y_0 \beta _4 \alpha _2+\alpha _2-y_0 \alpha _1 \alpha _4-x_0 \alpha _4 \beta _1+y_0 \alpha _4 \beta _2-y_0 \alpha _4 \beta _3 &=0,\\
-x_0 \alpha _1 \alpha _2+x_0 \beta _3 \alpha _2+\alpha _3-y_0 \alpha _1 \alpha _4-x_0 \alpha _4 \beta _1-x_0 \alpha _3 \beta _2-y_0 \alpha _3 \beta _4 &=0,\\
-x_0 \alpha _2^2-x_0 \alpha _3 \alpha _2-y_0 \alpha _4 \alpha _2+x_0 \beta _4 \alpha _2+x_0 \alpha _1 \alpha _4+\alpha _4-2 x_0 \alpha _4 \beta _2-y_0 \alpha _4 \beta _4 &=0,\\
-y_0 \beta _3^2+y_0 \alpha _1 \beta _3-x_0 \beta _1 \beta _3-y_0 \beta _2 \beta _3-x_0 \alpha _1 \beta _1-2 y_0 \alpha _3 \beta _1+\beta _1+y_0 \beta _1 \beta _4&=0, \\
-y_0 \alpha _4 \beta _1-x_0 \beta _4 \beta _1-x_0 \alpha _1 \beta _2-y_0 \alpha _3 \beta _2+\beta _2+y_0 \alpha _2 \beta _3-y_0 \beta _3 \beta _4 &=0,\\
-x_0 \alpha _2 \beta _1+x_0 \alpha _3 \beta _1-y_0 \alpha _4 \beta _1-x_0 \beta _4 \beta _1-x_0 \alpha _1 \beta _3+\beta _3-y_0 \beta _3 \beta _4 &=0,\\
-y_0 \beta _4^2-x_0 \beta _2 \beta _4+\beta _4+x_0 \alpha _4 \beta _1-x_0 \alpha _2 \beta _2-y_0 \alpha _4 \beta _2-x_0 \alpha _2 \beta _3&=0. \\
\end{array}\end{equation}
\subsection{$Char(\mathbb{F})\neq 2,3$}
For $A_{1}$ case we have the system of equations
\[\begin{array}{rr} \label{81}
-x_0 \alpha _1^2+y_0 \alpha _2 \alpha _1+\alpha _1-y_0-2 y_0 \alpha _2-x_0 \beta _1-x_0 \alpha _2 \beta _1+y_0 \alpha _4 \beta _1&=0, \\
y_0 \alpha _2^2-x_0 \alpha _1 \alpha _2+\alpha _2-y_0 \alpha _4-y_0 \alpha _1 \alpha _4-x_0 \alpha _4 \beta _1 &=0,\\
y_0 \alpha _2^2+y_0 \alpha _2+x_0 \alpha _2-x_0 \alpha _1 \alpha _2+\alpha _2+x_0 \alpha _1-y_0 \alpha _1 \alpha _4-x_0 \alpha _4 \beta _1+1 &=0,\\
-3 x_0 \alpha _2^2-x_0 \alpha _2+3 x_0 \alpha _1 \alpha _4+\alpha _4 &=0,\\
-3 y_0 \alpha _1^2+4 y_0 \alpha _1-y_0-2 y_0 \beta _1-x_0 \beta _1-3 y_0 \alpha _2 \beta _1+\beta _1 &=0,\\
x_0 \alpha _1^2+y_0 \alpha _1-y_0 \alpha _2 \alpha _1-\alpha _1+2 y_0 \alpha _2+x_0 \alpha _2 \beta _1-y_0 \alpha _4 \beta _1 &=0,\\
x_0 \alpha _1^2-x_0 \alpha _1-y_0 \alpha _2 \alpha _1-\alpha _1+y_0 \alpha _2+x_0 \beta _1+x_0 \alpha _2 \beta _1-y_0 \alpha _4 \beta _1+1 &=0,\\
-y_0 \alpha _2^2-x_0 \alpha _2+x_0 \alpha _1 \alpha _2-\alpha _2+y_0 \alpha _1 \alpha _4+x_0 \alpha _4 \beta _1 &=0.
\end{array}. \]
Combining the first, third, seventh and eighth equations of the system of equations above we get $y_0=2.$ The third and eighth equations imply $x_0\alpha_1=-1-2\alpha_2$, the second and eighth equations imply $x_0 \alpha _2+2 \alpha _4=0$, the  first and sixth equations imply $x_0\beta_1=2(\alpha_1-1).$

Therefore, if $\beta_1\neq 0$ then we have $$x_0=\frac{2(\alpha_1-1)}{\beta_1}, \
 \alpha_2=\frac{x_0\alpha_1+1}{-2}=-\frac{2\alpha_1(\alpha_1-1)+\beta_1}{2\beta_1}, \ \alpha_4=-\frac{x_0\alpha_2}{2}=-\frac{2\alpha_2(\alpha_1-1)}{2\beta_1}$$ and in this case $e_q=\left(
		\begin{array}{l}
		x_0 \\
		y_0
		\end{array}
		\right)=\left(
		\begin{array}{c}
		\frac{2(\alpha_1-1)}{\beta_1} \\
		2
		\end{array}
		\right)$ satisfies (8.1), i.e., $e_q$ is a left quasiunit.
		
If $\beta_1=0$ then for (8.1) to have a solution $\alpha_1$ must be $1$. In this case $x_0=-1-2\alpha_2$, $\alpha_4=\frac{\alpha_2(1+2\alpha_2)}{2}$. So only
$A_1(\alpha_1,-\frac{2\alpha_1(\alpha_1-1)+\beta_1}{2\beta_1},  -\frac{2\alpha_2(\alpha_1-1)}{2\beta_1}, \beta_1)$, where $\beta_1 \neq 0$, and $A_1(1,\alpha_2, \frac{\alpha_2(1+2\alpha_2)}{2},0)$ have left quasiunits $\frac{2(\alpha_1-1)}{\beta_1}e_1+2e_2$ and $ -(1+2\alpha_2)e_1+2e_2$, respectively.

In $A_{2}$ case the system (8.1) is written as follows
\[\begin{array}{rr}
-x_0 \alpha _1^2+\alpha _1+y_0 \beta _1  &=0, \\
\beta _2 y_0-y_0-x_0 \beta _1  &=0, \\
-y_0 \alpha _1-x_0 \beta _1 &=0,  \\
x_0 \alpha _1-2 x_0 \beta _2+1  &=0, \\
-2 y_0 \alpha _1^2+3 y_0 \alpha _1+y_0 \beta _2 \alpha _1-y_0-x_0 \beta _1+\beta _1-y_0 \beta _2  &=0,  \\
-y_0 \beta _1-x_0 \alpha _1 \beta _2+\beta _2  &=0, \\
x_0 \alpha _1^2-x_0 \alpha _1-\alpha _1-y_0 \beta _1+1  &=0, \\
x_0 \beta _1-y_0 \beta _2  &=0.
\end{array} \]
The second and eighth equations imply $y_0=0,$ the first and seventh equations imply $x_0=\frac{1}{\alpha _1}$, where $\alpha _1 \neq 0$, therefore, $\beta_1 =0$ and  $\beta_2=\alpha_1.$ This means that in this case only $A_2(\alpha_1,0,\alpha_1)$, where $\alpha_1\neq 0$, has left quasiunit, which is $\frac{1}{\alpha_1}e_1$.

In $A_{3}$ case the corresponding system is given as follows
\[\begin{array}{rr}
-2 y_0-x_0 \beta _1 &=0,  \\
y_0+1  &=0, \\
y_0+x_0-x_0 \beta _2+1  &=0, \\
-3 x_0  &=0, \\
-3 \beta _1 y_0-\beta _2 y_0-y_0-x_0 \beta _1+\beta _1  &=0, \\
-\beta _2 y_0+2 y_0+x_0 \beta _1+\beta _2  &=0, \\
y_0+x_0 \beta _1+1  &=0, \\
-y_0-x_0-1  &=0,\\
\end{array}\]
 which is inconsistent.

For $A_{4}$ case we have the system of equations
\[\begin{array}{rr}
\alpha _1-x_0 \alpha _1^2  &=0, \\
-2 y_0 \alpha _1^2+3 y_0 \alpha _1+y_0 \beta _2 \alpha _1-y_0-y_0 \beta _2 &=0,  \\
\beta _2-x_0 \alpha _1 \beta _2  &=0, \\
x_0 \alpha _1^2-x_0 \alpha _1-\alpha _1+1  &=0.\\
\end{array}\]
The system above implies that in this case only the following algebras have left quasiunits
 \begin{itemize}
 \item $A_{4}(1,\beta_2)$, the left quasiunits are $e_1+te_2$, where $t\in \mathbb{F}$,
 \item $A_{4}(\alpha_1,\beta_2)$, where $\alpha_1\neq 0,1,\beta_2\neq 2\alpha_1-1$, the left quasiunit is $\frac{1}{\alpha_1}e_1$,
 \item $A_{4}(\alpha_1,\beta_2)$, where $\alpha_1\neq 0,1,\beta_2= 2\alpha_1-1$, the left quasiunits are $\frac{1}{\alpha_1}e_1+te_2$, where $t\in \mathbb{F}$.
\end{itemize}
For $A_{5}$ we have the following system of equations
\[\begin{array}{rr}
\alpha _1-x_0 \alpha _1^2 &=0,  \\
1-x_0  &=0, \\
-2 x_0 \alpha _1^2+x_0 \alpha _1+2 \alpha _1-1  &=0, \\
x_0 \alpha _1^2-x_0 \alpha _1-\alpha _1+1 &=0.
\end{array}\]
It is easy to see that the only
$A_5(1)$ has left quasiunits $e_1+te_2$, $t\in \mathbb{F}$.

For $A_{6}$ the system is given as follows
\[\begin{array}{rr}
-x_0 \alpha _1^2+\alpha _1+y_0 \beta _1  &=0, \\
-\alpha _1 y_0+y_0-x_0 \beta _1  &=0, \\
-y_0 \alpha _1-x_0 \beta _1  &=0, \\
3 \alpha _1 x_0-2 x_0+1  &=0, \\
-3 y_0 \alpha _1^2+y_0 \alpha _1+\beta _1 &=0,  \\
x_0 \alpha _1^2-x_0 \alpha _1-\alpha _1-y_0 \beta _1+1 &=0,
\end{array}\]
and again it is easy to see that only $A_6(\frac{1}{2}, 0)$ has left quasiunit given by $2e_1$.

For $A_{7}$ the system consists of the equations
$
-x_0 \beta _1 =0, \
y_0+1  =0, \
y_0-x_0+1 =0,  \
-3 x_0  =0, \
\beta _1-3 y_0 \beta _1 =0, \
-y_0+x_0 \beta _1+1  =0 $
and it is inconsistent.

For $A_{8}$ we have the system consisting of the equations
$
\alpha _1-x_0 \alpha _1^2=0,\
y_0 \alpha _1-3 y_0 \alpha _1^2 =0, \
x_0 \alpha _1^2-x_0 \alpha _1-\alpha _1+1 =0$
and only $A_8(\alpha_1)$, where $\alpha_1 \neq 0, \frac{1}{3}$ and $A_8(\frac{1}{3})$ have left quasiunits given by $\frac{1}{\alpha_1}e_1$ and $3e_1+te_2$, $t\in \mathbb{F}$, respectively.

For $A_{9}$ the fifth equation of (\ref{SEQ}) is a contradiction.

For $A_{10}$ the system consists of the equations
$y_0+1=0, \
-3 x_0  =0$
and therefore $A_{10}$ has a left quasiunit $-e_2$.

For $A_{11}$ the system consists of the equations
$-x_0=0,\
y_0+1 =0, \
1-3 y_0 =0, $
so $A_{11}$ has no left quasiunit.

For $A_{12}$ the fifth equation of (\ref{SEQ}) is a contradiction.

\subsection{$Char(\mathbb{F})=2. $}
For $A_{1,2}$ we have the system :
\begin{equation}
\begin{array}{ll}
 -x_0 \alpha _1^2+y_0 \alpha _2 \alpha _1+\alpha _1-y_0-x_0 \beta _1-x_0 \alpha _2 \beta _1+y_0 \alpha _4 \beta _1&=0, \\
 y_0 \alpha _2^2-x_0 \alpha _1 \alpha _2+\alpha _2-y_0 \alpha _4-y_0 \alpha _1 \alpha _4-x_0 \alpha _4 \beta _1 &=0,\\
 y_0 \alpha _2^2+y_0 \alpha _2+x_0 \alpha _2-x_0 \alpha _1 \alpha _2+\alpha _2+x_0 \alpha _1-y_0 \alpha _1 \alpha _4-x_0 \alpha _4 \beta _1+1 &=0,\\
 - x_0 \alpha _2^2-x_0 \alpha _2+ x_0 \alpha _1 \alpha _4+\alpha _4 &=0,\\
 -y_0 \alpha _1^2-x_0 \beta _1- y_0 \alpha _2 \beta _1+\beta _1 &=0,\\
 x_0 \alpha _1^2+y_0 \alpha _1-y_0 \alpha _2 \alpha _1-\alpha _1+x_0 \alpha _2 \beta _1-y_0 \alpha _4 \beta _1 &=0,\\
 x_0 \alpha _1^2-x_0 \alpha _1-y_0 \alpha _2 \alpha _1-\alpha _1+y_0 \alpha _2+x_0 \beta _1+x_0 \alpha _2 \beta _1-y_0 \alpha _4 \beta _1+1 &=0,\\
 -y_0 \alpha _2^2-x_0 \alpha _2+x_0 \alpha _1 \alpha _2-\alpha _2+y_0 \alpha _1 \alpha _4+x_0 \alpha _4 \beta _1 &=0\\
\end{array} \end{equation}
 From the first, seventh ,third and eighth equations we get $y_0=0.$ Now, applying third and eighth equations we have $x_0\alpha_1=-1\Longrightarrow x_0=\frac{1}{\alpha_1}$ where $\alpha_1 \neq 0$ then one has $\beta_1=0.$ Using second and eighth equations we get $x_0 \alpha _2=0 \Longrightarrow \alpha _2=0.$ then
 \[\left( \begin{array}{c}
      \frac{1}{\alpha_1} \\
      0
    \end{array}\right) \mbox{ where } \alpha_1 \neq 0,\ \alpha_2=0,\ \beta_1=0.\]

For $A_{2,2}$ we have the system :
\begin{equation}
\begin{array}{ll}
 -x_0 \alpha _1^2+\alpha _1+y_0 \beta _1  &=0, \\
 \beta _2 y_0-y_0-x_0 \beta _1  &=0, \\
 -y_0 \alpha _1-x_0 \beta _1 &=0,  \\
 x_0 \alpha _1+1  &=0, \\
  y_0 \alpha _1+y_0 \beta _2 \alpha _1-y_0-x_0 \beta _1+\beta _1-y_0 \beta _2  &=0,  \\
 -y_0 \beta _1-x_0 \alpha _1 \beta _2+\beta _2  &=0, \\
 x_0 \alpha _1^2-x_0 \alpha _1-\alpha _1-y_0 \beta _1+1  &=0, \\
 x_0 \beta _1-y_0 \beta _2  &=0\\
\end{array} \end{equation}
From the second and eighth equations, the first and seventh equations, respectively  we get $y_0=0,$ $x_0=\frac{1}{\alpha _1}$ where $\alpha _1 \neq 0$ then $\beta_1 =0$ therefore
\[\left( \begin{array}{c}
      \frac{1}{\alpha_1} \\
      0
    \end{array}\right) \mbox{ where } \alpha_1\neq 0,\ \beta_1=0.\]

For $A_{3,2}(\alpha_1,\beta_2)$ we have the system :
\begin{equation}
\begin{array}{ll}
 \alpha _1+\alpha _1 y_0-\alpha _1^2 x_0 &=0,\\
 -y_0-\alpha _1 x_0+1 &=0,\\
 -y_0-\beta _2 x_0+x_0+1&=0,\\
 -x_0 &=0,\\
 \alpha _1 \beta _2 y_0+\alpha _1 y_0-\beta _2 y_0-y_0 &=0,\\
 \beta _2-\beta _2 y_0-\alpha _1 \beta _2 x_0 &=0,\\
 -\alpha _1+\alpha _1 y_0-y_0+\alpha _1^2 x_0-\alpha _1 x_0+1 &=0,\\
 -y_0+\alpha _1 x_0-x_0+1  &=0.\\
\end{array}\end{equation}
We have $x_0=0$ then  $y_0=1.$ Due to the fifth equation we get $\alpha_1=1$ or $\beta_2=1,$ then $A_{3,2}(1,\beta_2)$ and $A_{3,2}(\alpha_1,1)$ have a left quasiunit given by $e_2.$

For $A_{4,2}$ we have the system :
\begin{equation}
\begin{array}{ll}
 \alpha _1-x_0 \alpha _1^2  &=0, \\
  y_0 \alpha _1+y_0 \beta _2 \alpha _1-y_0-y_0 \beta _2 &=0,  \\
 \beta _2-x_0 \alpha _1 \beta _2  &=0, \\
 x_0 \alpha _1^2-x_0 \alpha _1-\alpha _1+1  &=0\\
\end{array}\end{equation}
Applying the first and fourth equations we get $x_0=\frac{1}{\alpha_1}$ the second equation implies
$y_0( \alpha _1+ \beta _2 \alpha _1-1-\beta _2 )=0 \Longrightarrow y_0( \alpha _1(1+ \beta _2 )-(1+\beta _2 ))=y_0(\alpha_1-1)(\beta _2 +1)=0$ therefore
\[\left( \begin{array}{c}
      \frac{1}{\alpha_1} \\
      0
    \end{array}\right) \mbox{ where } \alpha_1 \neq 0,1, \ \beta _2 \neq 1,\]
    \[\left( \begin{array}{c}
      1 \\
      y_0
    \end{array}\right) \mbox{ where } \alpha_1 =1,\]
\[\left( \begin{array}{c}
      \frac{1}{\alpha_1} \\
      y_0
    \end{array}\right) \mbox{ where } \alpha_1 \neq 0,\ 1,\beta _2 =1.\]
For $A_{5,2}$ we have the system :
\begin{equation}
\begin{array}{ll}
 \alpha _1-x_0 \alpha _1^2 &=0,  \\
 1-x_0  &=0, \\
 x_0 \alpha _1-1  &=0, \\
 x_0 \alpha _1^2-x_0 \alpha _1-\alpha _1+1 &=0 \\
\end{array}\end{equation}
Then
\[\left( \begin{array}{c}
      1 \\
      y_0
    \end{array}\right) \mbox{ where } \alpha_1 =1.\]

For $A_{6,2}$ we have the system :
\begin{equation}
\begin{array}{ll}
 -x_0 \alpha _1^2+\alpha _1+y_0 \beta _1  &=0, \\
 -\alpha _1 y_0+y_0-x_0 \beta _1  &=0, \\
 -y_0 \alpha _1-x_0 \beta _1  &=0, \\
  \alpha _1 x_0+1  &=0, \\
 - y_0 \alpha _1^2+y_0 \alpha _1+\beta _1 &=0,  \\
 x_0 \alpha _1^2-x_0 \alpha _1-\alpha _1-y_0 \beta _1+1 &=0,  \\
\end{array}\end{equation}
Applying the second and third equations we have $y_0=0$ so using the fifth equation we get $\beta_1=0$ and from the first and sixth equations we have $x_0=\frac{1}{\alpha_1}$ where $\alpha_1\neq 0$
Hence
\[\left( \begin{array}{c}
      \frac{1}{\alpha_1} \\
      0
    \end{array}\right) \mbox{ where } \beta_1=0,\ \alpha_1 \neq 0.\]

For $A_{7,2}$
\begin{equation}
\begin{array}{ll}
 -x_0 \alpha _1^2+y_0 \alpha _1+\alpha _1 &=0,\\
 y_0-x_0 \alpha _1+1  &=0,\\
 y_0-x_0-x_0 \alpha _1+1  &=0,\\
 - x_0 &=0, \\
 y_0 \alpha _1- y_0 \alpha _1^2  &=0,\\
 x_0 \alpha _1^2-y_0 \alpha _1-x_0 \alpha _1-\alpha _1-y_0+1  &=0,\\
\end{array}
\end{equation}
\[\left( \begin{array}{c}
      0 \\
      1
    \end{array}\right) \mbox{ where } \alpha_1 =0,\ 1.\]
For $A_{8,2}$ we have the system :
\begin{equation}
\begin{array}{ll}
 \alpha _1-x_0 \alpha _1^2 &=0,  \\
 y_0 \alpha _1- y_0 \alpha _1^2  &=0, \\
 x_0 \alpha _1^2-x_0 \alpha _1-\alpha _1+1  &=0, \\
\end{array}\end{equation}
From the first and third equations we have $x_0=\frac{1}{\alpha_1}$ where $\alpha_1 \neq 0$ by the second equation we get $y_0\alpha_1(1-\alpha_1)=0$ then
\[\left( \begin{array}{c}
      \frac{1}{\alpha_1} \\
      0
    \end{array}\right) \mbox{ where }\alpha_1 \neq 0, 1.\]
    \[\left( \begin{array}{c}
      1 \\
      y_0
    \end{array}\right) \mbox{ where }\alpha_1 = 1.\]

For $A_{9,2}$ the fifth equation in the system (\ref{SEQ}) not hold so there is no quasiunit in this case.\\

For $A_{10,2}$ we have the system :
\begin{equation}
\begin{array}{ll}
 y_0+1  &=0, \\
 - x_0  &=0 \\
\end{array}\end{equation}
\[\left( \begin{array}{c}
      0 \\
      1
    \end{array}\right).\]

For $A_{11,2},$ $A_{12,2}$ the fifth equation in the system (\ref{SEQ}) not hold so there is no quasiunit in these cases.

\subsection{$Char(\mathbb{F})=3. $}
For $A_{1,3}$ we have the system :
\begin{equation}
\begin{array}{ll}
 -x_0 \alpha _1^2+y_0 \alpha _2 \alpha _1+\alpha _1-y_0-2 y_0 \alpha _2-x_0 \beta _1-x_0 \alpha _2 \beta _1+y_0 \alpha _4 \beta _1&=0, \\
 y_0 \alpha _2^2-x_0 \alpha _1 \alpha _2+\alpha _2-y_0 \alpha _4-y_0 \alpha _1 \alpha _4-x_0 \alpha _4 \beta _1 &=0,\\
 y_0 \alpha _2^2+y_0 \alpha _2+x_0 \alpha _2-x_0 \alpha _1 \alpha _2+\alpha _2+x_0 \alpha _1-y_0 \alpha _1 \alpha _4-x_0 \alpha _4 \beta _1+1 &=0,\\
-x_0 \alpha _2+\alpha _4 &=0,\\
y_0 \alpha _1-y_0-2 y_0 \beta _1-x_0 \beta _1+\beta _1 &=0,\\
 x_0 \alpha _1^2+y_0 \alpha _1-y_0 \alpha _2 \alpha _1-\alpha _1+2 y_0 \alpha _2+x_0 \alpha _2 \beta _1-y_0 \alpha _4 \beta _1 &=0,\\
 x_0 \alpha _1^2-x_0 \alpha _1-y_0 \alpha _2 \alpha _1-\alpha _1+y_0 \alpha _2+x_0 \beta _1+x_0 \alpha _2 \beta _1-y_0 \alpha _4 \beta _1+1 &=0,\\
 -y_0 \alpha _2^2-x_0 \alpha _2+x_0 \alpha _1 \alpha _2-\alpha _2+y_0 \alpha _1 \alpha _4+x_0 \alpha _4 \beta _1 &=0\\
\end{array} \end{equation}
 From the first, seventh ,third and eighth equations we get $y_0=2.$ Now, applying third and eighth equations; the first and sixth equations respectively we have $x_0\alpha_1=-1-2\alpha_2$ and $x_0\beta_1=2(\alpha_1-1).$ Using the fourth equation we get $x_0 \alpha _2- \alpha _4=0.$\\
 \underline{If $\alpha_2 =0$} then $\alpha_4=0,$ $x_0=-\frac{1}{\alpha_1}$ where $\alpha_1 \neq 0$ ($\alpha_1$ cannot be zero in this case) and $\beta_1=-2\alpha_1(\alpha_1-1)$ then
 \[\left( \begin{array}{c}
      -\frac{1}{\alpha_1} \\
      2
    \end{array}\right) \mbox{ where } \alpha_1 \neq 0,\ \alpha_2=\alpha_4=0,\ \beta_1=-2\alpha_1(\alpha_1-1)\]
 \underline{If $\alpha_2 \neq 0$} then $x_0=\frac{\alpha_4}{\alpha_2}$ so\\
  $x_0\beta_1=2(\alpha_1-1) \Longrightarrow \alpha_1=\frac{\alpha_2- \alpha_4\beta_1}{\alpha_2}$ and\\
  $x_0\alpha_1=-1-2\alpha_2 \Longrightarrow 2\alpha_2^3 +\alpha_2^2 +\alpha_2\alpha_4 -\alpha_4^2 \beta_1=0.$ Therefore
  \[\left( \begin{array}{c}
      \frac{\alpha_4}{\alpha_2} \\
      2
    \end{array}\right) \mbox{ where } \alpha_1 =\frac{\alpha_2- \alpha_4\beta_1}{\alpha_2},\ 2\alpha_2^3 +\alpha_2^2 +\alpha_2\alpha_4 -\alpha_4^2 \beta_1=0\]

For $A_{2,3}$ we have the system :
\begin{equation}
\begin{array}{ll}
 -x_0 \alpha _1^2+\alpha _1+y_0 \beta _1  &=0, \\
 \beta _2 y_0-y_0-x_0 \beta _1  &=0, \\
 -y_0 \alpha _1-x_0 \beta _1 &=0,  \\
 x_0 \alpha _1-2 x_0 \beta _2+1  &=0, \\
 -2 y_0 \alpha _1^2+y_0 \beta _2 \alpha _1-y_0-x_0 \beta _1+\beta _1-y_0 \beta _2  &=0,  \\
 -y_0 \beta _1-x_0 \alpha _1 \beta _2+\beta _2  &=0, \\
 x_0 \alpha _1^2-x_0 \alpha _1-\alpha _1-y_0 \beta _1+1  &=0, \\
 x_0 \beta _1-y_0 \beta _2  &=0\\
\end{array} \end{equation}
From the second and eighth equations, the first and seventh equations, respectively  we get $y_0=0,$ $x_0=\frac{1}{\alpha _1}$ where $\alpha _1 \neq 0$ then $\beta_1 =0,\ \beta_2=\alpha_1$ therefore
\[\left( \begin{array}{c}
      \frac{1}{\alpha_1} \\
      0
    \end{array}\right) \mbox{ where } \alpha_1\neq 0,\ \beta_1=0,\ \beta_2= \alpha_1.\]

For $A_{3,3}$ we have the system :
\begin{equation}
\begin{array}{ll}
 -2 y_0-x_0 \beta _1 &=0,  \\
 y_0+1  &=0, \\
 y_0+x_0-x_0 \beta _2+1  &=0, \\
 -\beta _2 y_0-y_0-x_0 \beta _1+\beta _1  &=0, \\
 -\beta _2 y_0+2 y_0+x_0 \beta _1+\beta _2  &=0, \\
 y_0+x_0 \beta _1+1  &=0, \\
 -y_0-x_0-1  &=0\\
\end{array}\end{equation}
Using the second, seventh and sixth equations we get a contradiction.\\

For $A_{4,3}$ we have the system :
\begin{equation}
\begin{array}{ll}
 \alpha _1-x_0 \alpha _1^2  &=0, \\
 -2 y_0 \alpha _1^2+y_0 \beta _2 \alpha _1-y_0-y_0 \beta _2 &=0,  \\
 \beta _2-x_0 \alpha _1 \beta _2  &=0, \\
 x_0 \alpha _1^2-x_0 \alpha _1-\alpha _1+1  &=0\\
\end{array}\end{equation}
Applying the first and fourth equations we get $x_0=\frac{1}{\alpha_1}$ where $\alpha_1 \neq 0,$ the second equation implies
$y_0(   \alpha _1^2+ \beta _2 \alpha _1-1-\beta _2 )=0 \Longrightarrow y_0( (\alpha _1^2 -1)+ \beta _2(\alpha_1-1))=y_0(\alpha_1-1)(\alpha_1+1+ \beta _2 )=0$ therefore
\[\left( \begin{array}{c}
      \frac{1}{\alpha_1} \\
      0
    \end{array}\right) \mbox{ where } \alpha_1 \neq 0, 1, \ \beta _2 \neq - \alpha _1-1,\]
     \[\left( \begin{array}{c}
     1 \\
      y_0
    \end{array}\right) \mbox{ where } \alpha_1 =1,\]

    \[\left( \begin{array}{c}
      \frac{1}{\alpha_1} \\
      y_0
    \end{array}\right) \mbox{ where } \alpha_1 \neq 0,1, \ \beta _2 = 0- \alpha _1-1.\]

For $A_{5,3}$ we have the system :
\begin{equation}
\begin{array}{ll}
 \alpha _1-x_0 \alpha _1^2 &=0,\\
 -x_0+1&=0, \\
 x_0 \alpha _1^2+x_0 \alpha _1-\alpha _1-1 &=0,\\
 x_0 \alpha _1^2-x_0 \alpha _1-\alpha _1+1&=0. \\
\end{array}
\end{equation}
Then
\[\left( \begin{array}{c}
      1 \\
      y_0
    \end{array}\right) \mbox{ where } \alpha_1 =1.\]

For $A_{6,3}$ we have the system :
\begin{equation}
\begin{array}{ll}
 -x_0 \alpha _1^2+\alpha _1+y_0 \beta _1  &=0, \\
 -\alpha _1 y_0+y_0-x_0 \beta _1  &=0, \\
 -y_0 \alpha _1-x_0 \beta _
 1  &=0, \\
-2 x_0+1  &=0, \\
y_0 \alpha _1+\beta _1 &=0,  \\
 x_0 \alpha _1^2-x_0 \alpha _1-\alpha _1-y_0 \beta _1+1 &=0,  \\
\end{array}\end{equation}
Applying the second and third equations we have $y_0=0$ so using the fifth equation we get $\beta_1=0$ and from the fourth equation we have $x_0=-1$ then $\alpha_1=-1$
Hence
\[\left( \begin{array}{c}
      -1 \\
      0
    \end{array}\right) \mbox{ where } \beta_1=0,\ \alpha_1 =-1.\]

For $A_{7,3}$ we have the system :
\begin{equation}
\begin{array}{ll}
 -x_0 \beta _1  &=0, \\
 y_0+1  &=0, \\
 y_0-x_0+1 &=0,  \\
 \beta _1&=0,  \\
 -y_0+x_0 \beta _1+1  &=0 \\
\end{array}\end{equation}
Using the first, second and fifth equations we get a contradiction.\\

For $A_{8,3}$ we have the system :
\begin{equation}
\begin{array}{ll}
 \alpha _1-x_0 \alpha _1^2 &=0,  \\
 y_0 \alpha _1  &=0, \\
 x_0 \alpha _1^2-x_0 \alpha _1-\alpha _1+1  &=0, \\
\end{array}\end{equation}
From the first and third equations we have $x_0=\frac{1}{\alpha_1}$ where $\alpha_1 \neq 0$ by the second equation we get $y_0=0$ then
\[\left( \begin{array}{c}
      \frac{1}{\alpha_1} \\
      0
    \end{array}\right) \mbox{ where }\alpha_1 \neq 0.\]

For $A_{9,3}$ the fifth equation in the system (\ref{SEQ}) not hold so there is no quasiunit in this case.\\

For $A_{10,3}$ we have the system :
\begin{equation}
\begin{array}{ll}
 y_0+1  &=0, \\
\end{array}\end{equation}
\[\left( \begin{array}{c}
      x_0 \\
      -1
    \end{array}\right).\]

For $A_{11,3},$ $A_{12,3}$ the fifth equation in the system (\ref{SEQ}) not hold so there is no quasiunit in these cases.\\

We summarize the results in Table \ref{T7} (See Section APPENDIX).
\begin{remark} In APPENDIX the following must be taken into account:
\begin{itemize}
  \item In Table \ref{T1} only the algebra $A_4(1,1)$ has infinitely many left ideals;
  \item In Table \ref{T2} only the algebras $A_{4,2}(1,1)$ and $A_{4,3}(1,1)$ have infinitely many left ideals;
  \item In Table \ref{T3} only $A_4(\frac{1}{2},1)$ has infinitely many right ideals;
  \item In Table \ref{T4} only $A_{8,2}(1)$ and $A_{4,3}(-1,0)$ have infinitely many right ideals.
\end{itemize}\end{remark}

%


\newpage
\section{APPENDIX}
	\begin{footnotesize}
		\begin{longtable}{ |c|c| c |c| }
			 \caption{Nontrivial left ideals, $Char(\mathbb{F})\neq 2,3.$ $ P=\left(\alpha _4 \beta _1-\alpha _2 \left(1-\alpha _1\right)\right)^2+\left(2 \alpha _1 \alpha _4-\alpha _2 \left(2 \alpha _2+1\right)\right) \left(\left(2 \alpha _2+1\right) \beta _1-2 \alpha _1 \left(1-\alpha _1\right)\right).$}\label{T1}\\
				\hline
				\textbf{Alg.} & $\begin{array}{cc}&\textbf{No}\\ &\textbf{left ideal}\end{array}$ &  $\begin{array}{cc}&\textbf{One}\\ &\textbf{left ideal}\end{array}$ &  $\begin{array}{cc}&\textbf{Two}\\ &\textbf{left ideals}\end{array}$  \\
				\hline
				\multirow{16}{*}{$A_{1}$}&$\begin{array}{l}\alpha_4=0,\ \alpha_2\neq -1/2,\ 0\\
				\beta_1\neq\frac{(1-\alpha_1)(2\alpha_1+\alpha_2+3\alpha_1\alpha_2)}{(2\alpha_2+1)^2}\end{array}$& $\begin{array}{l}\alpha_4=0,\ \alpha_2\neq -1/2,\ 0\\
				\beta_1=\frac{(1-\alpha_1)(2\alpha_1+\alpha_2+3\alpha_1\alpha_2)}{(2\alpha_2+1)^2}\end{array}$ & \multirow{1}{*}{$\begin{array}{l}\alpha_4\neq 0,\\ \alpha_1= 1+\frac{(2\alpha_2+1)^2}{4\alpha_4},\\ \beta_1=\frac{\alpha_2(2\alpha_2+1)}{2\alpha_4},\\
					\alpha_2=-\frac{4\alpha_4+1}{2}\end{array}$}
				\\
				\cline{2-3}
				& $\begin{array}{l}\alpha_4= 0,\
				\alpha_1\neq 1,\\
				\alpha_2=-1/2 \end{array}$ &$\begin{array}{l}\alpha_4\neq 0,\
				\alpha_1=1+ \frac{(2\alpha_2+1)^2}{4\alpha_4},\\ \beta_1=\frac{(2\alpha_2+1)(2\alpha_2^2+\alpha_2-4\alpha_4\alpha_1)}{4\alpha_4^2},	
				\end{array}$	&   \\
				\cline{2-3}
				& $ \begin{array}{l} P\neq 0,\ \alpha_4\neq 0,\\
				\alpha_1\neq 1+\frac{(2\alpha_2+1)^2}{4\alpha_4}, \frac{\alpha_2(2\alpha_2+1)}{2\alpha_4},
				\end{array}$& $\begin{array}{l} P= 0,\ \alpha_4\neq 0,\\
				\alpha_1\neq 1+ \frac{(2\alpha_2+1)^2}{4\alpha_4}, \frac{\alpha_2(2\alpha_2+1)}{2\alpha_4},
				\end{array}$ &   \\
				\cline{2-4}
				&$\begin{array}{l}\alpha_4\neq 0,\
				\alpha_1=1+ \frac{(2\alpha_2+1)^2}{4\alpha_4},\\ \beta_1\neq\frac{(2\alpha_2+1)(2\alpha_2^2+\alpha_2-4\alpha_4\alpha_1)}{4\alpha_4^2},	
				\end{array}$ &  $
				\begin{array}{l}\alpha_1= 1,\
				\alpha_2=-1/2,\\ \alpha_4=0,\
				\beta_1=2
				\end{array} $ &$\begin{array}{l}\alpha_4= 0,\ \alpha_1=1,\\
				\alpha_2=-1/2,\
				\beta_1\neq 2
				\end{array}$\\
				\cline{2-4}
				&$\begin{array}{l}\alpha_4\neq 0,\
				\alpha_1\neq 1+\frac{(2\alpha_2+1)^2}{4\alpha_4},\\
				\alpha_1=\frac{\alpha_2(2\alpha_2+1)}{2\alpha_4},\
				\beta_1\neq \frac{\alpha_2(1-\alpha_1)}{\alpha_4},
				\end{array}$ & $\begin{array}{l}\alpha_4= \alpha_2= 0\\
				\beta_1\neq 2\alpha_1(1-\alpha_1)\end{array}$ &$\begin{array}{l}\alpha_4=\alpha_2= 0\\
				\beta_1=2\alpha_1(1-\alpha_1)\end{array}$
				\\
				\hline
				\multirow{6}{*}{$A_{2}$}& $\begin{array}{l} \beta_1\neq 0,\ \alpha_1=1\end{array}$&	$\begin{array}{l} \alpha_1=1,\ \beta_1=0\end{array}$&
				\\
				\cline{2-3}
				&$\begin{array}{l} \alpha_1 \neq 1,\ \beta_2\neq \alpha_1,\\ \beta_1^2 +\alpha_1-1\neq 0\end{array}$& \multirow{2}{*}{$\begin{array}{l} \alpha_1 \neq 1,\ \beta_2\neq \alpha_1,\\ \beta_1^2 +\alpha_1-1= 0\end{array}$ } & $\begin{array}{l} \alpha_1 \neq 1,\\ \beta_2= \alpha_1,\\ \beta_1=0\end{array}$\\
				\cline{2-2}
				& $\begin{array}{l} \alpha_1 \neq 1,\ \beta_2= \alpha_1,\ \beta_1\neq 0\end{array}$& & \\
				\hline
				$A_{3}$  & $1-2\beta_2-4\beta_1\neq0 $& $1-2\beta_2-4\beta_1=0 $ &$-$
				\\
				\hline
				$A_{4}$  &$-$ &  $\alpha_1\neq 1$  & $\alpha_1=1,\ \beta_2\neq \alpha_1$
				\\
				\hline
				$ A_{5}$  &$-$ & $+$ &$-$
				\\
				\hline
				\multirow{2}{*}{$ A_{6}$ }& $\alpha_1=0,\ \beta_1\neq 0,$ &  $\alpha_1= \frac{1}{2},\ \beta_1=0,$ & \multirow{2}{*}{$\alpha_1 = \frac{1}{2},\ \beta_1=0$}
				\\
				\cline{2-3}
				& $\alpha_1\neq 0, \frac{1}{2},\ \beta_1^2\neq \alpha_1$& $\alpha_1\neq 0, \frac{1}{2},\ \beta_1^2= \alpha_1$& \\
				\hline
				$ A_{7}$  & $\beta_1\neq 0,$ &  $\beta_1=0,$ & $-$
				\\
				\hline
				$ A_{8}$  &$-$ & $\alpha_1\neq 0$  &$\alpha_1=0$
				\\
				\hline
				$A_{9}$  & $-$ & $+$ & $-$
				\\
				\hline
				$A_{10}$  & $-$ &  $+$  &$-$
				\\
				\hline
				$A_{11}$ & $+$ & $-$ &  $-$
				\\
				\hline
				$A_{12}$  &  $-$ & $+$ &$-$
				\\
				\hline
			\end{longtable}
\end{footnotesize}
\begin{footnotesize}
		\begin{table}[h]\label{T2}
			\caption{Nontrivial left ideals, $Char(\mathbb{F})= 2,3.$}
			\label{T2}
			\noindent\makebox[\linewidth]{\begin{tabular}{ |c|c| c |c| }
				\hline
				\textbf{Alg.} & $\begin{array}{cc}&\textbf{No}\\ &\textbf{left ideal}\end{array}$ & $\begin{array}{cc}&\textbf{One}\\ &\textbf{left ideal}\end{array}$ &  $\begin{array}{cc}&\textbf{Two}\\ &\textbf{left ideals}\end{array}$  \\
				\hline
				\multirow{4}{*}{$A_{1,2}$}&$\begin{array}{l}\alpha_4=0,\ \alpha_2 \neq0,\\ \alpha_2(1-\alpha_1)^2-\beta_1\neq 0\end{array}$& $\begin{array}{l}\alpha_4=0,\ \alpha_2\neq 0, \\ \alpha_2(1-\alpha_1)^2-\beta_1=0\end{array}$ & $\begin{array}{l}\alpha_4\neq 0,\\ \alpha_2= \beta_1=0\end{array}$
				\\
				\cline{2-4}
				& $\begin{array}{l}\alpha_2\alpha_4\neq 0\\
				\alpha_4^2\beta_1^2+\alpha_2^2(1-\alpha_1)^2+\alpha_2\beta_1\neq 0\end{array}$ & $\begin{array}{l}\alpha_2\alpha_4\neq 0\\
				\alpha_4^2\beta_1^2+\alpha_2^2(1-\alpha_1)^2+\alpha_2\beta_1=0\end{array}$& \multirow{2}{*}{$\begin{array}{l}\alpha_4= \alpha_2=0,\\ \beta_1=0\end{array}$}  \\
				\cline{2-3}
				& $ \alpha_4\beta_1\neq 0,\ \alpha_2=0$&$\alpha_4=0,\ \alpha_2=0,\ \beta_1\neq 0$&
				\\
				\hline
				$A_{2,2}$& $(\alpha_1^2+\beta_2^2)(1-\alpha_1)+\beta_1^2\neq 0$&	$(\alpha_1^2+\beta_2^2)(1-\alpha_1)+\beta_1^2=0$& $-$
				\\
				\hline
				$A_{3,2}$  & $\alpha_1 \neq 1 $& $\alpha_1=1, \ \beta_2=1$ &$\alpha_1=1, \ \beta_2\neq 1$
				\\
				\hline
				$A_{4,2}$  &$-$ & $\alpha_1\neq 1$ &  $\alpha_1=1,\ \alpha_1-\beta_2 \neq 0$
				\\
				\hline
				$ A_{5,2}$  &$- $ & $+$ &$-$
				\\
				\hline
				$ A_{6,2}$ & $\alpha_1+\beta_1^2\neq 0,$ &  $\alpha_1+\beta_1^2=0,$ &$-$
				\\
				\hline
				$ A_{7,2}$  & $\alpha_1\neq 0,$ &  $-$ & $\alpha_1=0$
				\\
				\hline
				$ A_{8,2}$  & $-$ &$\alpha_1\neq 0$  &$\alpha_1=0$
				\\
				\hline
				$A_{9,2}$  & $-$ & $+$ & $-$
				\\
				\hline
				$A_{10,2}$  & $-$ &  $+$  &$-$
				\\
				\hline
				$A_{11,2}$ & $+$ & $-$ &  $-$
				\\
				\hline
				$A_{12,2}$  &  $-$ & $+$ &$-$
				\\
				\hline
				\multirow{16}{*}{$A_{1,3}$}&$\begin{array}{l}\alpha_4=0,\ \alpha_2\neq 0,1\\
				\beta_1\neq\frac{(1-\alpha_1)(2\alpha_1+\alpha_2)}{(2\alpha_2+1)^2}\end{array}$& $\begin{array}{l}\alpha_4=0,\ \alpha_2\neq 0,1\\
				\beta_1=\frac{(1-\alpha_1)(2\alpha_1+\alpha_2)}{(2\alpha_2+1)^2}\end{array}$ &\multirow{2}{*}{ $\begin{array}{l}\alpha_4\neq 0,\\ \alpha_1= 1+\frac{(2\alpha_2+1)^2}{\alpha_4},\\ \beta_1=\frac{\alpha_2(2\alpha_2+1)}{2\alpha_4},\\
					\alpha_2=-(\alpha_4+1)\end{array}$}
				\\
				\cline{2-3}
				& $\begin{array}{l}\alpha_4= 0,\
				\alpha_1\neq 1,\
				\alpha_2=1 \end{array}$ &$\begin{array}{l}\alpha_4\neq 0,\
				\alpha_1=1+ \frac{(2\alpha_2+1)^2}{\alpha_4},\\ \beta_1=\frac{(2\alpha_2+1)(2\alpha_2^2+\alpha_2-\alpha_4\alpha_1)}{\alpha_4^2},	
				\end{array}$	&    \\
				\cline{2-4}
				& $ \begin{array}{l} P\neq 0,\ \alpha_4\neq 0,\\
				\alpha_1\neq 1+\frac{(2\alpha_2+1)^2}{\alpha_4}, \frac{\alpha_2(2\alpha_2+1)}{2\alpha_4}
				\end{array}$& $\begin{array}{l} P=0,\ \alpha_4\neq 0,\\
				\alpha_1\neq 1+ \frac{(2\alpha_2+1)^2}{\alpha_4}, \frac{\alpha_2(2\alpha_2+1)}{2\alpha_4}
				\end{array}$ & \multirow{2}{*}{$\begin{array}{l}\alpha_4= 0,\\ \alpha_1=1,\\
					\alpha_2=1,\
					\beta_1\neq 2
					\end{array}$}  \\
				\cline{2-3}
				&$\begin{array}{l}\alpha_4\neq 0,\
				\alpha_1=1+ \frac{(2\alpha_2+1)^2}{\alpha_4},\\ \beta_1\neq\frac{(2\alpha_2+1)(2\alpha_2^2+\alpha_2-\alpha_4\alpha_1)}{\alpha_4^2},	
				\end{array}$ &  $
				\begin{array}{l}\alpha_1= 1,\
				\alpha_2=1,\\ \alpha_4=0,\
				\beta_1=2
				\end{array} $ & \\
				\cline{2-4}
				&\small{$\begin{array}{l}\alpha_4\neq 0,\
				\alpha_1\neq 1+\frac{(2\alpha_2+1)^2}{\alpha_4},\\
				\alpha_1=\frac{\alpha_2(2\alpha_2+1)}{2\alpha_4},\
				\beta_1\neq \frac{\alpha_2(1-\alpha_1)}{\alpha_4},
				\end{array}$} & $\begin{array}{l}\alpha_4= \alpha_2= 0\\
				\beta_1\neq 2\alpha_1(1-\alpha_1)\end{array}$ &$\begin{array}{l}\alpha_4=\alpha_2= 0\\
				\beta_1=2\alpha_1(1-\alpha_1)\end{array}$
				\\
				\hline
				\multirow{6}{*}{$A_{2,3}$}& $\beta_1\neq 0,\ \alpha_1=1$&	$\alpha_1=1,\ \beta_1=0$&
				\\
				\cline{2-3}
				&$\begin{array}{l} \alpha_1 \neq 1,\ \beta_2\neq \alpha_1,\\ \beta_1^2 +\alpha_1-1\neq 0\end{array}$& \multirow{2}{*}{$\begin{array}{l} \alpha_1 \neq 1,\ \beta_2\neq \alpha_1,\\ \beta_1^2 +\alpha_1-1= 0\end{array}$ } & $\begin{array}{l} \beta_2=\alpha_1 \neq 1,\ \beta_1=0\end{array}$\\
				\cline{2-2}
				& $\begin{array}{l} \beta_2=\alpha_1 \neq 1,\ \beta_1\neq 0\end{array}$& &  \\
				\hline
				$A_{3,3}$  & $1-2\beta_2-\beta_1\neq0 $& $1-2\beta_2-\beta_1=0 $ &$-$
				\\
				\hline
				$A_{4,3}$  &$-$ &  $\alpha_1\neq 1$  & $\alpha_1=1,\ \beta_2\neq \alpha_1$
				\\
				\hline
				$ A_{5,3}$  &$-$ & $+$ &$-$
				\\
				\hline
				\multirow{2}{*}{$ A_{6,3}$ }& $\alpha_1=0,\ \beta_1\neq 0,$ &  $\alpha_1= -1,\ \beta_1=0,$ & \multirow{2}{*}{$\alpha_1 = -1,\ \beta_1=0$}
				\\
				\cline{2-3}
				& $\alpha_1\neq 0, -1,\ \beta_1^2\neq \alpha_1$& $\alpha_1\neq 0, -1,\ \beta_1^2= \alpha_1$&  \\
				\hline
				$ A_{7,3}$  & $\beta_1\neq 0,$ &  $\beta_1=0,$ & $-$
				\\
				\hline
				$ A_{8,3}$  &$-$ & $\alpha_1\neq 0$  &$\alpha_1=0$
				\\
				\hline
				$A_{9,3}$  & $+$ & $-$ & $-$
				\\
				\hline
				$A_{10,3}$  & $-$ &  $+$  &$-$
				\\
				\hline
				$A_{11,3}$ & $-$ & $+$ &  $-$
				\\
				\hline
				$A_{12,3}$  &  $-$ & $+$ &$-$
				\\
				\hline
			\end{tabular}}\end{table}
\end{footnotesize}

	\begin{center}\begin{small}
		\begin{table}[h] \label{T3}
			\caption{Nontrivial right ideals, $Char(\mathbb{F})\neq 2,3.$\newline \scriptsize $P=(\alpha _4 \beta _1+\alpha _1 (\alpha _2+1))^2+ \alpha _2(\alpha _2+1)(-4 \alpha _2 \beta _1-4 \alpha _1^2+2 \alpha _1)+\alpha _4 (4 \alpha _2 \alpha _1 \beta _1-2 \alpha _2 \beta _1+4 (\alpha _1-1) \alpha _1^2+\alpha _1).$}
			\label{T3}
			\noindent\makebox[\linewidth]{\begin{tabular}{ |c|c| c |c| }
				\hline
				\textbf{Alg.} & $\begin{array}{cc}&\textbf{No}\\ &\textbf{right ideal}\end{array}$  &  $\begin{array}{cc}&\textbf{One}\\ &\textbf{right ideal}\end{array}$  &  $\begin{array}{cc}&\textbf{Two}\\ &\textbf{right ideals}\end{array}$   \\
				\hline
				\multirow{16}{*}{$A_{1}$}&$\begin{array}{l}\alpha_4=0,\\ \alpha_2\neq 0,-1\\
				\beta_1\neq\frac{\alpha_1(2\alpha_2+\alpha_1-3\alpha_1\alpha_2)}{4\alpha_2^2}\end{array}$& $\begin{array}{l}\alpha_4=0,\\ \alpha_2\neq 0,-1\\
				\beta_1=\frac{\alpha_1(2\alpha_2+\alpha_1-3\alpha_1\alpha_2)}{4\alpha_2^2}\end{array}$ & \multirow{5}{*}{$\begin{array}{l}\alpha_4\neq 0,\\  \alpha_1=\frac{1}{2}+ \frac{\alpha_2(\alpha_2+1)}{\alpha_4},\\ \alpha_2\neq \frac{-\alpha_4}{2}\\
					\beta_1= -\frac{\alpha_1(\alpha_2+1)}{\alpha_4}
					\end{array}$}
				\\
				\cline{2-3}
				& $\begin{array}{l}\alpha_4=\alpha_2= 0,\\
				\alpha_1\neq 0\\
				\end{array}$ &$
				\begin{array}{l}\alpha_4= \alpha_1=\alpha_2=0,\\
				\beta_1=-\frac{1}{4}
				\end{array} $&  \multirow{12}{*}{$\begin{array}{l}\alpha_4=\alpha_1=\alpha_2=0,\\
					\beta_1\neq -\frac{1}{4}
					\end{array}$}  \\
				\cline{2-3}
				& $ \begin{array}{l}\alpha_4\neq 0,\\
				\alpha_1\neq \frac{\alpha_2^2}{\alpha_4},\ \frac{1}{2}+ \frac{\alpha_2(\alpha_2+1)}{\alpha_4},\\
				P\neq 0
				\end{array}$& $\begin{array}{l}\alpha_4\neq 0,\\
				\alpha_1\neq \frac{\alpha_2^2}{\alpha_4}, \frac{1}{2}+\frac{\alpha_2(\alpha_2+1)}{\alpha_4},\\
				P=0
				\end{array}$ & \\
				\cline{2-4}
				&$\begin{array}{l}\alpha_4\neq 0,\\
				\alpha_1=\frac{\alpha_2^2}{\alpha_4},\\ \beta_1\neq\frac{\alpha_2^2(\alpha_2+2\alpha_4-4\alpha_4\alpha_1+1)}{\alpha_4^2},	
				\end{array}$ &$\begin{array}{l}\alpha_4\neq 0,\\
				\alpha_1=\frac{\alpha_2^2}{\alpha_4},\\ \beta_1=\frac{\alpha_2^2(\alpha_2+2\alpha_4-4\alpha_4\alpha_1+1)}{\alpha_4^2}	
				\end{array}$ & \\
				\cline{2-4}
				&$\begin{array}{l}\alpha_4\neq 0,\\
				\alpha_1=\frac{1}{2}+\frac{\alpha_2(\alpha_2+1)}{\alpha_4},\\
				 \alpha_2\neq \frac{-\alpha_4}{2},\\
				\beta_1\neq -\frac{\alpha_1(\alpha_2+1)}{\alpha_4}
				\end{array}$ &  $\begin{array}{l}\alpha_4=0,\\ \alpha_2=-1,\\
				\beta_1\neq\frac{\alpha_1(2\alpha_1-1)}{2}\end{array}$ & $\begin{array}{l}\alpha_4=0,\\ \alpha_2=-1,\\
				\beta_1=\frac{\alpha_1(2\alpha_1-1)}{2}\end{array}$
				\\
				\hline
				\multirow{6}{*}{$A_{2}$}& $\begin{array}{l} \beta_1\neq 0,\\ \beta_2=0\end{array}$&	$\begin{array}{l} \beta_2\neq 0,\\ \beta_1=0\end{array}$&
				\\
				\cline{2-3}
				&$\begin{array}{l} \alpha_1 \neq \frac{1}{2},\\ \beta_2\neq 0,\\ \beta_1^2 -\beta_2(2\alpha_1-1)^2\neq 0\end{array}$& \multirow{4}{*}{$\begin{array}{l} \alpha_1 \neq \frac{1}{2},\\ \beta_2\neq 0,\\ \beta_1^2 -\beta_2(2\alpha_1-1)^2= 0\end{array}$ } & $\begin{array}{l} \alpha_1 =\frac{1}{2},\\ \beta_2\neq0,\\ \beta_1=0\end{array}$\\
				\cline{2-2}
				& $\begin{array}{l} \alpha_1 =\frac{1}{2},\\ \beta_2 \neq 0,\\ \beta_1\neq 0\end{array}$& &  \\
				\hline
				$A_{3}$  & $\beta_2^2-2\beta_2-4\beta_1\neq0 $& $\beta_2^2-2\beta_2-4\beta_1=0 $ &$-$
				\\
				\hline
				$A_{4}$ &$-$ &$\beta_2 \neq 0$ &  $\beta_2=0,\ \alpha_1\neq \frac{1}{2}$
				\\
				\hline
				$ A_{5}$  &$- $ & $+$ &$-$
				\\
				\hline
				\multirow{2}{*}{$ A_{6}$ }& $\alpha_1=1,\ \beta_1\neq 0,$ &  $\alpha_1= 1,\ \beta_1=0,$ & \multirow{4}{*}{$\alpha_1 = \beta_1=0$}
				\\
				\cline{2-3}
				& $\begin{array}{l} \alpha_1\neq 1, 0,\\ \beta_1^2 - 4\alpha_1^2+4\alpha_1^3\neq 0\end{array}$& \multirow{2}{*}{$\begin{array}{l} \alpha_1\neq 1, 0,\\ \beta_1^2 - 4\alpha_1^2+4\alpha_1^3= 0\end{array}$}& \\
				\cline{2-2}
				& $\alpha_1=0,\ \beta_1\neq 0$ & &
				\\
				\hline
				$ A_{7}$  & $\beta_1\neq \frac{1}{4}$ &  $\beta_1=\frac{1}{4}$ & $-$
				\\
				\hline
				$ A_{8}$  &$-$ & $\alpha_1\neq 1$  &$\alpha_1=1$
				\\
				\hline
				$A_{9}$  & $-$ & $+$ & $-$
				\\
				\hline
				$A_{10}$  & $-$ &  $+$  &$-$
				\\
				\hline
				$A_{11}$ & $+$ & $-$ &  $-$
				\\
				\hline
				$A_{12}$  &  $-$ & $+$ &$-$
				\\
				\hline
			\end{tabular}}
			\end{table}
\end{small}\end{center}
\begin{tiny}
		\begin{table}[h] \label{T4}
			\caption{Nontrivial right ideals,\ \ $Char(\mathbb{F})= 2,3.$}
						\label{T4}
			\noindent\makebox[\linewidth]{\begin{tabular}{ |c|c| c |c| c| }
				\hline
				\textbf{Alg.} & $\begin{array}{cc}&\textbf{No}\\ &\textbf{right ideal}\end{array}$  &  $\begin{array}{cc}&\textbf{One}\\ &\textbf{right ideal}\end{array}$  &  $\begin{array}{cc}&\textbf{Two}\\ &\textbf{right ideals}\end{array}$ \\
				\hline
				\multirow{3}{*}{$A_{1,2}$}&$\begin{array}{l}\alpha_4=0,\ \alpha_1\neq 0,\ \alpha_2\neq 1\end{array}$& $\begin{array}{l}\alpha_4\neq 0,\\ \alpha_1=-\frac{(\alpha_2+1)^2(1-\alpha_1)^2-\beta_1^2\alpha_4^2}{\alpha_4}\end{array}$ & $\begin{array}{l}\alpha_4=\alpha_1=0,\\ \alpha_2\neq 1\end{array}$
				\\
				\cline{2-4}
				& $\begin{array}{l}\alpha_4\neq 0\\
				\alpha_1\neq -\frac{(\alpha_2+1)^2(1-\alpha_1)^2-\beta_1^2\alpha_4^2}{\alpha_4}\end{array}$& $\begin{array}{l}\alpha_4=0,\ \alpha_1\neq 0,\ \alpha_2=1\end{array}$&$\begin{array}{l}\alpha_4=\alpha_1= 0,\\ \alpha_2=1\end{array}$
				\\
				\hline
				$A_{2,2}$& $\beta_1-\beta_2^2\neq 0$&	$\beta_1-\beta_2^2=0$& $-$
				\\
				\hline
				$A_{3,2}$  & $\beta_2 \neq 0 $& $-$ &$\beta_2=0$
				\\
				\hline
				$A_{4,2}$  &$-$& $\beta_2\neq 0$ &  $\beta_2=0$
				\\
				\hline
				$ A_{5,2}$  &$-$ & $+$ &$-$
				\\
				\hline
				$ A_{6,2}$ & $\beta_1\neq 0,$ &  $\beta_1=0,$ &$-$
				\\
				\hline
				$ A_{7,2}$  & $\alpha_1\neq 0,$ &  $\alpha_1=0$ & $-$
				\\
				\hline
				$ A_{8,2}$  & $-$& $\alpha_1\neq 1$  &$-$
				\\
				\hline
				$A_{9,2}$  & $-$ & $+$ & $-$
				\\
				\hline
				$A_{10,2}$  & $-$ &  $+$  &$-$
				\\
				\hline
				$A_{11,2}$ & $+$ & $-$ &  $-$
				\\
				\hline
				$A_{12,2}$  &  $-$ & $+$ &$-$
				\\
				\hline
				\multirow{16}{*}{$A_{1,3}$}&$\begin{array}{l}\alpha_4=0,\ \alpha_2\neq 0,-1\\
				\beta_1\neq\frac{\alpha_1(2\alpha_2+\alpha_1)}{\alpha_2^2}\end{array}$& $\begin{array}{l}\alpha_4=0,\ \alpha_2\neq 0,-1\\
				\beta_1=\frac{\alpha_1(2\alpha_2+\alpha_1)}{\alpha_2^2}\end{array}$ & \multirow{5}{*}{$\begin{array}{l}\alpha_4\neq 0,\\ \alpha_1=\frac{1}{2}+ \frac{\alpha_2(\alpha_2+1)}{\alpha_4},\\ \alpha_2\neq \alpha_4\\
					\beta_1= -\frac{\alpha_1(\alpha_2+1)}{\alpha_4}
					\end{array}$}
				\\
				\cline{2-3}
				& $\alpha_1\neq 0,\ \alpha_4=\alpha_2= 0$ &$
				\alpha_4= \alpha_1=\alpha_2=0,\				\beta_1=-\frac{1}{4}
				 $&  \multirow{10}{*}{$\begin{array}{l}\alpha_4=\alpha_1=\alpha_2=0,\\
					\beta_1\neq -1
					\end{array}$}  \\
				\cline{2-3}
				& $ \begin{array}{l} P\neq 0,\ \alpha_4\neq 0,\\
				\alpha_1\neq \frac{\alpha_2^2}{\alpha_4},\ -1+ \frac{\alpha_2(\alpha_2+1)}{\alpha_4}
				\end{array}$& $\begin{array}{l}P=0,\ \alpha_4\neq 0,\\
				\alpha_1\neq \frac{\alpha_2^2}{\alpha_4}, -1+\frac{\alpha_2(\alpha_2+1)}{\alpha_4}
				\end{array}$ &  \\
				\cline{2-4}
				&$\begin{array}{l}\alpha_4\neq 0,\
				\alpha_1=\frac{\alpha_2^2}{\alpha_4},\\ \beta_1\neq\frac{\alpha_2^2(\alpha_2+2\alpha_4-\alpha_4\alpha_1+1)}{\alpha_4^2},	
				\end{array}$ &$\begin{array}{l}\alpha_4\neq 0,\
				\alpha_1=\frac{\alpha_2^2}{\alpha_4},\\ \beta_1=\frac{\alpha_2^2(\alpha_2+2\alpha_4-\alpha_4\alpha_1+1)}{\alpha_4^2}	
				\end{array}$ & \\
				\cline{2-4}
				&$\begin{array}{l}\alpha_2\neq \alpha_4,\ \alpha_4\neq 0,\\
				\alpha_1=-1+\frac{\alpha_2(\alpha_2+1)}{\alpha_4},\\
				 				\beta_1\neq -\frac{\alpha_1(\alpha_2+1)}{\alpha_4}
				\end{array}$ &  $\begin{array}{l}\alpha_4=0,\\ \alpha_2=-1,\\
				\beta_1\neq -\alpha_1(2\alpha_1-1)\end{array}$ & $\begin{array}{l}\alpha_4=0,\\ \alpha_2=-1,\\
				\beta_1=-\alpha_1(2\alpha_1-1)\end{array}$
				\\
				\hline
				\multirow{6}{*}{$A_{2,3}$}& $\begin{array}{l} \beta_1\neq 0,\ \beta_2=0\end{array}$&	$\begin{array}{l} \beta_2\neq 0,\ \beta_1=0\end{array}$&
				\\
				\cline{2-3}
				&$\begin{array}{l} \alpha_1 \neq -1,\ \beta_2\neq 0,\\ \beta_1^2 -\beta_2(2\alpha_2-1)^2\neq 0\end{array}$& $\begin{array}{l} \alpha_1 \neq -1,\\ \beta_2\neq 0,\\ \beta_1^2 -\beta_2(2\alpha_2-1)^2= 0\end{array}$  &$\begin{array}{l} \alpha_1 =-1,\\ \beta_2\neq0,\\ \beta_1=0\end{array}$\\
				\cline{2-2}
				& $\begin{array}{l} \alpha_1 =-1,\ \beta_2 \neq 0,\ \beta_1\neq 0\end{array}$& &  \\
				\hline
				$A_{3,3}$  & $\beta_2^2-2\beta_2-\beta_1\neq0 $& $\beta_2^2-2\beta_2-\beta_1=0 $ &$-$
				\\
				\hline
				$A_{4,3}$ &$-$ &$\beta_2 \neq 0$ &  $\beta_2=0,\ \alpha_1\neq -1$
				\\
				\hline
				$ A_{5,3}$  &$- $ & $+$ &$-$
				\\
				\hline
				\multirow{2}{*}{$ A_{6,3}$ }& $\alpha_1=1,\ \beta_1\neq 0,$ &  $\alpha_1= 1,\ \beta_1=0,$ & \multirow{4}{*}{$\alpha_1 = \beta_1=0$}
				\\
				\cline{2-3}
				& $\begin{array}{l} \alpha_1\neq 0,1,\ \beta_1^2 - \alpha_1^2+\alpha_1^3\neq 0\end{array}$& \multirow{2}{*}{$\begin{array}{l} \alpha_1\neq 1, 0,\\ \beta_1^2 - \alpha_1^2+\alpha_1^3= 0\end{array}$}&  \\
				\cline{2-2}
				& $\alpha_1=0,\ \beta_1\neq 0$ & &
				\\
				\hline
				$ A_{7,3}$  & $\beta_1\neq 1$ &  $\beta_1=1$ & $-$
				\\
				\hline
				$ A_{8,3}$  &$-$ & $\alpha_1\neq 1$  &$\alpha_1=1$
				\\
				\hline
				$A_{9,3}$  & $-$ & $+$ & $-$
				\\
				\hline
				$A_{10,3}$  & $-$ &  $+$  &$-$
				\\
				\hline
				$A_{11,3}$ & $-$ & $+$ &  $-$
				\\
				\hline
				$A_{12,3}$  &  $-$ & $+$ &$-$
				\\
				\hline
			\end{tabular}}\end{table}\end{tiny}

\begin{center}\begin{small}
 		\begin{table}[h] \label{T5}
 			\caption{Two-sided ideals, $Char(\mathbb{F})\neq 2,3$.}
 			\label{T5}
 			\begin{tabular}{ |c|c| c |c|  }
 				\hline
 				\textbf{Algebra} & \textbf{No Ideal} &  \textbf{One Ideal} &  \textbf{Two Ideals} \\
 				\hline
 				$A_{1}$&$\begin{array}{l}\alpha_4\neq -\alpha_1-2\alpha_2\  \mbox{or}\\
 				\beta_1\neq 2\alpha_1+\alpha_2\end{array}$ & $\begin{array}{l}\alpha_4=-\alpha_1-2\alpha_2\\
 				\beta_1= 2\alpha_1+\alpha_2
 				\end{array}$&$-$
 				\\
 				\hline
 				\multirow{6}{*}{$A_{2}$}& $\begin{array}{l}\alpha_1=1\\ \beta_1\neq 0,\\ \beta_2= 0\end{array}$&	$\begin{array}{l} \alpha_1=1,\\ \beta_1=\beta_2=0\end{array}$&
 				\\
 				\cline{2-3}
 				&$\begin{array}{l} \alpha_1 \neq \frac{1}{2}, 1,\\ \beta_1^2 +\alpha_1-1\neq 0\\ \beta_2=1-\alpha_1\end{array}$& \multirow{2}{*}{$\begin{array}{l} \alpha_1 \neq \frac{1}{2}, 1,\\ \beta_2=1-\alpha_1,\\ \beta_1^2 +\alpha_1-1= 0\end{array}$ } & $\begin{array}{l} \alpha_1=\beta_2= \frac{1}{2},\\ \beta_1=0\end{array}$\\
 				\cline{2-2}
 				& $\begin{array}{l} \alpha_1=\beta_2= \frac{1}{2},\\ \beta_1\neq 0\end{array}$& & \\
 				\cline{2-2}
 				& $ \beta_2\neq 1-\alpha_1$& &  \\
 				\hline
 				\multirow{2}{*}{$A_{3}$}  & $\beta_1\neq-\frac{1}{4},\ \beta_2=1$& \multirow{2}{*}{$\beta_1=-\frac{1}{4},\ \beta_2=1$}& \multirow{2}{*}{$-$}\\
 				\cline{2-2}
 				& $ \beta_1\neq 1$& & 
 				\\
 				\hline
 				\multirow{2}{*}{$A_{4}$}  &\multirow{2}{*}{$-$} &  $\alpha_1\neq 1,\ \beta_2=1-\alpha_1$&\multirow{2}{*}{$\alpha_1=1,\ \beta_2= 0$} \\
 				\cline{3-3}
 				& &$ \beta_2\neq 1-\alpha_1$
 				&
 				\\
 				\hline
 				$ A_{5}$  &$-$ & $+$ &$-$
 				\\
 				\hline
 				$ A_{6}$ & $+$ & $-$ & $-$
 				\\
 				\hline
 				$ A_{7}$  & $+$ & $-$& $-$
 				\\
 				\hline
 				$ A_{8}$  &$-$ & $+$  &$-$
 				\\
 				\hline
 				$A_{9}$  & $-$ & $+$ & $-$
 				\\
 				\hline
 				$A_{10}$  & $-$ &  $+$  &$-$
 				\\
 				\hline
 				$A_{11}$ & $+$ & $-$ &  $-$   				\\
 				\hline
 				$A_{12}$  &  $-$ & $+$ &$-$
 				\\
 				\hline
 			\end{tabular}\end{table}
 \end{small}\end{center}

 \begin{center}\begin{small}
 		\begin{table}[h]\label{T6}
 			\caption{Two-sided ideals, $Char(\mathbb{F})= 2,3$.}
 			\label{T6}
 			\begin{tabular}{ |c|c|c| c |c|  }
 				\hline
 				&\textbf{Algebra} & \textbf{No Ideal} &  \textbf{One Ideal} &  \textbf{Two Ideals} \\
 				\cline{2-5}
 				\multirow{15}{*}{\begin{turn}{90}$Char\left(\mathbb{F}\right)= 2
$\end{turn}}&$A_{1,2}$&$\alpha_4\neq \alpha_1$ or $\beta_1\neq \alpha_2$ & $\alpha_4=\alpha_1,\ \beta_1=\alpha_2$ & $-$
 				\\
 				\cline{2-5}
 				&$A_{2,2}$& $\begin{array}{l}
 				1-\alpha_1+\beta_1^2\neq 0 \\ \beta_2=1-\alpha_1\end{array}$&$\begin{array}{l}1-\alpha_1+\beta_1^2= 0\\\beta_2=1-\alpha_1 \end{array}$&\multirow{2}{*}{$-$ } \\
 				\cline{3-3}
 				& &$\beta_2\neq 1-\alpha_1$&	&
 				\\
 				\cline{2-5}
 				& $A_{3,2}$  & $\alpha_1 \neq 1 $ or $\beta_2 \neq 0$& $-$ &$\alpha_1 = 1, \ \beta_2 = 0$
 				\\
 				\cline{2-5}
 				& $A_{4,2}$  &\multirow{2}{*}{$-$} & $\alpha_1\neq 1,\ \beta_2=1-\alpha_1$&\multirow{2}{*}{$\alpha_1=1,\ \beta_2 = 0$} \\
 				\cline{4-4}
                & & &$\beta_2\neq 1-\alpha_1$ &
 				\\
 				\cline{2-5}
 				&$ A_{5,2}$  &$- $ & $+$ &$-$
 				\\
 				\cline{2-5}
 				&$ A_{6,2}$ & $+$ & $-$ &$-$
 				\\
 				\cline{2-5}
 				&$ A_{7,2}$  & + &  $-$ &$ -$
 				\\
 				\cline{2-5}
 				&$ A_{8,2}$  & $-$ &$+$&$-$
 				\\
 				\cline{2-5}
 				&$A_{9,2}$  & $-$ & $+$ & $-$
 				\\
 				\cline{2-5}
 				&$A_{10,2}$  & $-$ &  $+$  &$-$
 				\\
 				\cline{2-5}
 				& $A_{11,2}$ & $+$ & $-$ &  $-$
 				\\
 				\cline{2-5}
 				&$A_{12,2}$  &  $-$ & $+$ &$-$
 				\\
 				\hline
 				\multirow{20}{*}{\begin{turn}{90}$Char\left(\mathbb{F}\right)= 3
$\end{turn}}&$A_{1,3}$&$\begin{array}{l}\alpha_4\neq -\alpha_1-2\alpha_2\  \mbox{or}\\
 				\beta_1\neq 2\alpha_1+\alpha_2\end{array}$ & $\begin{array}{l}\alpha_4=-\alpha_1-2\alpha_2\\
 				\beta_1= 2\alpha_1+\alpha_2
 				\end{array}$&$-$
 			\\
 				\cline{2-5}
 				& \multirow{6}{*}{$A_{2,3}$}& $\begin{array}{l}\alpha_1=1\\ \beta_1\neq 0,\\ \beta_2= 0\end{array}$&	$\begin{array}{l} \alpha_1=1,\\ \beta_1=\beta_2=0\end{array}$&
 				\\
 				\cline{3-4}
 				& &$\begin{array}{l} \alpha_1 \neq \frac{1}{2}, 1,\\ \beta_1^2 +\alpha_1-1\neq 0\\ \beta_2=1-\alpha_1\end{array}$& \multirow{4}{*}{$\begin{array}{l} \alpha_1 \neq \frac{1}{2}, 1,\\ \beta_2=1-\alpha_1,\\ \beta_1^2 +\alpha_1-1= 0\end{array}$ } & $\begin{array}{l} \alpha_1=\beta_2= \frac{1}{2},\\ \beta_1=0\end{array}$\\
 				\cline{3-3}
 				& & $\begin{array}{l} \alpha_1=\beta_2= \frac{1}{2},\\ \beta_1\neq 0\end{array}$& &  \\
 				\cline{3-3}
 				& & $ \beta_2\neq 1-\alpha_1$& &  \\
 				\cline{2-5}
 				& \multirow{2}{*}{$A_{3,3}$}  & $\beta_1\neq -1,\ \beta_2=1$ & \multirow{2}{*}{ $ \beta_1=- 1,\ \beta_2=1$}&\multirow{2}{*}{$ -$} \\
 				\cline{3-3}	
 				& &$ \beta_2\neq 1 $
 				&
 				&
 				 				\\
 				\cline{2-5}
 				& \multirow{2}{*}{$A_{4,3}$}  &\multirow{2}{*}{$-$ }&  $\alpha_1\neq 1,\ \beta_2=1-\alpha_1$& \multirow{2}{*}{$\alpha_1=1,\ \beta_2= 0$}\\
 				\cline{4-4}
 				& & &$ \beta_2\neq 1-\alpha_1$
 				&
 				\\
 				\cline{2-5}
 				& $ A_{5,3}$  &$- $ & $+$ &$-$
 				\\
 				\cline{2-5}
 				& $ A_{6,3}$& $+$ &$-$ &$-$
 				\\
 				\cline{2-5}
 				&$ A_{7,3}$  &+ &  $-$ & $-$
 				\\
 				\cline{2-5}
 				&$ A_{8,3}$  &$-$ & $+$ &$-$
 				\\
 				\cline{2-5}
 				&$A_{9,3}$  & $+$ & $-$ & $-$
 				\\
 				\cline{2-5}
 				& $A_{10,3}$  & $-$ &  $+$  &$-$
 				\\
 				\cline{2-5}
 				& $A_{11,3}$ & $-$ & $+$ &  $-$
 				\\
 				\cline{2-5}
 				& $A_{12,3}$  &  $-$ & $+$ &$-$
 				\\
 				\hline
 			\end{tabular}\end{table}
 \end{small}\end{center}

	\begin{center}\begin{small}
		\begin{table}[h]\label{T7}
			\caption{Left quasiunits,\ $Char(\mathbb{F})\neq 2,3$,\ $Char(\mathbb{F})= 2$ and $Char(\mathbb{F})= 3$.}
			\label{T7}
			\begin{tabular}{|c|c|c|}
				\hline
				&\textbf{Algebras} & \textbf{Left quasiunits} \\
						\cline{2-3}
				\multirow{11}{*}{\begin{turn}{90}$Char\left(\mathbb{F}\right)\neq 2,3
$\end{turn}}& $A_1\left(\alpha_1,-\frac{2\alpha_1(\alpha_1-1)+\beta_1}{2\beta_1},  -\frac{2\alpha_2(\alpha_1-1)}{2\beta_1}, \beta_1\right)$,  $\beta_1 \neq 0$
				&  $\frac{2(\alpha_1-1)}{\beta_1}e_1+2e_2$\\
				\cline{2-3}
				&$A_1\left(1,\alpha_2, \frac{\alpha_2(1+2\alpha_2)}{2},0\right)$&$ -(1+2\alpha_2)e_1+2e_2$
				\\
				\cline{2-3}
				&$A_2(\alpha_1,0,\alpha_1)$,\ $\alpha_1\neq 0$& $\frac{1}{\alpha_1}e_1$\\
				\cline{2-3}
				&$A_{4}(1,\beta_2)$& $e_1+te_2$,\  $t\in \mathbb{F}$\\
				\cline{2-3}
				&$A_{4}(\alpha_1,\beta_2)$,\ $\alpha_1\neq 0,1,\beta_2\neq 2\alpha_1-1$&
				$\frac{1}{\alpha_1}e_1$\\
				\cline{2-3}
				&$A_{4}(\alpha_1,2\alpha_1-1)$,\ $\alpha_1\neq 0,1$& $\frac{1}{\alpha_1}e_1+te_2$,\ $t\in \mathbb{F}$\\
				\cline{2-3}
				&$A_5(1)$ & $e_1+te_2$,\ $t\in \mathbb{F}$\\	
				\cline{2-3}
				&$A_6\left(\frac{1}{2}, 0\right)$ & $2e_1$\\
				\cline{2-3}
				&$A_8(\alpha_1)$,\ $\alpha_1 \neq 0, \frac{1}{3}$& $\frac{1}{\alpha_1}e_1$\\
				\cline{2-3}
				&$A_8\left(\frac{1}{3}\right)$ &  $3e_1+te_2$,\ $t\in \mathbb{F}$\\					
				\cline{2-3}
				&$A_{10}$ & $-e_2$\\

							\hline
				\multirow{12}{*}{\begin{turn}{90}$Char\left(\mathbb{F}\right)= 2
$\end{turn}}&$A_{1,2}(\alpha_1,0,\alpha_4,0)$,\ $\alpha_1\neq 0$ & $
				-\frac{1}{\alpha_1}e_1$ \\
				\cline{2-3}
				&$A_{2,2}(\alpha_1,0,\beta_2)$,\ $\alpha_1\neq 0$ & $
				\frac{1}{\alpha_1}e_1$ \\
				\cline{2-3}
				&$A_{3,2}(1,\beta_2)$ & $e_2$ \\
				\cline{2-3}
&$A_{3,2}(\alpha_1,1)$ & $e_2$ \\
				\cline{2-3}
				&$ A_{4,2}(1,\beta_2) $ & $e_1+te_2$,\ $t\in \mathbb{F}$ \\
				\cline{2-3}
				&$ A_{4,2}(\alpha_1,\beta_2),\ \alpha_1\neq 0,1,\ \beta_2\neq 1 $ &
				$\frac{1}{\alpha _1}e_1$ \\
				\cline{2-3}
				&$ A_{4,2}(\alpha_1,1),\ \alpha_1\neq 0,1$ &
				$\frac{1}{\alpha _1}e_1+te_2,\ t\in \mathbb{F}$ \\
				\cline{2-3}				
				&$A_{5,2}(1)$ & $e_1+te_2,\ t\in \mathbb{F}$ \\
				\cline{2-3}
				&$A_{6,2}(\alpha_1,0)$\ $\alpha_1\neq 0$ & $
				\frac{1}{\alpha _1}e_1$ \\
				\cline{2-3}
				&$A_{7,2}(0)$,\ $A_{7,2}(1)$ & $e_2$ \\
				\cline{2-3}
				&$A_{8,2}(\alpha_1),\ \alpha_1\neq 0,1 $& $\frac{1}{\alpha _1}e_1$\\
				\cline{2-3}
				&$A_{8,2}(1)$& $e_1+te_2,\ t\in \mathbb{F}$\\
				\cline{2-3}
				&$A_{10,2}$ & $e_2 $\\
				\hline
				\multirow{11}{*}{\begin{turn}{90}$Char\left(\mathbb{F}\right)=3 $\end{turn}}&$A_{1,3}\left(\alpha_1,\frac{\alpha_1(1-\alpha_1)}{\beta_1}-\frac{1}{2},\frac{\alpha_1(1-\alpha_1)^2}{\beta^2_1}-\frac{1-\alpha_1}{2\beta_1},\beta_1\right)$,\ $\beta_1\neq 0$&  $	\frac{1-\alpha_1}{\beta_1}e_1+2e_2$ \\
				\cline{2-3}
				&$A_{1,3}\left(1,\alpha_2,\frac{\alpha_2(2\alpha_2+1)}{2},0\right)$ &  $(\alpha_2 -1)e_1+2e_2$ \\
				\cline{2-3}
				&$A_{2,3}(\alpha_1,0,\alpha_1)$,\ $\alpha_1\neq 0$ & $
				\frac{1}{\alpha_1}e_1$ \\
				\cline{2-3}
				&$ A_{4,3}(\alpha_1,-1-\alpha_1), \alpha_1\neq 0,1 $ & $
				\frac{1}{\alpha _1}e_1+te_2,\ t\in \mathbb{F}$ \\
				\cline{2-3}
				&$ A_{4,3}(\alpha_1,\beta_2), \alpha_1\neq 0,1,\ \beta_2\neq -1-\alpha_1 $ & $
				\frac{1}{\alpha _1}e_1$ \\
				\cline{2-3}
				&$ A_{4,3}(1,\beta_2)$ & $
				e_1+te_2,\ t\in \mathbb{F}$ \\
				\cline{2-3}
				&$A_{5,3}(1)$ & $e_1+te_2,\ t\in \mathbb{F}$\\
				\cline{2-3}
				&$A_{6,3}(-1,0)$ & $-e_1$ \\
				\cline{2-3}
				&$A_{8,3}(\alpha _1),\ \alpha _1\neq 0$ & $
				\frac{1}{\alpha _1}e_1$ \\
				\cline{2-3}
				&$A_{10,3}$ & $te_1-e_2,\ t\in \mathbb{F}$ \\
				\hline
							\end{tabular}\newline \end{table}
\end{small}\end{center}

%

%
%
%
%

\end{document}